\documentclass{article}
\usepackage [totalwidth=400pt, totalheight=600pt] {geometry}
\usepackage{graphicx}
\usepackage{amsmath}
\usepackage{amsfonts}
\usepackage{amssymb}
\usepackage{caption2}
\usepackage[nooneline, figbotcap]{subfigure}
\usepackage [section]{placeins}

\newtheorem{theorem}{Theorem}

\newtheorem{conjecture}[theorem]{Conjecture}
\newtheorem{corollary}[theorem]{Corollary}

\newtheorem{definition}[theorem]{Definition}
\newtheorem{example}[theorem]{Example}

\newtheorem{lemma}[theorem]{Lemma}

\newtheorem{proposition}[theorem]{Proposition}
\newtheorem{remark}[theorem]{Remark}

\newenvironment{proof}[1][Proof]{\textbf{#1.} }{\ \rule{0.5em}{0.5em}}

\subfiguretopcaptrue
\begin{document}

\title{A set of moves for Johansson representation of 3-manifolds. An
outline.\footnote{This research has been partially supported by a predoctoral
grant from the U.N.E.D. (1999).}}
\author {Rub\'en Vigara
\\Departamento de Matem\'aticas Fundamentales \\U.N.E.D.,
Spain\\rvigara@mat.uned.es}
\date{July 2004}
\maketitle
\begin{abstract}
A Dehn sphere $\Sigma$ \cite{Papa} in a closed 3-manifold $M$ is a
2-sphere immersed in $M$ with only double curve and triple point
singularities. The sphere $\Sigma$ fills $M$ \cite{Montesinos} if
it defines a cell-decomposition of $M$. The inverse image in
$S^{2}$ of the double curves of $\Sigma$ is the Johansson diagram
of $\Sigma$ \cite{Johansson1} and if $\Sigma$ fills $M$ it is
possible to reconstruct $M$ from the diagram. A Johansson
representation of $M$ is the Johansson diagram of a filling Dehn
sphere of $M$. In \cite{Montesinos} it is proved that every closed
3-manifold has a Johansson representation coming from a
nulhomotopic filling Dehn sphere. In this paper a set of moves for
Johansson representations of 3-manifolds is given. In a
forthcoming paper \cite{RHomotopies} it is proved that this set of
moves suffices for relating different Johansson representations of
the same 3-manifold coming from nulhomotopic filling Dehn spheres.
The proof of this result is outlined here.(Math. Subject
Classification: 57N10, 57N35)
\end{abstract}

\section{Introduction.\label{SECTION Introduction}}

Through the whole paper all 3-manifolds are assumed to be closed, that is,
compact connected and without boundary, and all surfaces are assumed to be
compact and without boundary. A surface may have more than one connected
component. We will denote a 3-manifold by $M$ and a surface by $S$.

Let $M$ be a 3-manifold.

\begin{figure}[htb]
\centering \subfigure []{ \label{fig1a}
\includegraphics[ width=0.35\textwidth]{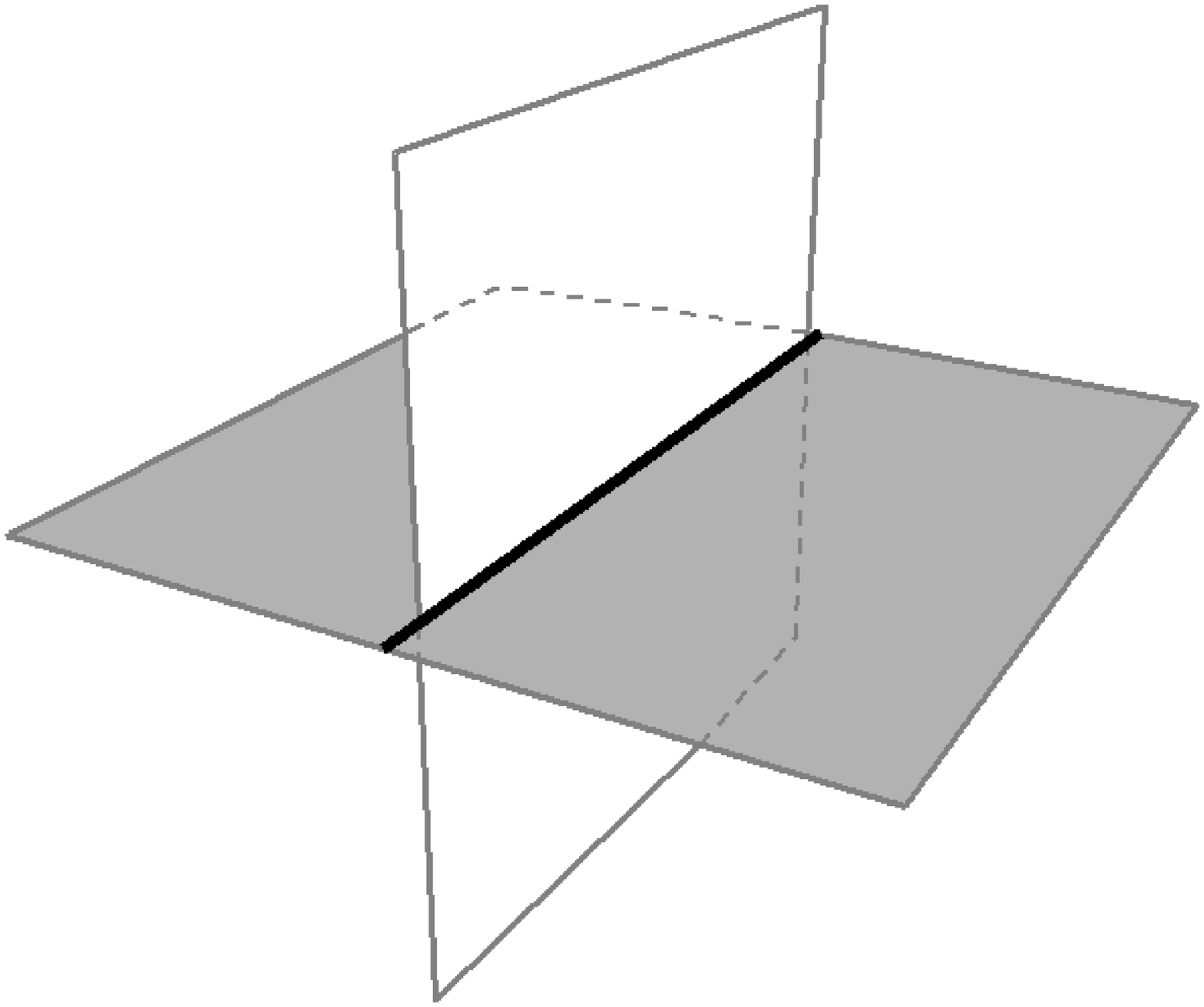}} \hfill
\subfigure[]{ \label{fig1b}
\includegraphics[
width=0.35\textwidth]{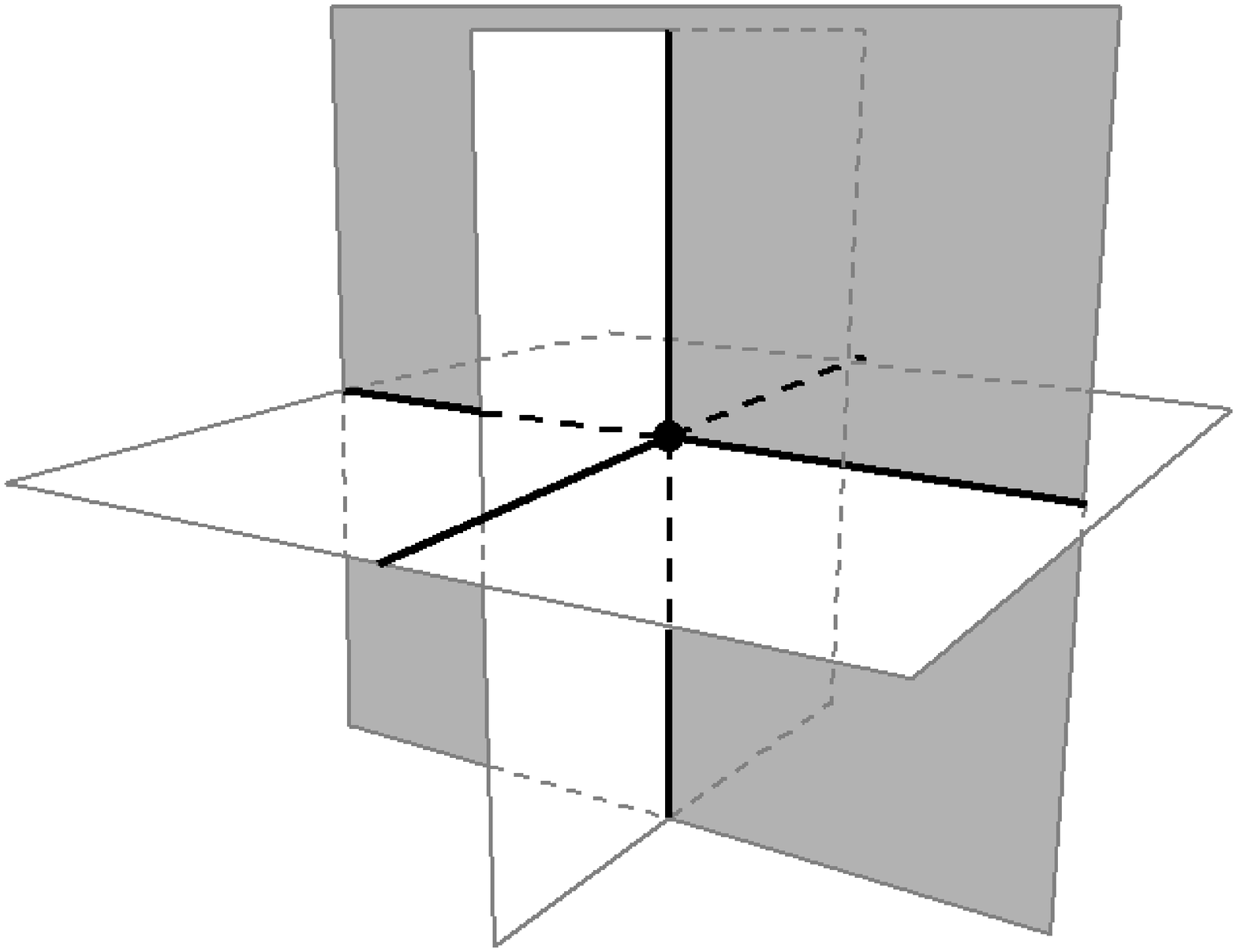}} \caption { }
\end{figure}

A subset $\Sigma\subset M$ is a \textit{Dehn surface} in $M$
\cite{Papa} if there exists a surface $S$ and a transverse
immersion $f:S\rightarrow M$ such that $\Sigma=f\left(  S\right)
$. In this situation we say that $f$ \textit{parametrizes}
$\Sigma$. If $S$ is a 2-sphere then $\Sigma$ is a \textit{Dehn
sphere}. For a Dehn surface $\Sigma\subset M$, its singularities
are divided into \textit{double points }(Figure \ref{fig1a}), and
\textit{triple points} (Figure \ref{fig1b})\textit{, }and they are
arranged along \textit{double curves }(see section \ref{SECTION
Preliminaries} below for definitions). A Dehn surface
$\Sigma\subset M$ \textit{fills} $M$ \cite{Montesinos} if it
defines a cell-decomposition of $M$ in which the 0-skeleton is the
set of triple points of $\Sigma$; the 1-skeleton is the set of
double and triple points of $\Sigma$; and the 2-skeleton is
$\Sigma$ itself. Filling Dehn spheres of 3-manifolds are defined
in \cite{Montesinos} following ideas of W. Haken (see
\cite{Haken1}). In \cite{Montesinos} it is proved the following
Theorem (see also \cite{Anewproof}):

\begin{theorem}
[\cite{Montesinos}]\label{ThmMontesinos}Every closed orientable 3-manifold has
a nulhomotopic filling Dehn sphere.
\end{theorem}

A filling Dehn sphere is \textit{nulhomotopic} if one (and hence any) of its
parametrizations is nulhomotopic, that is, homotopic to a constant map.

Let $\Sigma\subset M$ be a filling Dehn sphere and
$f:S^{2}\rightarrow M$ a transverse immersion parametrizing
$\Sigma$. In this case we say that $f$ is a \textit{filling
immersion}. The inverse image by $f$ in $S^{2}$ of the set of
double and triple points of $\Sigma$ is the \textit{singular set}
of $f$. The singular set of $f$, together with the information of
how its points are identified by $f$ in $M$, is the
\textit{Johansson diagram} of $\Sigma$ in the notation of
\cite{Montesinos}. As it is stated in \cite{Montesinos}, for a
given diagram in $S^{2}$ it is possible to know if it is the
Johansson diagram for a filling Dehn sphere $\Sigma$ in some
3-manifold $M$. If this is the case, it is possible also to
reconstruct such $M$ from the diagram. Thus, Johansson diagrams
are a suitable way for representing all closed, orientable
3-manifolds and it is interesting to further study them. For a
3-manifold $M$, we say that a Johansson diagram of a filling Dehn
sphere of $M$ is a \textit{Johansson representation} of $M$ (see
\cite{Montesinos}). In \cite{Montesinos} an algorithm is given for
obtaining a Johansson representation of a closed orientable
3-manifold $M$ from any Heegaard diagram of $M$. A simpler
algorithm is given in \cite{Anewproof}. In both papers, the
Johansson representations obtained come from nulhomotopic filling
Dehn spheres of $M$.

We will deal here with the problem of deciding how different Johansson
representations of the same 3-manifold are related to each other. With this
problem in mind, we study how different filling Dehn spheres of the same
3-manifold are related to each other. In the forthcoming paper
\cite{RHomotopies}, the following Theorem is proved.

\begin{theorem}
\label{MAINtheorem}Let $M$ be a closed 3-manifold. Let
$f,g:S^{2}\rightarrow M$ be two nulhomotopic filling immersions.
Then, there is a finite sequence of filling immersions
$f=f_{0},f_{1},...,f_{n}=g$ such that for each $i=0,...,n-1$ the
immersions $f_{i}$ and $f_{i+1}$ differ by an ambient isotopy of
$S^{2}$, or by an ambient isotopy of $M$, or by one of the moves
depicted in Figure 2.
\end{theorem}

\begin{figure}[htb]
\centering \subfigure []{ \label{fig2a}
\includegraphics[ width=0.5\textwidth]{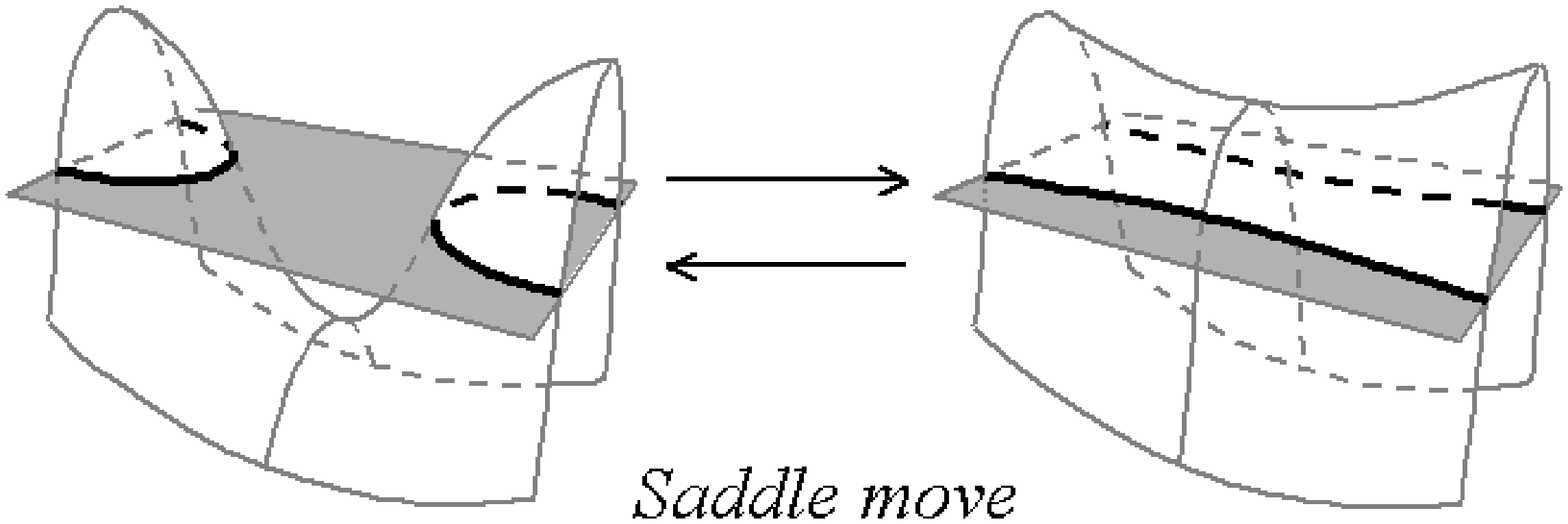}}\\
\subfigure[]{ \label{fig2b}
\includegraphics[
width=0.5\textwidth]{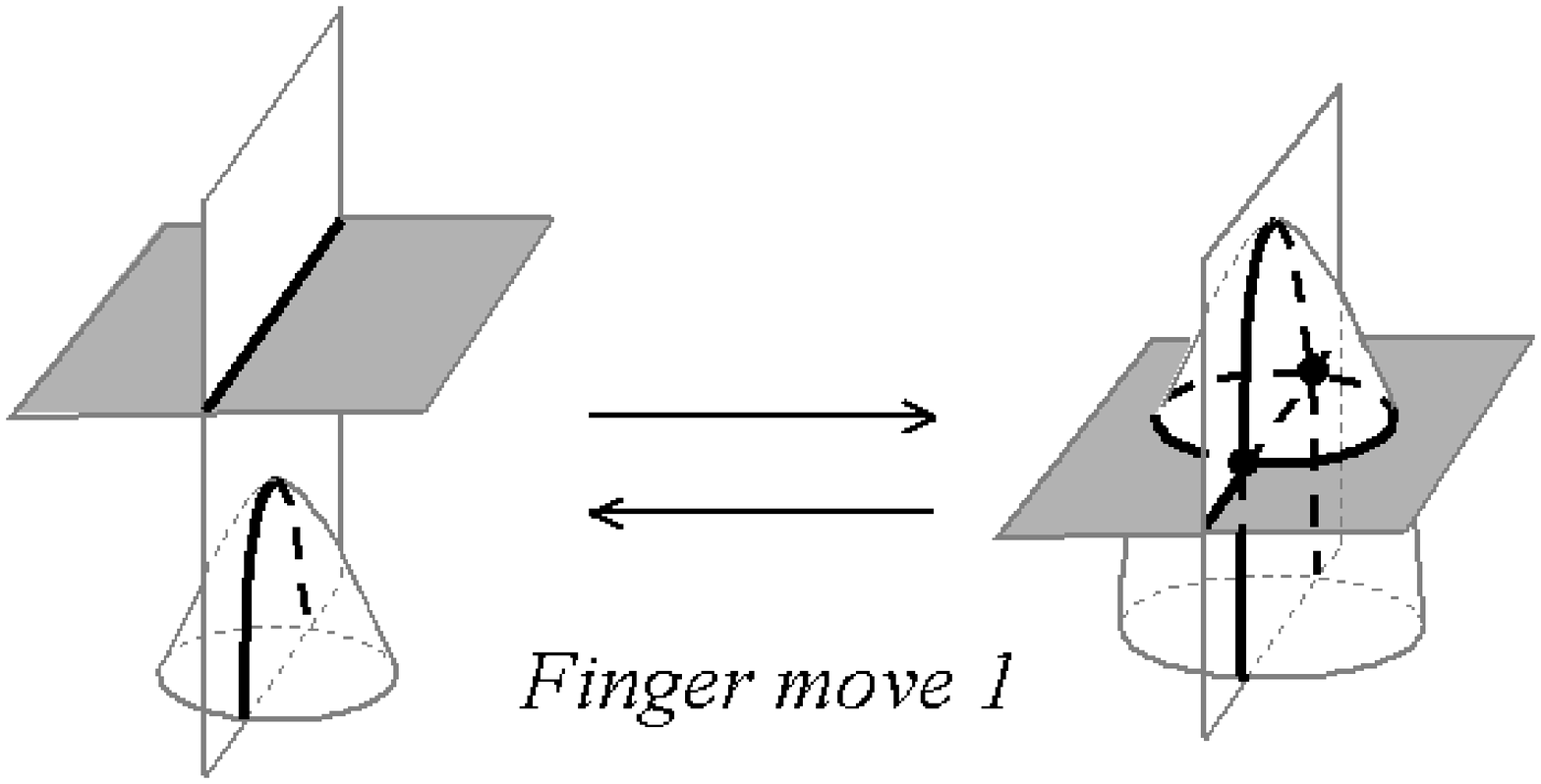}}\\
\subfigure[]{ \label{fig2c}
\includegraphics[
width=0.5\textwidth]{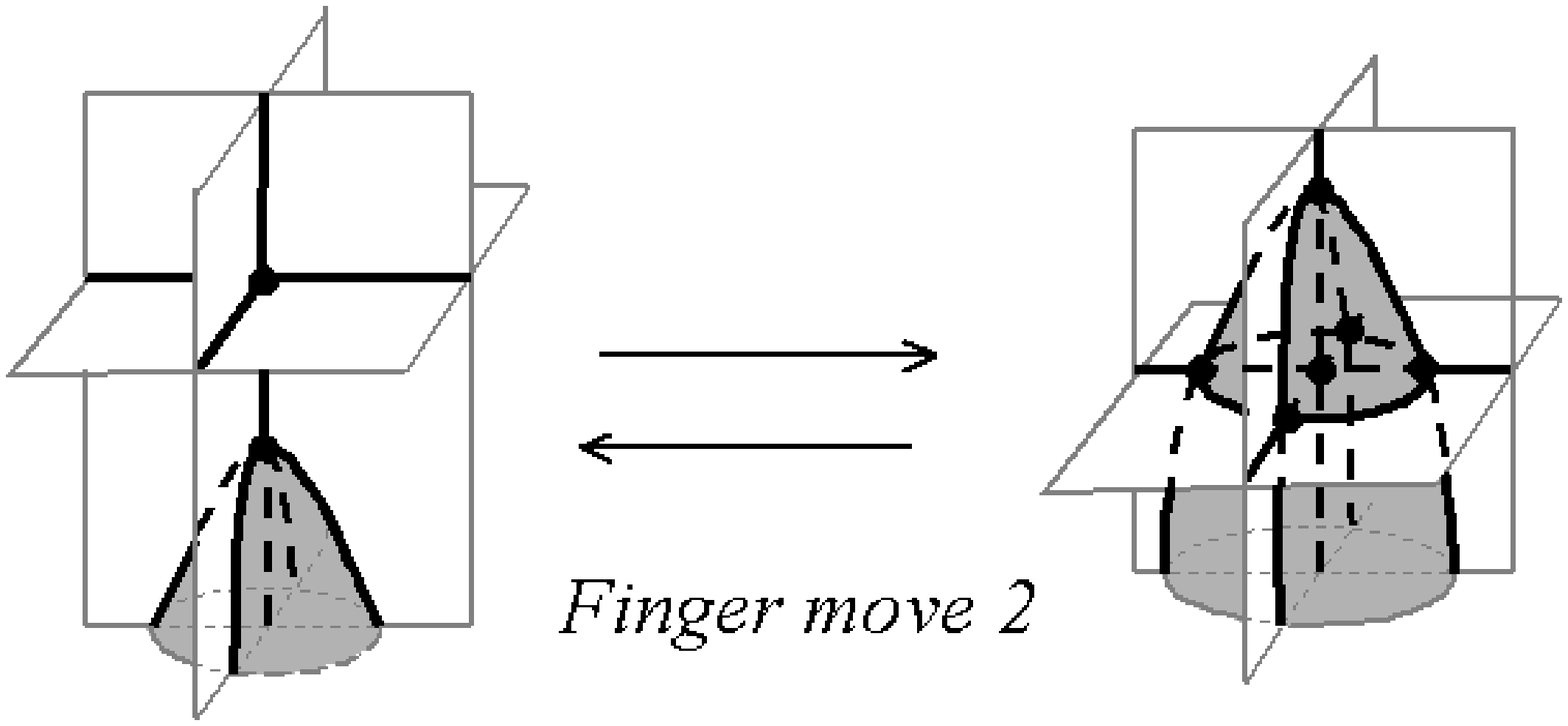}} \caption { }
\end{figure}

This theorem gives a complete set of moves for relating Johansson
representations of the same 3-manifold coming from nulhomotopic
filling Dehn spheres (see Corollary \ref{CORtoMAINtheorem}).

The detailed proof of Theorem \ref{MAINtheorem} is quite long, and
uses both smooth and combinatorial techniques. In this paper we
will give an outline of this proof. The paper is organized as
follows.

In section \ref{SECTION Preliminaries}, we give some preliminary definitions
about Dehn surfaces and cell complexes. Most of section \ref{SECTION
FillingHomotopy} and sections \ref{SECTION Pushing Disks} to \ref{SECTION
Growing} introduce some partial results for giving a sketch of the proof of
Theorem \ref{MAINtheorem}. This sketch of the proof is given in section
\ref{SECTION ProofofMainTheorem}. A reader wishing to skip details can jump
directly from sections \ref{SECTION FillingHomotopy} to \ref{SECTION Diagrams}.

The proof of Theorem \ref{MAINtheorem} in \cite{RHomotopies} relies on three
Key Lemmas that we will state here without proof. In section \ref{SECTION
FillingHomotopy} we present some results about regular homotopies of
immersions of surfaces in 3-manifolds, and we introduce the concept of
\textit{filling-preserving moves} and \textit{filling homotopy} for filling
immersions. Key Lemma 1 is stated in section \ref{SECTION Pushing Disks},
where we define the modifications of immersions of surfaces in 3-manifolds by
\textit{pushing-disks}. These kind of modifications was defined in
\cite{HommaNagase}, and Key Lemma 1 assert that every regular homotopy can be
decomposed into pushing disks with some nice properties.

In section \ref{SECTION Spiral Piping} we introduce a surgery method for
modifying Dehn surfaces that will be useful later, and in section \ref{SECTION
What can be done} we present three examples of modifications of filling Dehn
surfaces that can be done using only the filling-preserving moves defined in
section \ref{SECTION FillingHomotopy}.

In section \ref{SECTION Shellability} we introduce some combinatorial tools
that will be essential in Key Lemmas 2 and 3: the concept of \textit{shelling}
of a cell complex and the concept of \textit{simplicial collapsing} for a
simplicial complex. These concepts will appear almost everywhere in sections
\ref{SECTION InflatingTriangulations} to \ref{SECTION ProofofMainTheorem}. In
the same section \ref{SECTION Shellability} we introduce smooth triangulations
of manifolds, that give us the theorethical basis for applying these
previously defined combinatorial concepts to our case.

We explain in section\ref{SECTION InflatingTriangulations} how any smooth
triangulation $T$ of a 3-manifold $M$ can be ''inflated'' to obtain a filling
Dehn sphere of $M$. This inflated filling Dehn sphere of $M$ has the property
that it is transverse to any Dehn sphere of $M$ that lies in the 2-skeleton of
$T$. When the triangulation $T$ of $M$ is ''sufficiently good'' with respect
to a filling Dehn surface $\Sigma$ of $M$ we can use it to obtain from
$\Sigma$ another filling surfaces ''as complicated as we want'' using only
filling-preserving moves. These constructions are used in Key Lemma 2, which
is stated in the same section \ref{SECTION InflatingTriangulations}. Key Lemma
2 gives a method for transforming any pushing disk (as defined in section
\ref{SECTION Pushing Disks}) in a pushing disk which can be performed used
only filling-preserving moves by putting many new sheets of the Dehn surface
in the middle of the pushing ball (also defined in section \ref{SECTION
Pushing Disks}). The basis for proving Key Lemma 2 in \cite{RHomotopies} is
that when a pushing disk is ''sufficiently good'', it can be performed using
only filling-preserving moves.

In section \ref{SECTION FillingPairs} we discuss briefly how two
filling Dehn spheres of the same 3-manifold can intersect each
other, and this discussion is used in section \ref{SECTION
Growing}, where we state Key Lemma 3. Key Lemma 3 assures that
when two filling Dehn surfaces intersect in a ''sufficiently
good'' way, the inflating constructions introduced in section
\ref{SECTION InflatingTriangulations} can be made simultaneously
for one of them and for the union of both.

All the constructions that we have mentioned above as
''sufficiently good'' are intimately related with the concept of
shelling.

In section \ref{SECTION Diagrams} we translate Theorem \ref{MAINtheorem} for
Johansson representations of 3-manifolds and we give some examples, and in
section \ref{SECTION Duplication} we explain briefly how we can obtain a
nulhomotopic Johansson representation of a 3-manifold $M$ from any Johansson
representation of $M$.

In the final section \ref{SECTION Miscelany} we make a brief discussion about
some related problems.

This paper is part of the Ph. D. Thesis of the author, which has been done
under the supervision of Prof. J. M. Montesinos. I'm very grateful to him for
all his valuable advices, specially for his suggestions and comments during
the writing of this paper and his careful reading of the previous versions of
this manuscript.

\section{Preliminaries.\label{SECTION Preliminaries}}

Because our starting point is Theorem \ref{ThmHassHughes} below, we will work
in the smooth category. Nevertheless, if one can check that the analogue of
Theorem \ref{ThmHassHughes} in the PL category is true (we don't know of any
reference), all our constructions have their translation to the PL case and so
Theorem \ref{MAINtheorem} would also be true in the PL case.

Thus, all the manifolds are assumed equipped with a smooth
structure and maps between two manifolds are assumed to be smooth
with respect to their respective smooth structures.

For the standard definitions of differential topology (immersions,
tranversality, etc.), see \cite{Hirsch} or \cite{Guillemin}, for
example. For a general treatment about PL topology we refer to
\cite{Hudson}, for example.

For a subset $X$ of a manifold, we denote the interior, the closure and the
boundary of $X$ by $int(X)$, $cl(X)$ and $\partial X$ respectively.

Let $A$ and $B$ be two sets. For a map $f:A\rightarrow B$ the \textit{singular
values} or \textit{singularities} of $f$ are the points $x\in B$ with
$\#\left\{  f^{-1}(x)\right\}  >1$, and the \textit{singular points} of $f $
are the inverse image points by $f$ of the singularities of $f$. The
\textit{singular set} $S(f)$ of $f$ is the set of singular points of $f$ in
$A$, and the \textit{singularity set} $\bar{S}(f)$ of $f$ is the set of
singularities of $f$ in $B$. Of course $f(S(f))=\bar{S}(f)$. This is a
notation similar but slightly different to that of \cite{A.Shima2}.

From now on, $M$ will denote a 3-manifold in the conditions
indicated at the beginning of section \ref{SECTION Introduction}.

Let $\Sigma$ be a Dehn surface in $M$.

Let $S$ be a surface and $f:S\rightarrow M$ a transverse immersion
parametrizing $\Sigma$. In this case we say that the surface $S$
is the \textit{domain} of $\Sigma$. For any $x\in M$ it is
$\#\left\{  f^{-1} (x)\right\}  \leq3$ \cite{Hempel}. The
singularities of $f$ are divided into double points of $f$, with
$\#\left\{  f^{-1}(x)\right\}  =2$, and triple points of $f$, with
$\#\left\{  f^{-1}(x)\right\}  =3$. A small neighbourhood of a
double or a triple point looks like in Figures \ref{fig1a} and
\ref{fig1b} respectively. The singularity set $\bar{S}(f)$ of $f$,
the set of triple points of $f$, and the domain $S$ (up to
homeomorphism) do not depend upon the parametrization $f$ of
$\Sigma$ we have chosen. We define the singularity set of
$\Sigma$, and we denote it by $\bar{S}(\Sigma)$, as the
singularity set of any parametrization of $\Sigma$. A
\textit{double curve} of $\Sigma$ is the image of an immersion
$\bar{\gamma}:S^{1}\rightarrow M$ contained in the singularity set
of $\Sigma$ \cite{A.Shima2}. The singularity set of $\Sigma$ is
the union of the double curves of $\Sigma$. Because $S$ is
compact, $\Sigma$ has a finite number of double curves. Following
\cite{A.Shima2}, we denote by $T(\Sigma)$ the set of triple points
of $\Sigma$. The Dehn surface $\Sigma$ is \textit{embedded} if its
singularity set is empty. A \textit{standardly embedded sphere} in
$M$ is a 2-sphere embedded in $M$ that bounds a 3-ball in $M$.

A \textit{component} of $\Sigma$ is the image by $f$ of a connected component
of the domain $S$. Note that the components of $\Sigma$ may not coincide with
the connected components of $\Sigma$.

A Dehn surface $\Sigma$ in $M$ fills $M$ if it defines a cell-decomposition of
$M$ as it has been indicated in section \ref{SECTION Introduction}. This
definition generalizes to general surfaces a definition given in
\cite{Montesinos} for Dehn spheres.

The following Proposition gives an equivalent definition of filling Dehn surface.

\begin{proposition}
\label{PROPfillsMifandonlyIF}$\Sigma$ fills $M$ if and only if

\begin{enumerate}
\item $M-\Sigma$ is a disjoint union of open 3-balls,

\item $\Sigma-\bar{S}(\Sigma)$ is a disjoint union of open 2-disks, and

\item $\bar{S}(\Sigma)-T(\Sigma)$ is a disjoint union of open intervals.
\end{enumerate}
\end{proposition}

The following statements and definitions that we will give now for
cell complexes are valid also for simplicial complexes. We
consider the cells of a cell complex as \textit{open} cells. If
$K$ is a cell complex, and $\epsilon,\epsilon^{\prime}$ are two
cells of $K$, we will denote $\epsilon<\epsilon^{\prime}$ when
$\epsilon$ is a face of $\epsilon^{\prime}$, that is, when it is
$cl(\epsilon)\subset cl(\epsilon^{\prime})$. The cells $\epsilon$
and $\epsilon^{\prime}$ are \textit{incident} if $\epsilon
<\epsilon^{\prime}$ or $\epsilon^{\prime}<\epsilon$, and they are
\textit{adjacent} if $cl(\epsilon)\cap
cl(\epsilon^{\prime})\neq\varnothing$. For a cell $\epsilon$ of
$K$, we define the (open) \textit{star} of $\epsilon$ as the union
of all the cells $\epsilon^{\prime}$ of $K$ with $\epsilon
<\epsilon^{\prime}$. The star of $\epsilon$ is denoted by\textit{\
} $star(\epsilon)$.

If $\epsilon$ is a cell of the cell complex $K$, and $P$ is a vertex (0-cell)
of $\epsilon$, we say that $\epsilon$ is \textit{self-adjacent} at $P$ if a
regular neighbourhood of $P$ in $K$ intersects $\epsilon$ in more than one
connected component. Otherwise we say that $\epsilon$ is \textit{regular at
}$P$. We say that $\epsilon$ is \textit{regular} if it is regular at every
vertex of $\epsilon$. The complex $K$ is \textit{regular at }$P$ if every cell
of $K$ incident with $P$ is regular at $P$, and $K$ is \textit{regular} if
every cell of $K$ is regular (compare \cite{Massey}). A filling Dehn surface
$\Sigma$ of $M$ is \textit{regular (regular at a triple point)} if the cell
decomposition of $M$ that defines $\Sigma$ is regular (at this triple point).

If $\Sigma$ is a filling Dehn surface, then a connected component of
$M-\Sigma$ is called a \textit{region} of $M-\Sigma$, and a connected
component of $\Sigma-\bar{S}(\Sigma)$ is sometimes called a \textit{face} of
$\Sigma$.

\section{Filling homotopy.\label{SECTION FillingHomotopy}}

An \textit{ambient isotopy} of a manifold $N$ is a map $\varsigma
:N\times\left[  0,1\right]  \rightarrow N$ such that
$\varsigma_{t} =\varsigma(\cdot,t)$ is a diffeomorphism for each
$t\in\left[  0,1\right]  $ and $\varsigma_{0}=id_{N}$. Two
immersions $f,g:S\rightarrow M$ are \textit{ambient isotopic in
}$M$ if there is an ambient isotopy $\bar {\varsigma}$ of $M$ with
$\bar{\varsigma}\circ f=g$. The same immersions are
\textit{ambient isotopic in }$S$ if there is an ambient isotopy
$\varsigma$ of $S$ with $f\circ\varsigma=g$. We generally say that
$f$ and $g$ are \textit{ambient isotopic} if they are related by
ambient isotopies of $S$ and ambient isotopies of $M$.

Two immersions $f,g:S\rightarrow M$ from a surface $S$ into the 3-manifold $M
$, are \textit{regularly homotopic} if there is an homotopy $H:S\times\left[
0,1\right]  \rightarrow M$ with $H(\cdot,0)=f$ and $H(\cdot,1)=g$, and such
that $H(\cdot,t)$ is an immersion for each $t\in\left[  0,1\right]  $. The
homotopy $H$ defines a smooth path of immersions from $S$ into $M$ having $f$
and $g$ as its endpoints. If $f$ and $g$ are regularly homotopic, they are
indeed homotopic. The converse is not true in general. Nevertheless, an
immediate corollary of Theorem 1.1 in \cite{HassHugues} or Theorem 6 in
\cite{Li Banghe} is the next Theorem

\begin{theorem}
\label{ThmHassHughes}Two immersions from the 2-sphere $S^{2}$ into a
3-manifold are regularly homotopic if and only if they are homotopic.
\end{theorem}

In particular, two parametrizations of nulhomotopic filling Dehn
spheres of $M$ must be regularly homotopic.

\begin{figure}[b]
\centering \subfigure []{ \label{fig3a}
\includegraphics[
width=0.5\textwidth]{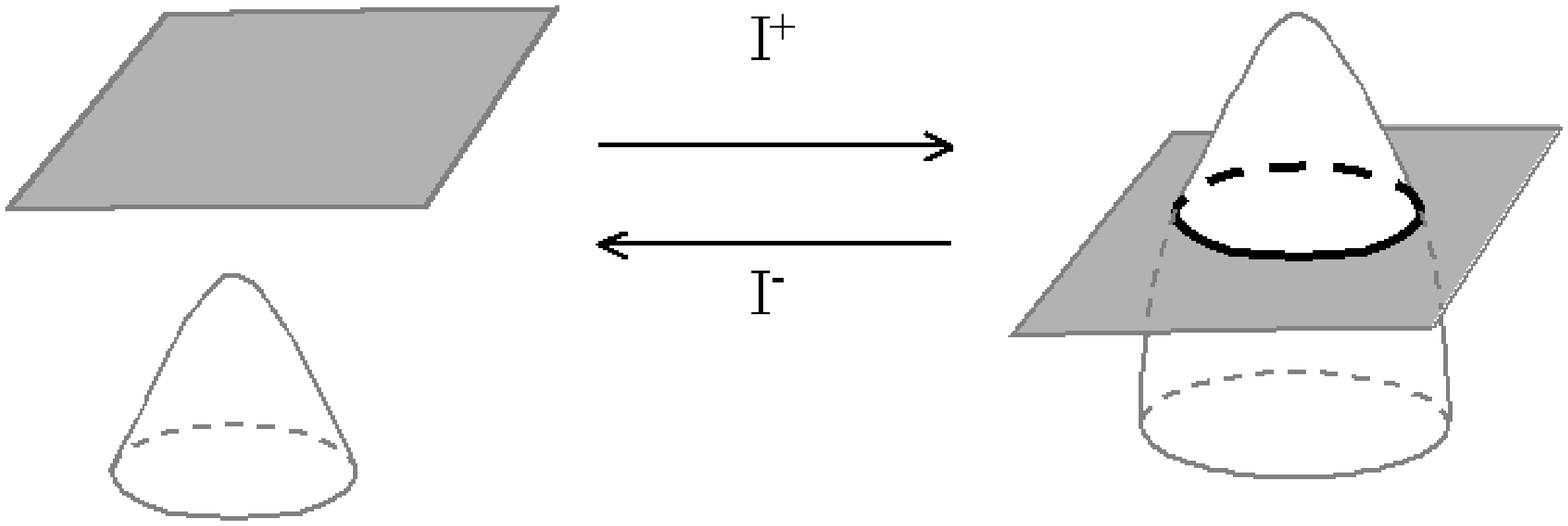}} \caption { }
\end{figure}

In \cite{HommaNagase} is introduced a set of \textit{elementary
deformations} for immersions of surfaces in 3-manifolds. This set
of moves is composed by the \textit{saddle move} (which is called
\textit{elementary deformation of type VI}, in the notation of
\cite{HommaNagase}) of Figure \ref{fig2a}, together with the moves
depicted in Figure 3 (Figures \ref{fig3a}, \ref{fig3b} and
\ref{fig3c} have been taken from \cite{Yashiro}). We will call
these elementary deformations the \textit{Homma-Nagase moves}. In
\cite{HommaNagase2} it is proved the following Theorem

\begin{theorem}
\label{THM hommaNagase}Two transverse immersions from a closed surface $S$
into a 3-manifold $M$ are regularly homotopic if and only if we can deform
them into one another by a finite sequence of Homma-Nagase moves, together
with ambient isotopies of $M$.
\end{theorem}

The proof of this Theorem in \cite{HommaNagase2} is in the PL category. A
proof of the smooth version of this result is indicated in \cite{Roseman}. An
equivalent result, also in the differentiable case is Theorem 3.1 of
\cite{HassHugues}.

\addtocounter{figure}{-1}
\begin{figure} [htb]
\addtocounter{subfigure}{1} \centering \subfigure[]{ \label{fig3b}
\includegraphics[
width=0.5\textwidth]{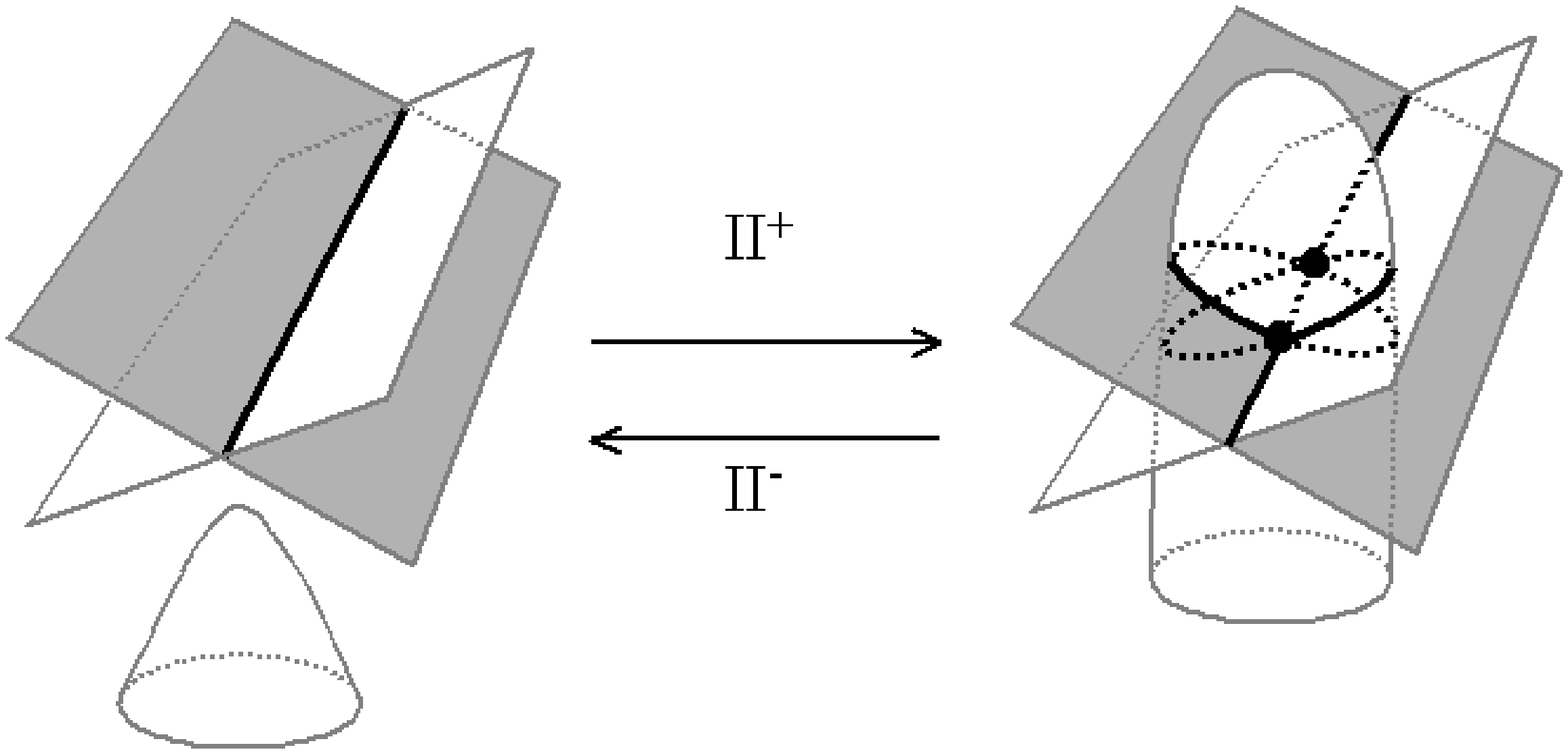}}\\
\subfigure[]{ \label{fig3c}
\includegraphics[
width=0.5\textwidth]{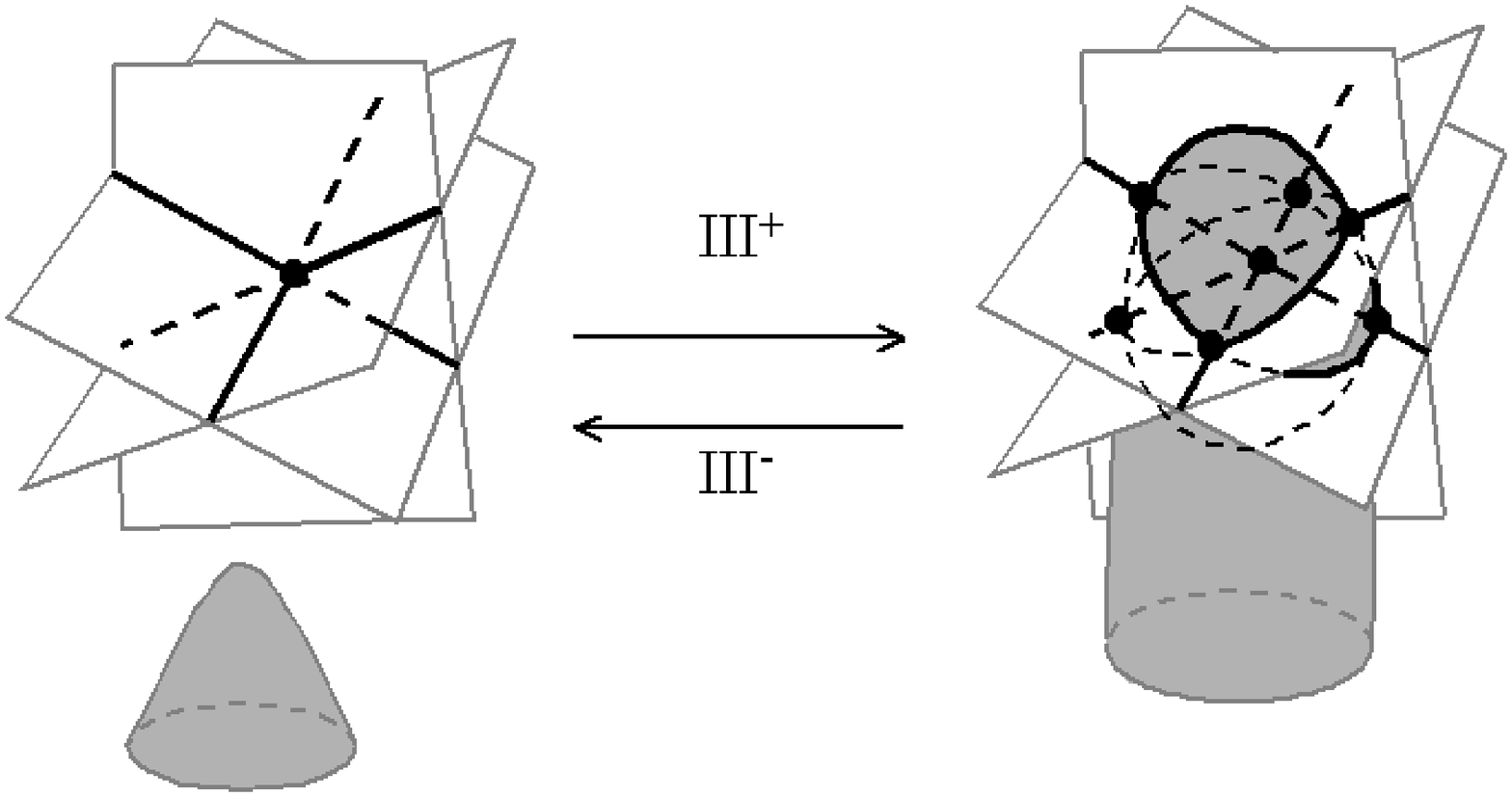}} \caption { }
\end{figure}

We will propose another set of moves (Haken moves), which is the
result of substituting in the Homma-Nagase set of moves the moves
of type II and III by the \textit{finger moves} 1 and 2 depicted
in Figures \ref{fig2b} and \ref{fig2c} respectively. The following
Lemma can be easily proved:

\begin{lemma}
\label{EQUIVALENCEHOmma_NAGASE-HAKEN}The Homma-Nagase set of moves and the
Haken set of moves are equivalent.
\end{lemma}

For proving this Lemma it must be shown that each Homma-Nagase
move can be obtained using Haken moves (and ambient isotopies) and
viceversa. Thus, in Theorem \ref{THM hommaNagase} we can
substitute the Homma-Nagase moves by the Haken moves.

The Haken moves fit better than the Homma-Nagase moves when we are dealing
with filling Dehn surfaces. In the Haken set of moves, the Homma-Nagase move
of type $I$ is called \textit{finger move 0}. For $i=0,1,2$ a finger move $i$
is a finger move $+i$ when it happens from left to right in the figure, and it
is a finger move $-i$ if it happens in the opposite sense. A saddle move is
equivalent (symmetric) in both senses.

\begin{lemma}
\label{LemmaTYPE???ISFILLINGPRESERVING}Let $f,g:S\rightarrow M$ be two
immersions. Then:

\begin{enumerate}
\item  if $f$ and $g$ are related by a finger move 0, then one of them is not
a filling immersion;

\item  if $f$ and $g$ are related by a finger move 1 or 2, then $f$ is a
filling immersion if and only if $g$ is a filling immersion; and

\item  if $f$ and $g$ are related by a saddle move and $f$ is a filling
immersion, then $g$ is not necessarily a filling immersion.
\end{enumerate}
\end{lemma}

This Lemma can be proved by inspection, using the characterization
of filling immersions given by Proposition
\ref{PROPfillsMifandonlyIF}.

Lemma \ref{LemmaTYPE???ISFILLINGPRESERVING} inspired us the following
definition. If $f:S\rightarrow M$ is a filling immersion and we modify $f$ by
a Haken move, we say that that move is \textit{filling-preserving} if the
immersion $g$ we get after the move is again a filling immersion. With this
notation, Lemma \ref{LemmaTYPE???ISFILLINGPRESERVING} means that a finger move
0 cannot be filling-preserving, finger moves 1 and 2 are always
filling-preserving, and saddle moves are sometimes filling-preserving and
sometimes not. The next step is the following Definition:

\begin{definition}
\label{DEFfillinghomotopic}Let $f,g:S\rightarrow M$ be two filling immersions.
We say that $f$ and $g$ are \textit{filling homotopic} if there is a finite
sequence $f=f_{0},f_{1},...,f_{n}=g$ of immersions such that for each
$i=0,...,n-1$ the immersions $f_{i}$ and $f_{i+1}$ are ambient isotopic or
they are related by a filling-preserving move.
\end{definition}

Note that in the previous Definition, all the immersions of the sequence
$f_{0},...,f_{n}$ are filling immersions. With this notation, Theorem
\ref{MAINtheorem} can be restated as follows

\begin{theorem}
If $f,g:S^{2}\rightarrow M$ are nulhomotopic filling immersions, then they are
filling homotopic.
\end{theorem}

This gives a partial answer to the following Conjecture:

\begin{conjecture}
\label{CONJECture}Regularly homotopic filling immersions of arbitrary surfaces
are filling homotopic.
\end{conjecture}

Perhaps the proof of Theorem \ref{MAINtheorem} given in \cite{RHomotopies} and
sketched here can be adapted to a more general case but we still don't know
how to do this.

\section{Pushing disks.\label{SECTION Pushing Disks}}

Let $f,g:S\rightarrow M$ be two immersions. Assume that there is a closed disk
$D\subset S$ such that:

\begin{enumerate}
\item $f$ and $g$ agree in $S-D$;

\item $f\mid_{D}$ and $g\mid_{D}$ are both embeddings;

\item $f\left(  D\right)  $ and $g\left(  D\right)  $ intersect only in
$f\left(  \partial D\right)  =g\left(  \partial D\right)  $; and

\item $f\left(  D\right)  \cup g\left(  D\right)  $ bounds a 3-ball $B$ in $M$
(Figure 4).
\end{enumerate}

\begin{figure} [htb]
\centering \subfigure[]{ \label{fig4a}
\includegraphics[
width=0.45\textwidth]{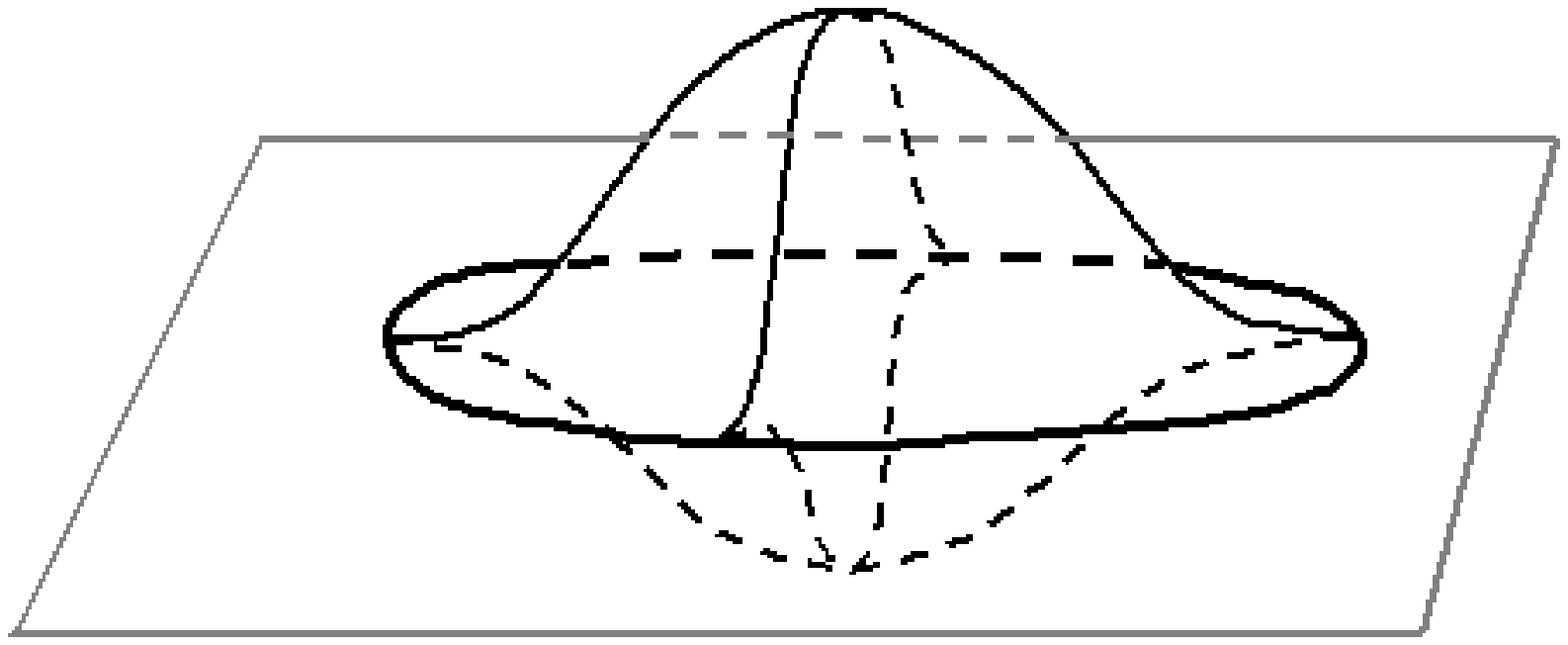}}\hfill \subfigure[]{ \label{fig4b}
\includegraphics[
width=0.45\textwidth]{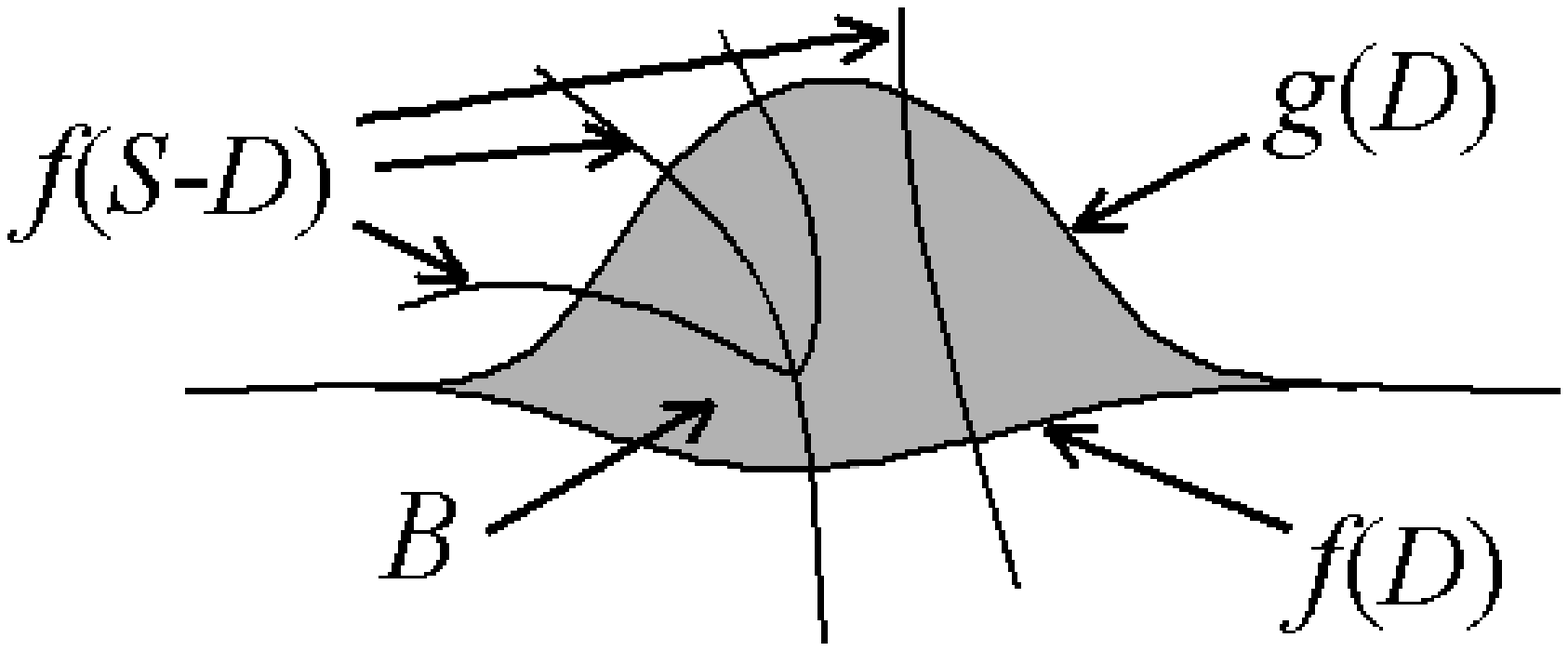}} \caption { }
\end{figure}

Then we say that $g$ is obtained from $f$ by pushing the disk $D$
through $B$ or along $B$ (see Figure 4). The pair $\left(
D,B\right)  $ is a \textit{pushing disk} (see \cite{HommaNagase}).
In the pushing disk $\left( D,B\right)  $, the disk $D$ is
\textit{the pushed disk}, $B$ is the \textit{pushing ball} and we
say also that $f\left(  \partial D\right) =g\left(  \partial
D\right)  $ is the \textit{equator} of $B$ and it is denoted by
$eq(B)$. If both $f$ and $g$ are transverse immersions, we say
that the pushing disk $\left(  D,B\right)  $ is
\textit{transverse}. In the pushing disk $\left(  D,B\right)  $,
the ''rest'' of the immersed surface, $f(S-D)$, may intersect $B$
in any manner (Figure \ref{fig4b}). If we are given the immersion
$f$ and the pushing disk $\left(  D,B\right)  $, then the
immersion $g$ is well defined up to an ambient isotopy of $S$.

We will say that two (transverse) immersions $f,g:S\rightarrow M$ are
\textit{regularly homotopic by (transverse) pushing disks} if there is a
finite sequence $f=f_{0},f_{1},...,f_{n}=g$ of (transverse) immersions such
that $f_{i}$ is obtained from $f_{i-1}$ by a pushing disk for $i=1,...,n$.

The first Step in the proof of Theorem \ref{MAINtheorem} is the following
Lemma, whose proof is in \cite{RHomotopies}.

\begin{lemma}
[Key Lemma 1]\label{KeyLemma1}Let $f,g:S\rightarrow M$ be two regularly
homotopic immersions, then:

\begin{enumerate}
\item [ A)]they are regularly homotopic by pushing disks.

\item[ B)] if they are transverse, then they are regularly homotopic by
transverse pushing disks.

\item[ C)] if $f$ and $g$ agree over a disk $D$ of $S$, such that the
restrictions $f\mid_{D}=g\mid_{D}$ are embeddings, then in cases A) and B) the
pushing disks can be chosen keeping $D$ fixed.
\end{enumerate}
\end{lemma}

Thus, we can decompose any regular homotopy into a finite sequence of pushing
disks. Note that the Homma-Nagase moves and the Haken moves are special kinds
of transverse pushing disks. However, Theorem \ref{THM hommaNagase} decomposes
a regular homotopy into transverse pushing disks \underline{and} ambient
isotopies of $M$. Disposing of this ambient isotopy is the hardest part in the
proof of Key Lemma 1 in \cite{RHomotopies}. In the same way as an immersion
behaves locally as an embedding, a regular homotopy behaves locally as an
isotopy. Using this, the proof of Key Lemma 1 will be obtained after a
detailed study of isotopies of embedded surfaces in 3-manifolds, and it is
mainly inspired in \cite{HudsonZeeman}.

\section{Spiral piping.\label{SECTION Spiral Piping}}

\begin{figure} [htb]
\centering \subfigure[]{ \label{fig5a}
\includegraphics[
width=0.35\textwidth]{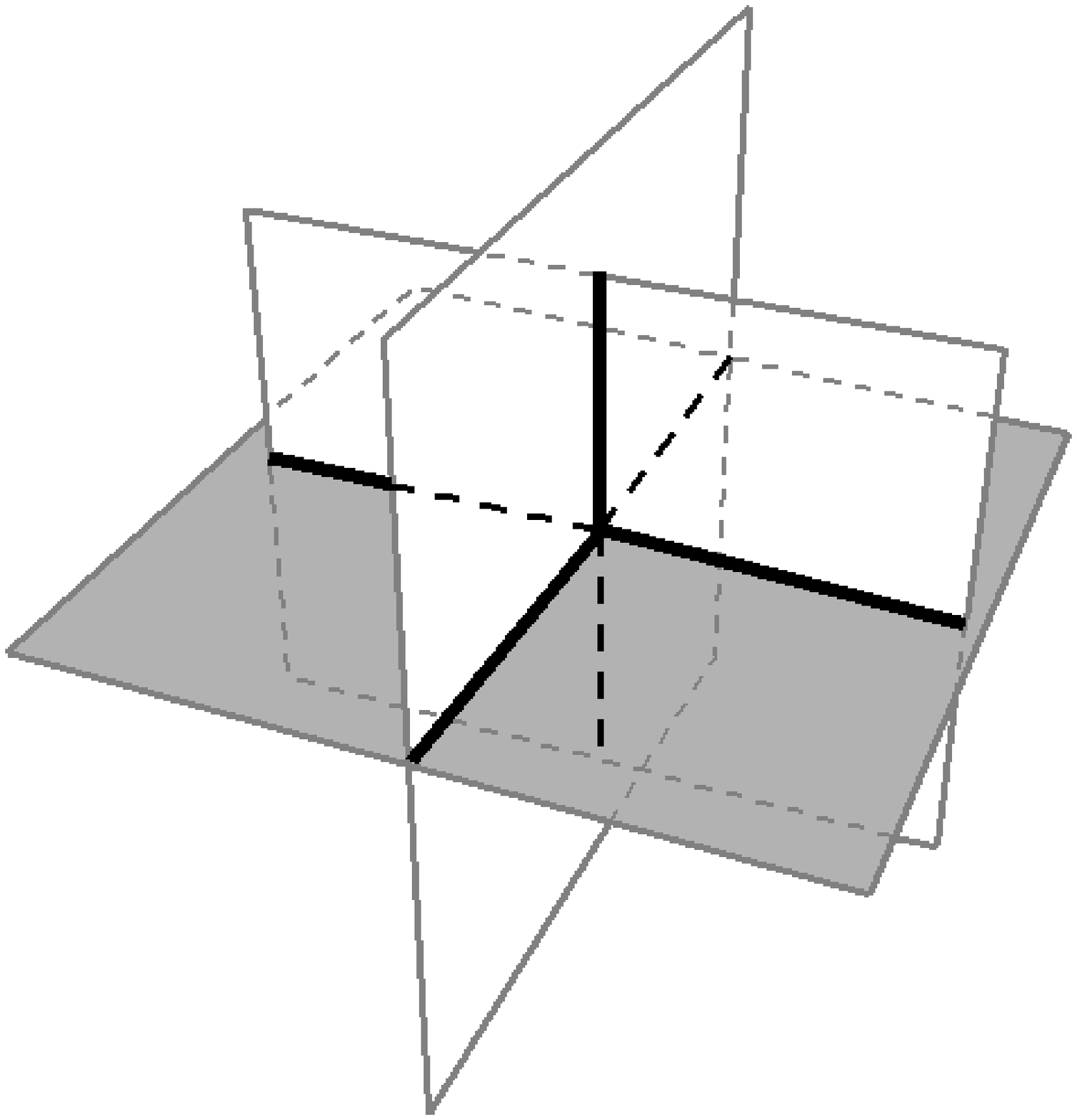}}\hfill \subfigure {
\includegraphics[ height=0.25\textheight, width=0.05\textwidth ] {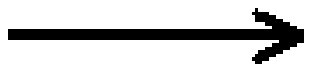} }\hfill
\addtocounter {subfigure}{-1} \subfigure[]{ \label{fig5b}
\includegraphics[
width=0.35\textwidth]{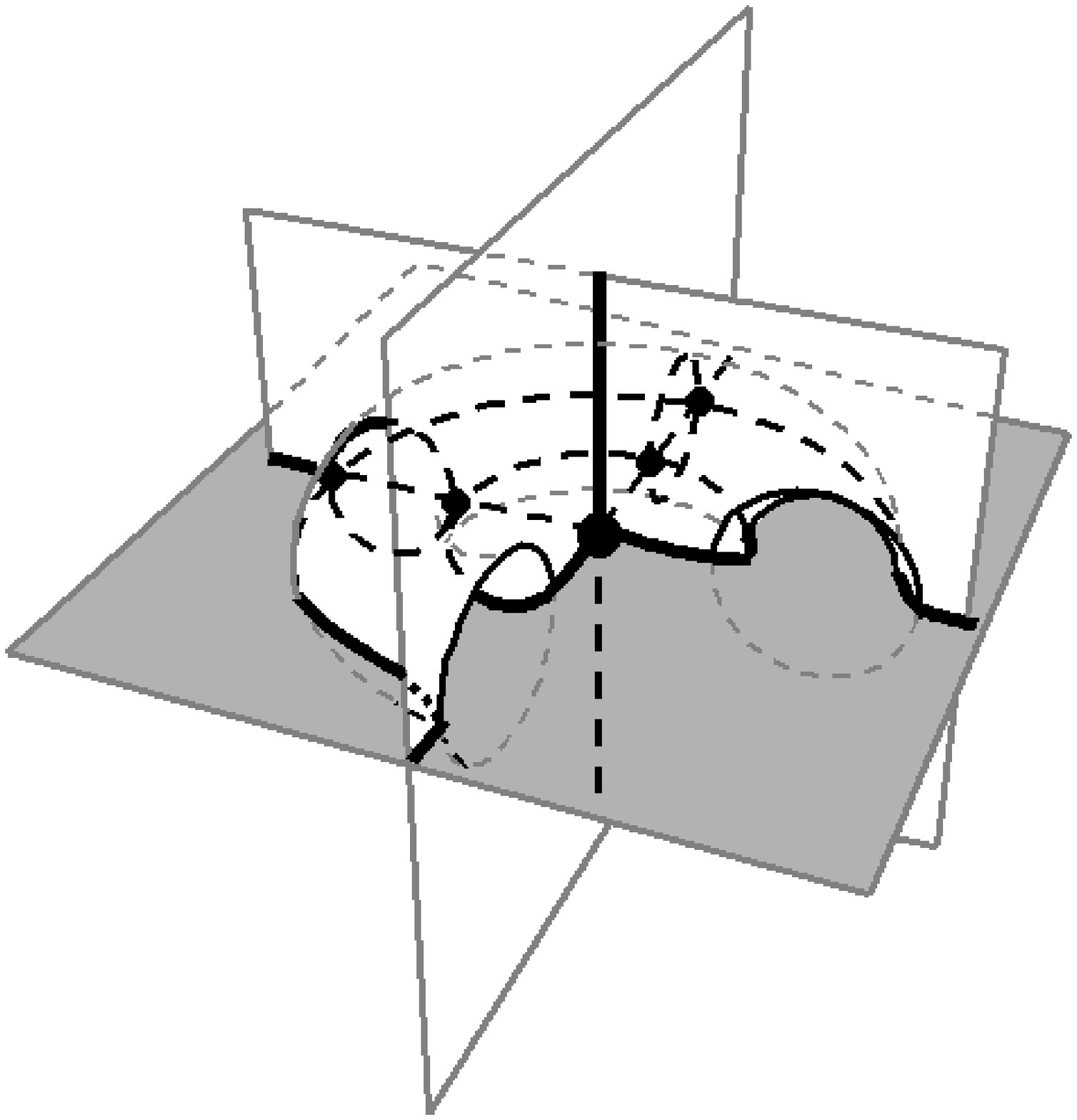}} \caption { }
\end{figure}

In \cite{Banchoff} it is explained how to modify Dehn surfaces by
\textit{surgery}, also called \textit{piping }(see
\cite{RourkeSanderson}, p. 67). We introduce now a special kind of
piping that will be useful later. Let $\Sigma$ be a Dehn surface
in $M$, and let $P$ be a triple point of $\Sigma$. If $P$ is the
triple point depicted in Figure \ref{fig5a}, consider the surface
$\Sigma^{\prime}$ that is exactly identical with $\Sigma$ except
in a neighbourhood of $P$ that can be as small as necessary. In
this neighbourhood of $P$, the Dehn surface $\Sigma^{\prime}$
looks like Figure \ref{fig5b}, and we say that $\Sigma^{\prime}$
is obtained from $\Sigma$ by a \textit{spiral piping around} $P$.

\begin{proposition}
\label{PROPSpiralPipingPreserveFillingness}In this situation, if $\Sigma$ is a
(regular) filling Dehn surface of $M$, then $\Sigma^{\prime}$ is a (regular)
filling Dehn surface of $M$.
\end{proposition}

See \cite{Anewproof} for more details.

If the two sheets of $\Sigma$ that become connected by the piping (the two
vertical sheets in Figure 5) belong to different components $\Sigma_{1}$ and
$\Sigma_{2}$ of $\Sigma$, then after performing the spiral piping these two
components of $\Sigma$ become a unique component $\Sigma_{1}\#\Sigma_{2}$ of
$\Sigma^{\prime}$.

If $S$ is the domain of $\Sigma$, and $S^{\prime}$ is the domain
of $\Sigma^{\prime}$ it is easy to check that $S^{\prime}$ is the
result of removing the interior of two small closed disks
$\delta_{1},\delta_{2}$ from $S$ and identify their boundaries in
an adequate way. If $\delta_{1}$ and $\delta_{2}$ belong to
different connected components $S_{1},S_{2}$ of $S$ respectively,
then $S^{\prime}$ is the result of substituting the union
$S_{1}\cup S_{2}$ in $S$ by the connected sum $S_{1}\#S_{2}$.

The following Definition and Theorem appear in \cite{Anewproof}.

\begin{definition}
\label{DEFfillingCOLLECTION}A Dehn surface $\Sigma\subset M$ that fills $M$ is
called a \textit{filling collection of spheres} in $M$ if its domain is a
disjoint union of a finite number of 2-spheres.
\end{definition}

\begin{theorem}
\label{THM fillingcollection implies fillingSPHERE}If $M$ has a filling
collection of spheres $\Sigma$, then $M$ has a filling Dehn sphere
$\Sigma^{\prime}$. If each component of $\Sigma$ is nulhomotopic, we can
choose $\Sigma^{\prime}$ to be nulhomotopic.
\end{theorem}

\begin{proof}
Let $\Sigma$ be a filling collection of spheres in $M$, and let $\Sigma
_{1},...,\Sigma_{m}$ be the different components of $\Sigma$.

The 2-skeleton of any cell decomposition of $M$ is connected
because $M$ is connected. Then, $\Sigma$ is connected. Thus, we
can assume that $\Sigma_{1},...,\Sigma_{m}$ are ordered in such a
way that $\Sigma_{1} \cup...\cup\Sigma_{k}$ is connected for every
$k\in\left\{  1,...,m\right\} $. In particular, $\Sigma_{k}$
intersects $\Sigma_{1}\cup...\cup\Sigma_{k-1}$ for all
$k\in\left\{  2,...,m\right\}  $.

Because $\Sigma_{1}\cap\Sigma_{2}$ is nonempty, the intersection
$\Sigma _{1}\cap\Sigma_{2}$ contains a double curve of $\Sigma$,
and because $\Sigma$ fills $M$, this double curve contains at
least one triple point $P$ of $\Sigma$. Connecting $\Sigma_{1}$
and $\Sigma_{2}$ near $P$ by a spiral piping, we obtain a new Dehn
sphere $\Sigma_{1}\#\Sigma_{2}$ such that $\left(
\Sigma_{1}\#\Sigma_{2}\right) \cup\Sigma_{3}\cup...\cup\Sigma_{m}$
still fills $M$.

Because $\Sigma_{3}$ intersects $\Sigma_{1}\cup\Sigma_{2}$, it
intersects $\Sigma_{1}\#\Sigma_{2}$. Where
$\Sigma_{1}\#\Sigma_{2}$ and $\Sigma_{3}$ intersect transversely
there is a triple point of $\Sigma$ (and therefore of $\left(
\Sigma_{1}\#\Sigma_{2}\right)  \cup...\cup\Sigma_{m}$). We can
perform another piping operation (as before) obtaining a new Dehn
sphere $\Sigma_{1}\#\Sigma_{2}\#\Sigma_{3}$, and such that the new
Dehn surface $\left(  \Sigma_{1}\#\Sigma_{2}\#\Sigma_{3}\right)
\cup\Sigma_{4} \cup...\cup\Sigma_{m}$ still fills $M$.

Inductively, for $k>3$, we obtain a Dehn sphere
$\Sigma_{1}\#\Sigma _{2}\#...\#\Sigma_{k}$ piping
$\Sigma_{1}\#\Sigma_{2}\#...\#\Sigma_{k-1}$ with $\Sigma_{k}$
around a triple point of $\Sigma$ lying in the intersection of
$\Sigma_{1}\#\Sigma_{2}\#...\#\Sigma_{k-1}$ and $\Sigma_{k}$, with
the property that $\left(  \Sigma_{1}\#\Sigma_{2}\#...\#\Sigma
_{k}\right)  \cup\Sigma_{k+1}\cup...\cup\Sigma_{m}$ still fills
$M$.

Repeating this operation we finally obtain a Dehn sphere $\Sigma^{\prime
}=\Sigma_{1}\#\Sigma_{2}\#...\#\Sigma_{m}$ that fills $M$.

If all the components of $\Sigma$ are nulhomotopic, this implies
that we can deform the Dehn sphere $\Sigma_{m}$ continuously to a
point. If $g_{m} :S^{2}\rightarrow M$ is an immersion
parametrizing $\Sigma^{\prime}$, we can use this deformation to
construct an homotopy between $g_{m}$ and an immersion $g_{m-1}$
parametrizing $\Sigma_{1}\#\Sigma_{2}\#...\#\Sigma_{m-1}$. In the
same way, we can construct an homotopy between $g_{m-1}$ and an
immersion $g_{m-2}$ parametrizing
$\Sigma_{1}\#\Sigma_{2}\#...\#\Sigma_{m-2}$. Repeating this
process, we finally obtain that $g_{m}$ is homotopic to an
immersion $g_{1}$ parametrizing $\Sigma_{1}$ and so $g_{m}$ is
nulhomotopic.
\end{proof}

\begin{figure} [htb]
\centering \subfigure[]{ \label{fig6a}
\includegraphics[
width=0.35\textwidth]{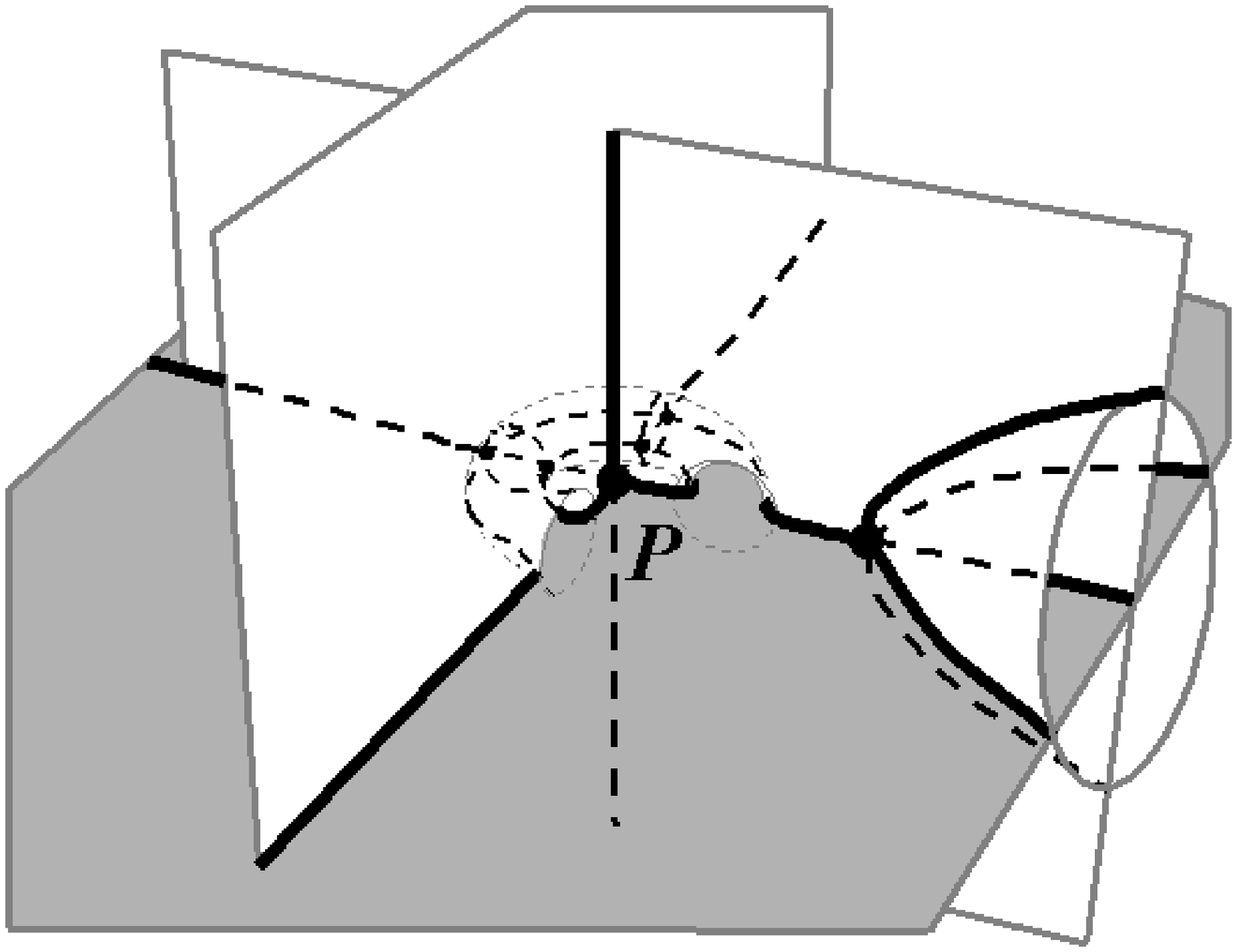}}\hfill \subfigure {
\includegraphics[ height=0.2\textheight, width=0.05\textwidth ] {fig35ab} }\hfill
\addtocounter {subfigure}{-1} \subfigure[]{ \label{fig6b}
\includegraphics[ width=0.35\textwidth]{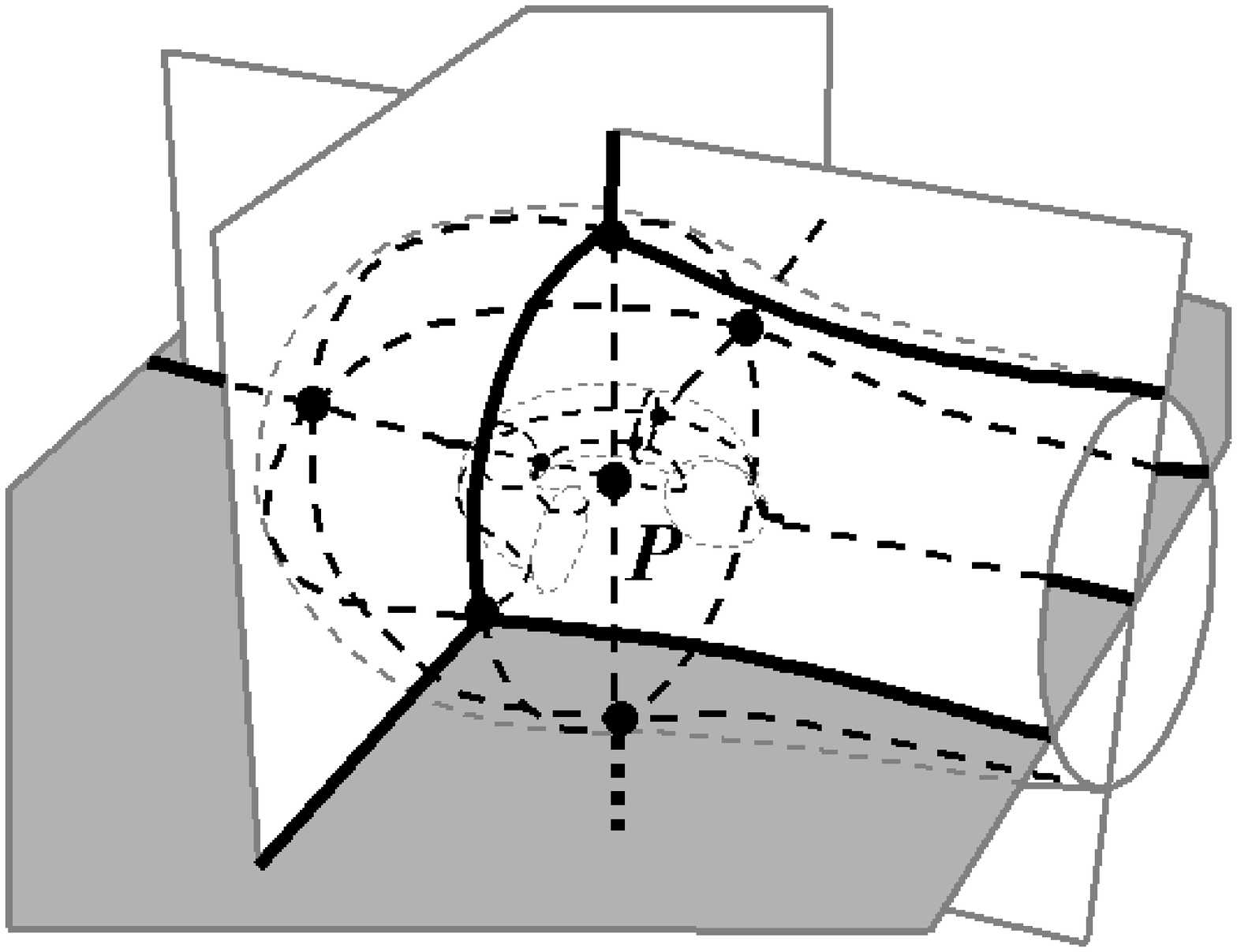}} \caption { }
\end{figure}

Another property of spiral pipings is that \textit{they do not
disturb filling homotopies}, as it is stated in the following
Lemma, that we will state without proof:

\begin{lemma}
\label{LEMMAIgnoringPipings}Let $f,g:S\rightarrow M$ be two
filling immersions such that $g$ is obtained from $f$ after a
finger move +2 through the triple point $P$ of $f$. Let
$S^{\prime},f^{\prime},g^{\prime}$ be the surface and immersions
that come from $S,f,g$ respectively after performing a spiral
piping around $P$. We assume that this spiral piping is as small
as necessary, in comparation with the finger move (Figure
\ref{fig6b}). Then, $f^{\prime}$ and $g^{\prime}$ are filling
homotopic.
\end{lemma}

In the situation explained in the previous Lemma we say that the immersions
$f^{\prime},g^{\prime}$ are related by a \textit{piping passing move} through
$P$.

\section{What can be done using filling-preserving moves.\label{SECTION What
can be done}}

We will give here three examples of operations to be performed in a filling
Dehn surface using only filling-preserving moves.

We start considering a filling Dehn surface $\Sigma$ of the 3-manifold $M$.

\subsection{Inflating a double point.\label{SUBSECTION InflatingDoublePoint}}
\vspace{10pt}
\begin{figure} [htb]
\centering
\includegraphics[ width=0.55\textwidth
]{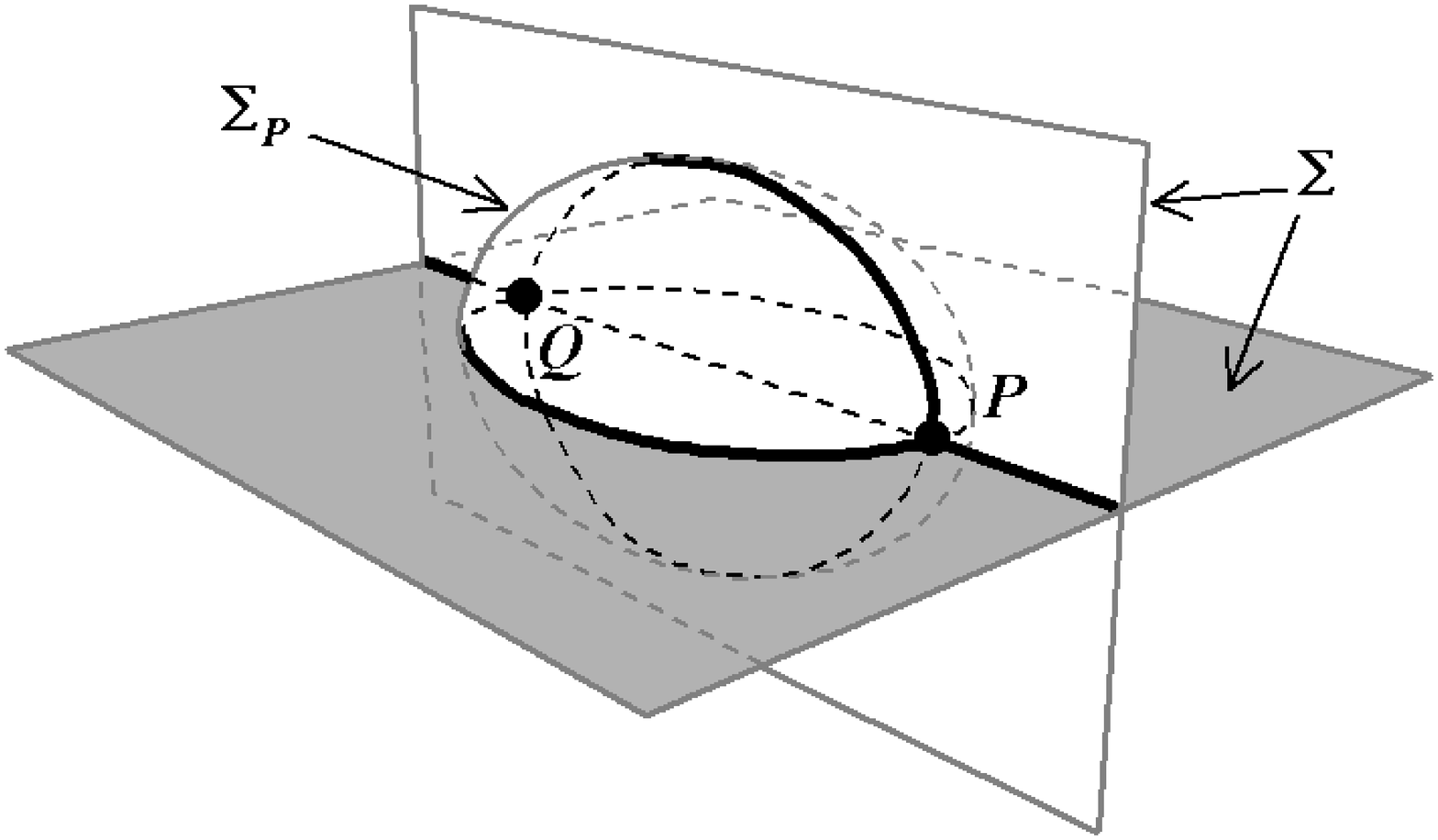} \caption { }\label{fig7}
\end{figure}

Let $P$ be a double point of $\Sigma$. Consider a standardly
embedded 2-sphere $\Sigma_{P}$ in $M$ as in Figure \ref{fig7}. The
sphere $\Sigma_{P}$ contains $P$, and its intersection with
$\Sigma$ is the union of two circles. These two circles intersect
themselves at $P$ and at other point $Q$, which are the unique
double points of $\Sigma$ lying in $\Sigma_{P}$. Note that the
union $\Sigma\cup\Sigma_{P}$ is a filling Dehn surface of $M $.
Consider a filling Dehn surface $\Sigma\#\Sigma_{P}$ that is the
result of modifying $\Sigma \cup\Sigma_{P}$ by a spiral piping
around $P$ (see section \ref{SECTION Spiral Piping}).

\begin{proposition}
\label{PropositionINFLATINGdoublepoint}We can choose the piping such that
$\Sigma$ is filling homotopic to $\Sigma\#\Sigma_{P}$.
\end{proposition}

\begin{proof}
We consider the filling Dehn surface $\Sigma\#\Sigma_{P}$ as in
Figure \ref{fig8a}. This surface is identical with
$\Sigma\cup\Sigma_{P}$ except in a small neighbourhood of $P$,
where it looks like Figure \ref{fig8b}.

\begin{figure}[htb]
\centering \subfigure []{ \label{fig8a}
\includegraphics[ width=0.3\textwidth]{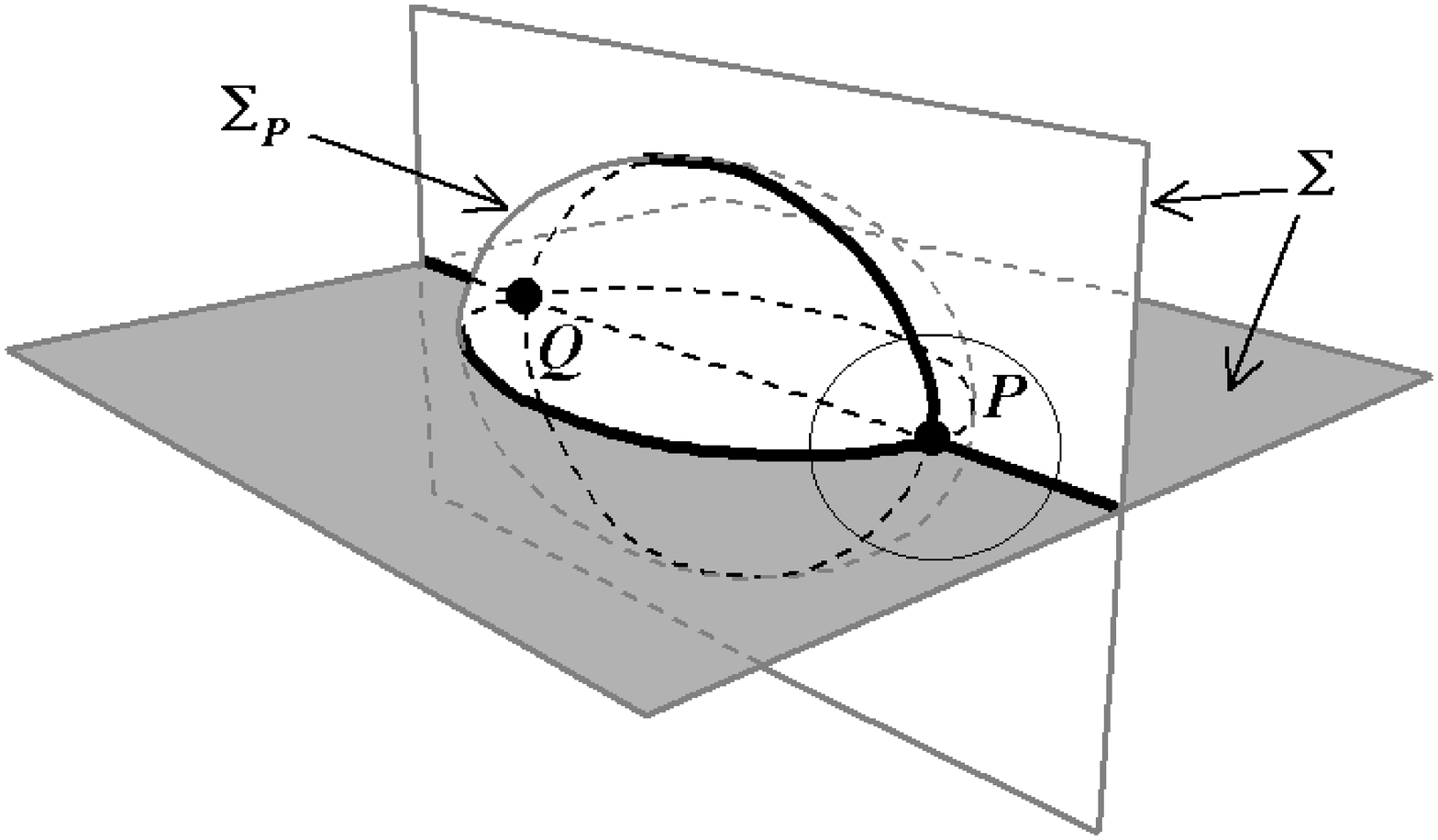}}\hfill
\subfigure[]{ \label{fig8b}
\includegraphics[ width=0.3\textwidth]{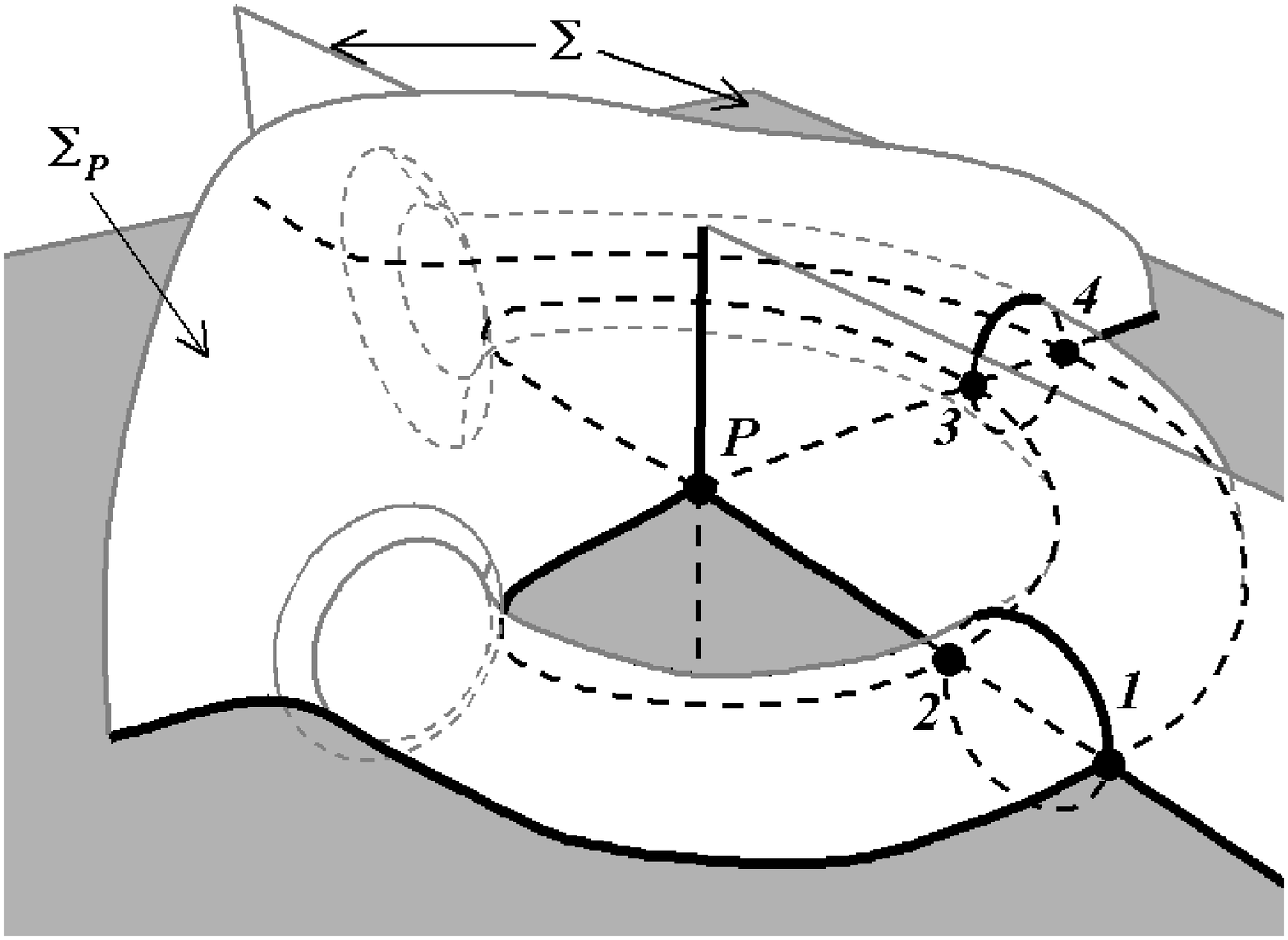}}\\
\subfigure []{ \label{fig8c}
\includegraphics[ width=0.3\textwidth]{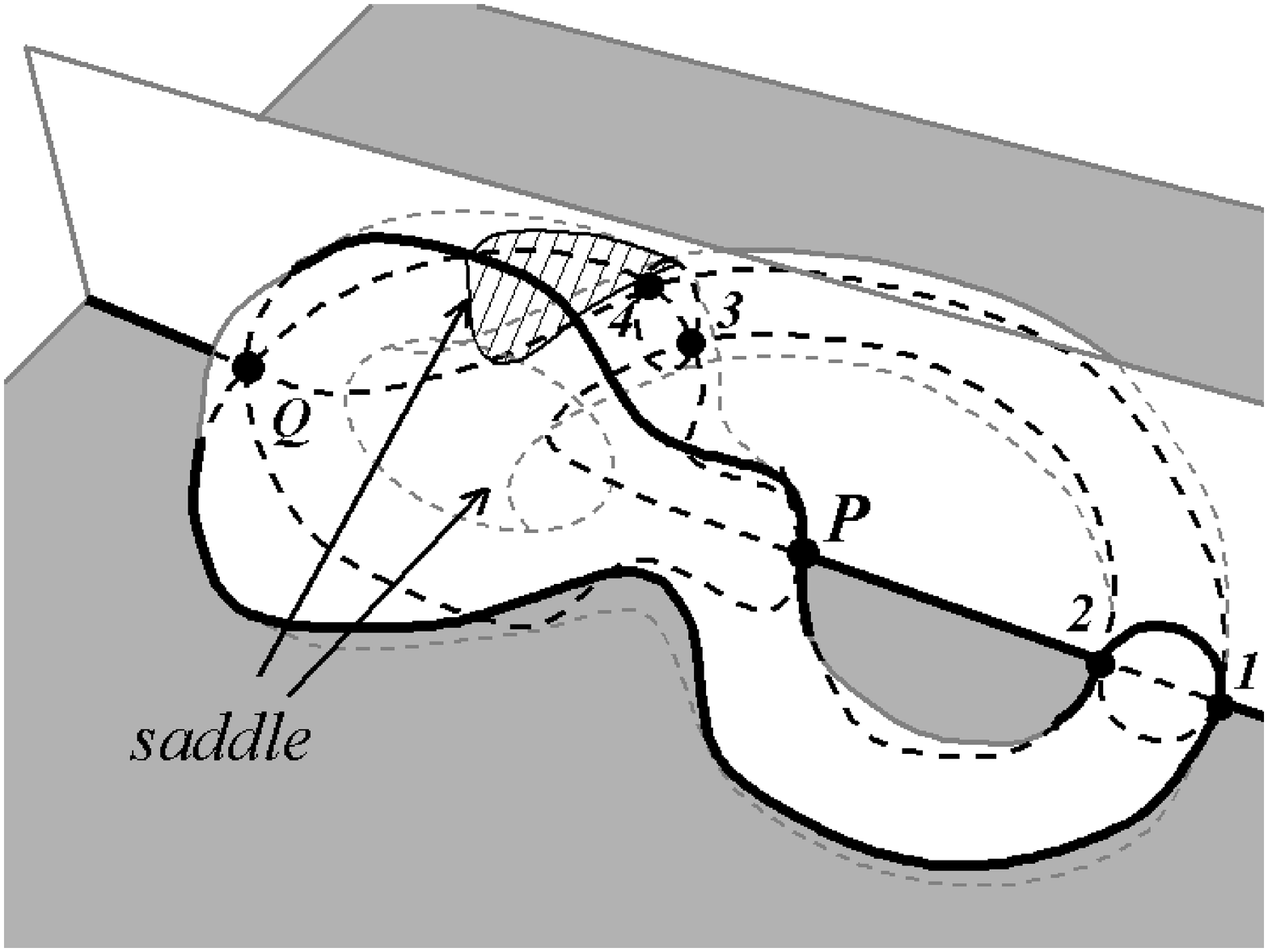}}\hfill
\subfigure[]{ \label{fig8d}
\includegraphics[ width=0.3\textwidth]{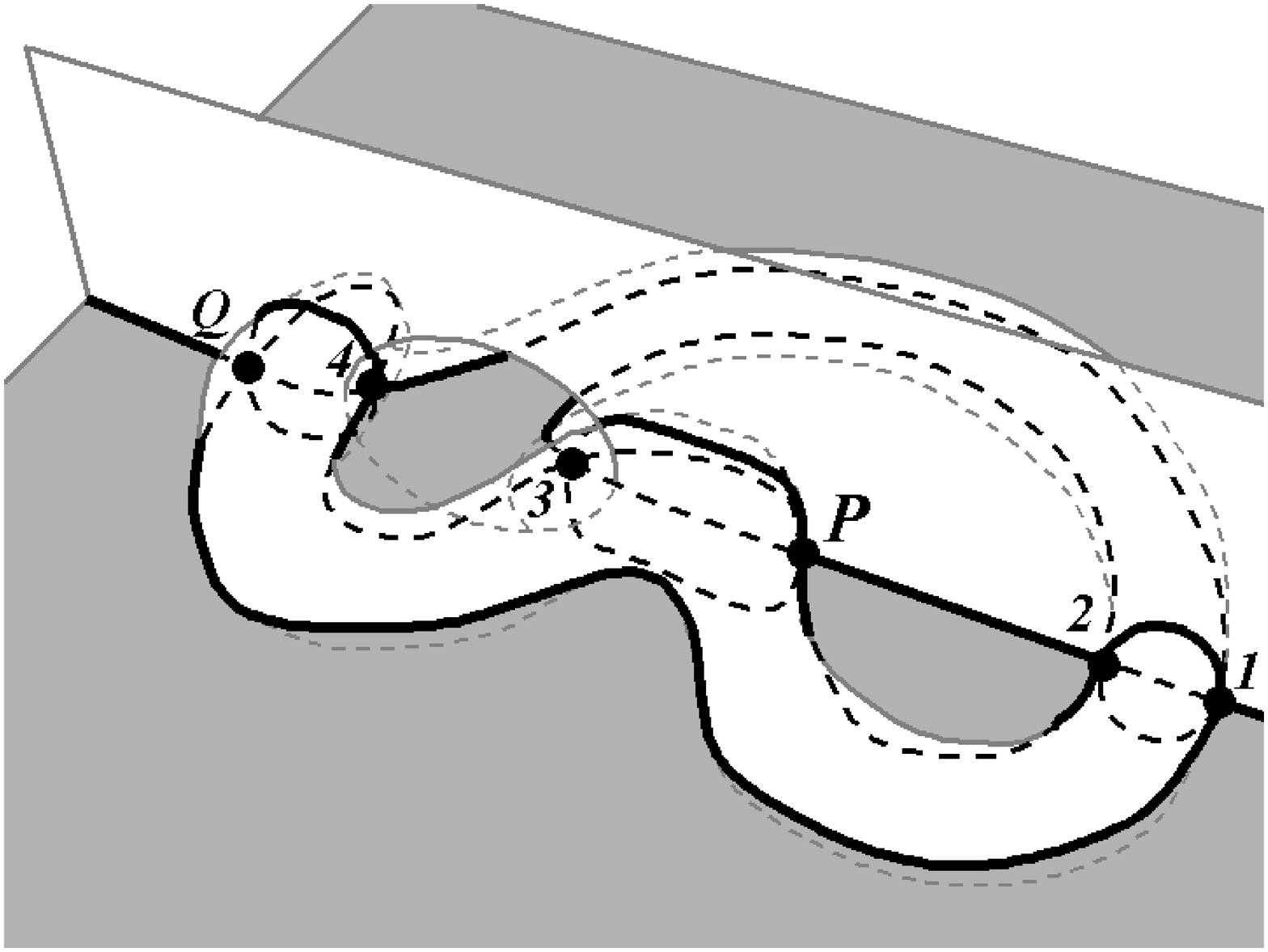}}\\
\subfigure []{ \label{fig8e}
\includegraphics[ width=0.3\textwidth]{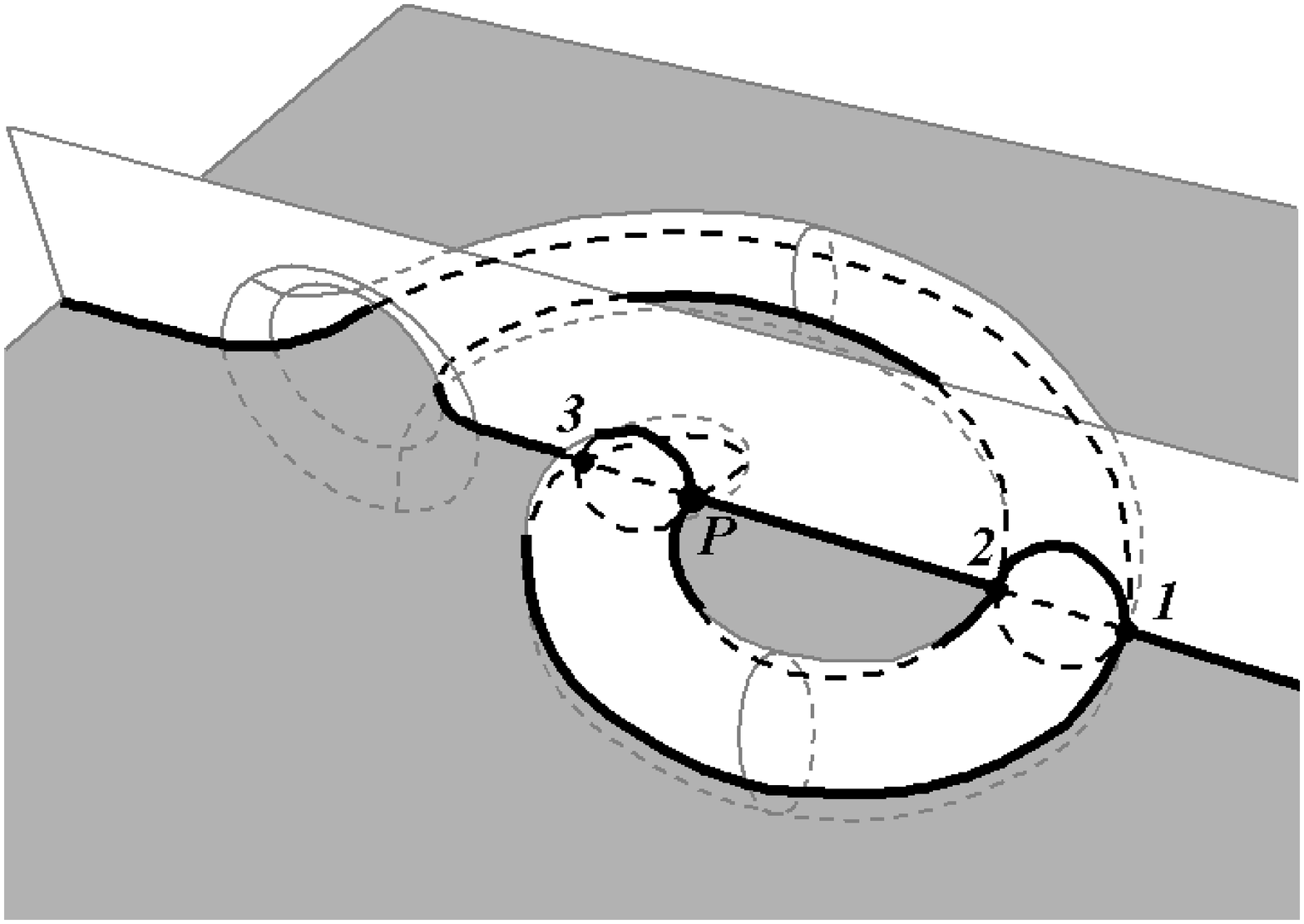}}\hfill
\subfigure[]{ \label{fig8f}
\includegraphics[ width=0.3\textwidth]{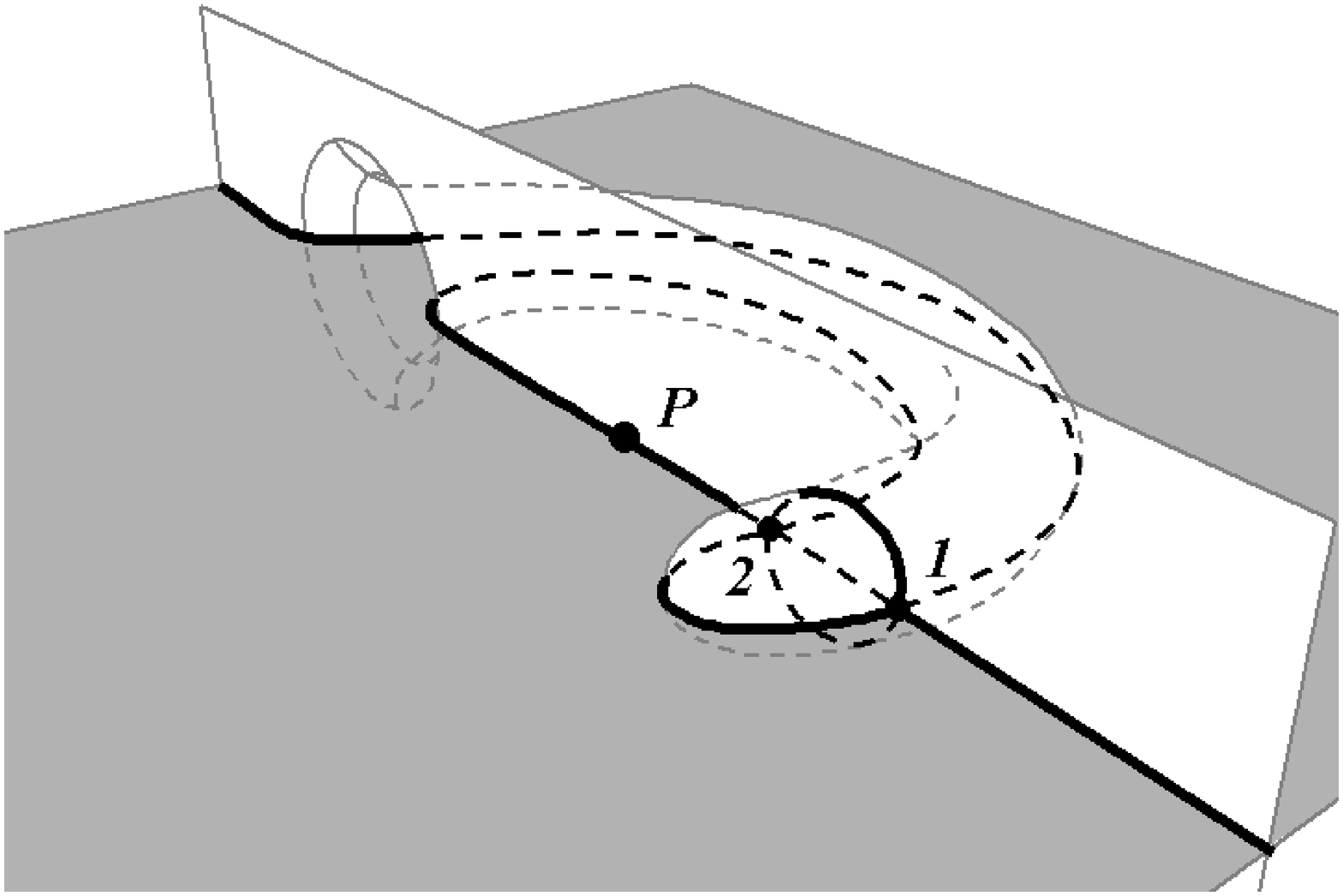}} \caption { }
\end{figure}

We can uninflate $\Sigma_{P}$ through $P$ by an ambient isotopy of
$M$, until we reach the situation depicted in Figure \ref{fig8c}.

We ''open the entrance of the tunnel'' using two
filling-preserving saddle moves, over and under the sheet of
$\Sigma$ containing the spiral piping (Figure \ref{fig8c}), and we
get the situation of Figure \ref{fig8d}. Now, after three
consecutive finger moves -1 we make $\Sigma_{P}$ desappear
completely (Figures \ref{fig8e} and \ref{fig8f}).
\end{proof}

In the previous Proposition, the statement should be ''every
parametrization of $\Sigma$ is filling homotopic to a
parametrization of $\Sigma\#\Sigma_{P} $'', which is a little
stronger that the chosen statement, but we use this language for
sake of simplicity.

We say that the Dehn surface $\Sigma\#\Sigma_{P}$ as in the previous
proposition is obtained from $\Sigma$ \textit{by inflating} $P$.

\subsection{Inflating 2-disks.\label{SUBSECTION InflatingDisks}}

\begin{figure}[htb]
\centering \subfigure []{ \label{fig9a}
\includegraphics[ width=0.35\textwidth]{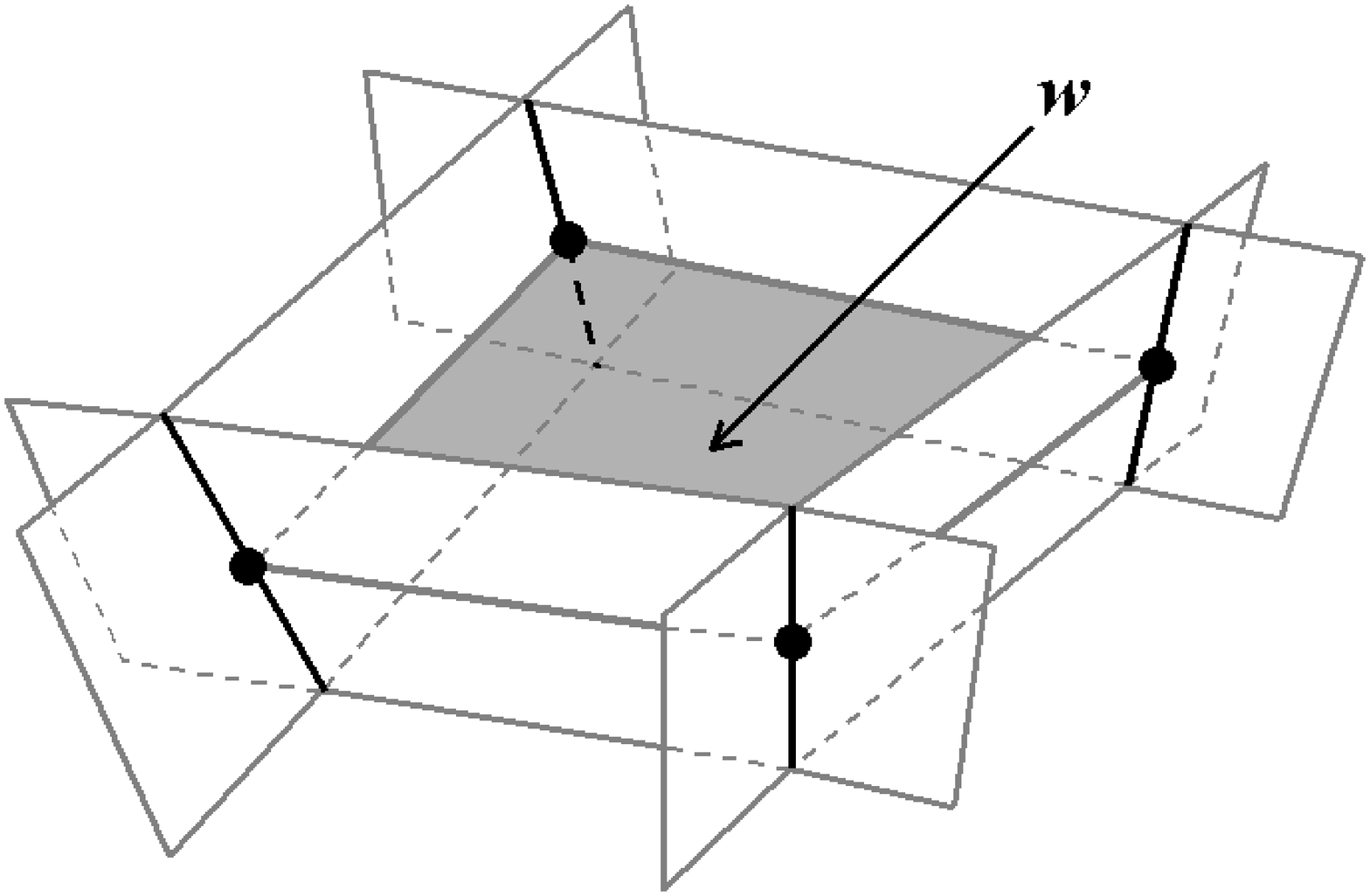}}\hfill \subfigure {
\includegraphics[ height=0.15\textheight, width=0.05\textwidth ] {fig35ab} }\hfill
\addtocounter {subfigure}{-1} \subfigure[]{ \label{fig9b}
\includegraphics[ width=0.35\textwidth]{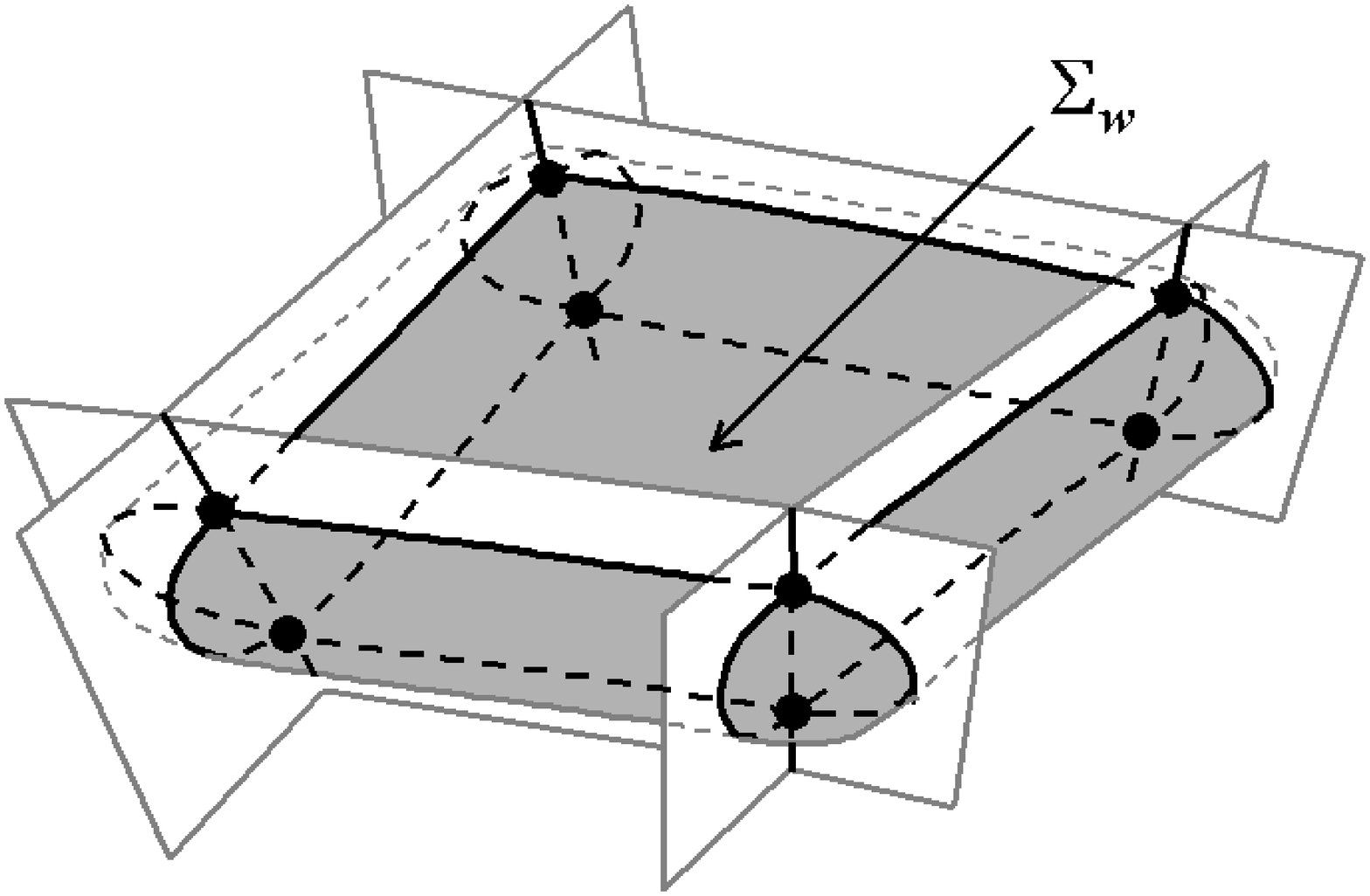}} \caption { }
\end{figure}

Let $R$ be a region of $M-\Sigma$. Consider a closed disk $w$ in
$M$ such that: (i) $(\partial w,int(w))\subset(\partial R,R)$;
(ii) $\partial w$ contains no triple point of $\Sigma$; and (iii)
$\partial w$ contains a finite number $n>0$ of double points of
$\Sigma$ ($n=4$ in Figure \ref{fig9a}). In this situation, we say
that $w$ is an $n$\textit{-gon} in $R$. We say that the double
points of $\Sigma$ in $\partial w$ are the \textit{vertices} of
$w$ and that the connected components of $\partial w-\left\{
\text{vertices of }w\right\}  $ are the \textit{edges} of $w$.

We can take a closed 3-ball $B_{w}\subset M$ with $w$ in its
interior as in Figure \ref{fig9b}. We take $B_{w}$ in such a way
that, if $\Sigma_{w}$ denote the 2-sphere bounding $B_{w}$, it
verifies that $\Sigma\cup\Sigma_{w}$ is a Dehn surface in $M$
(that is, $\Sigma_{w}$ intersects $\Sigma$ transversely). The
complement of $\Sigma$ in $B_{w}$ has $2n+1$ connected components.
The closure of one of the connected components of $B_{w}-\Sigma$
contains $w$ and it is a prism $w\times\left[  -1,1\right]  $ over
$w$ with $w\times\left\{  0\right\} \approx w$; its intersection
with $\Sigma$ is the ''boundary wall'' $\partial w\times\left[
-1,1\right]  $ and $w\times\left\{  -1\right\}  ,w\times \left\{
-1\right\}  $ are two $n$-gons in $R$ ''parallel'' to $w$. There
is a 2-gon prism (with one quadrangular face in $\partial
w\times\left[ -1,1\right]  $ and the other quadrangular face in
$\Sigma_{w}$) along each of the $n$ edges of $w$. Close to each
vertex of $w$ there appears a trihedron with two faces in $\Sigma$
and the other face in $S_{w}$. In this situation, we say that the
2-sphere $\Sigma_{w}$ is obtained by\textit{\ inflating }$w$. Note
that the new Dehn surface $\Sigma\cup\Sigma_{w}$ is another
filling Dehn surface of $M$.

Now, consider a point $P$ in $\Sigma_{w}$ which is a triple point of
$\Sigma\cup\Sigma_{w}$, and consider a filling Dehn surface $\Sigma
\#\Sigma_{w}$ that is obtained from $\Sigma\cup\Sigma_{w}$ by a spiral piping
around $P$. Then, we have the following Proposition:

\begin{proposition}
\label{PROPinflating2-cells}We can choose the piping such that $\Sigma$ is
filling homotopic to $\Sigma\#\Sigma_{w}$.
\end{proposition}

\begin{proof}
We start with the filling Dehn surface $\Sigma$. The point $P$ is
a double point of $\Sigma$. We inflate $P$ to obtain a 2-sphere
$\Sigma_{P}$ (piped with $\Sigma$), such that $\Sigma\#\Sigma_{P}$
is filling homotopic to $\Sigma$. We choose $\Sigma_{P}$ inside
$B_{w}$ (see Figure \ref{fig10a}). Now, we modify
$\Sigma\#\Sigma_{P}$ by $n-1$ consecutive finger moves +1 along
the edges of $w$ until we reach to the situation of Figure
\ref{fig10c}. We modify the Dehn surface of Figure \ref{fig10c} by
an ambient isotopy of $M$ until we reach the filling Dehn surface
$\Sigma^{\prime}$ of Figure \ref{fig10d}. Note that
$\Sigma^{\prime}$ differ from $\Sigma\#\Sigma_{w}$ by a saddle
move. This saddle move is filling preserving because both
$\Sigma^{\prime}$ and $\Sigma\#\Sigma_{w}$ are filling Dehn
surfaces of $M$. Thus, $\Sigma$ and $\Sigma\#\Sigma_{w}$ are
filling homotopic.

\begin{figure}[htb]
\centering \subfigure []{ \label{fig10a}
\includegraphics[ width=0.35\textwidth]{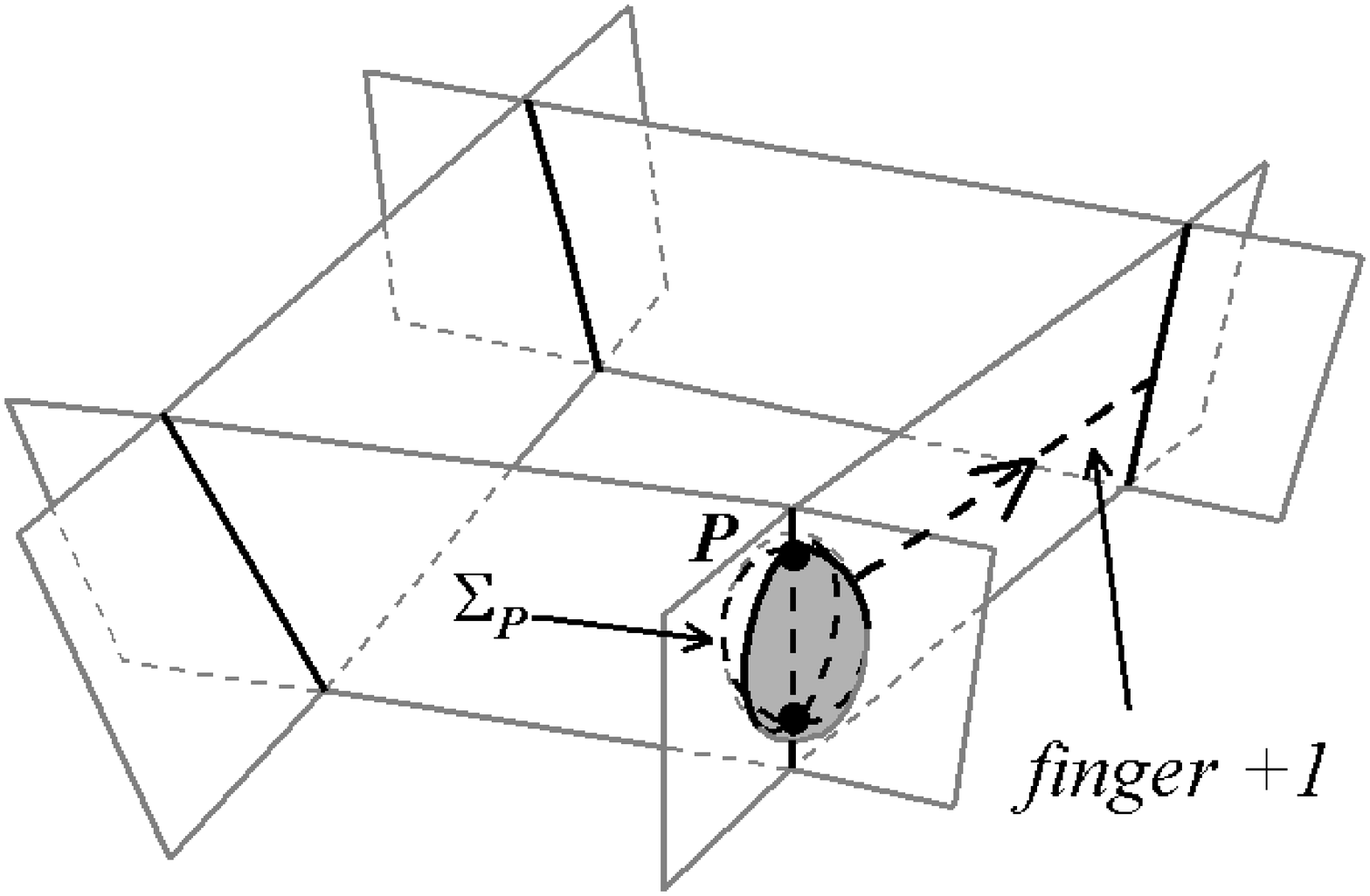}}\hfill
\subfigure[]{ \label{fig10b}
\includegraphics[ width=0.35\textwidth]{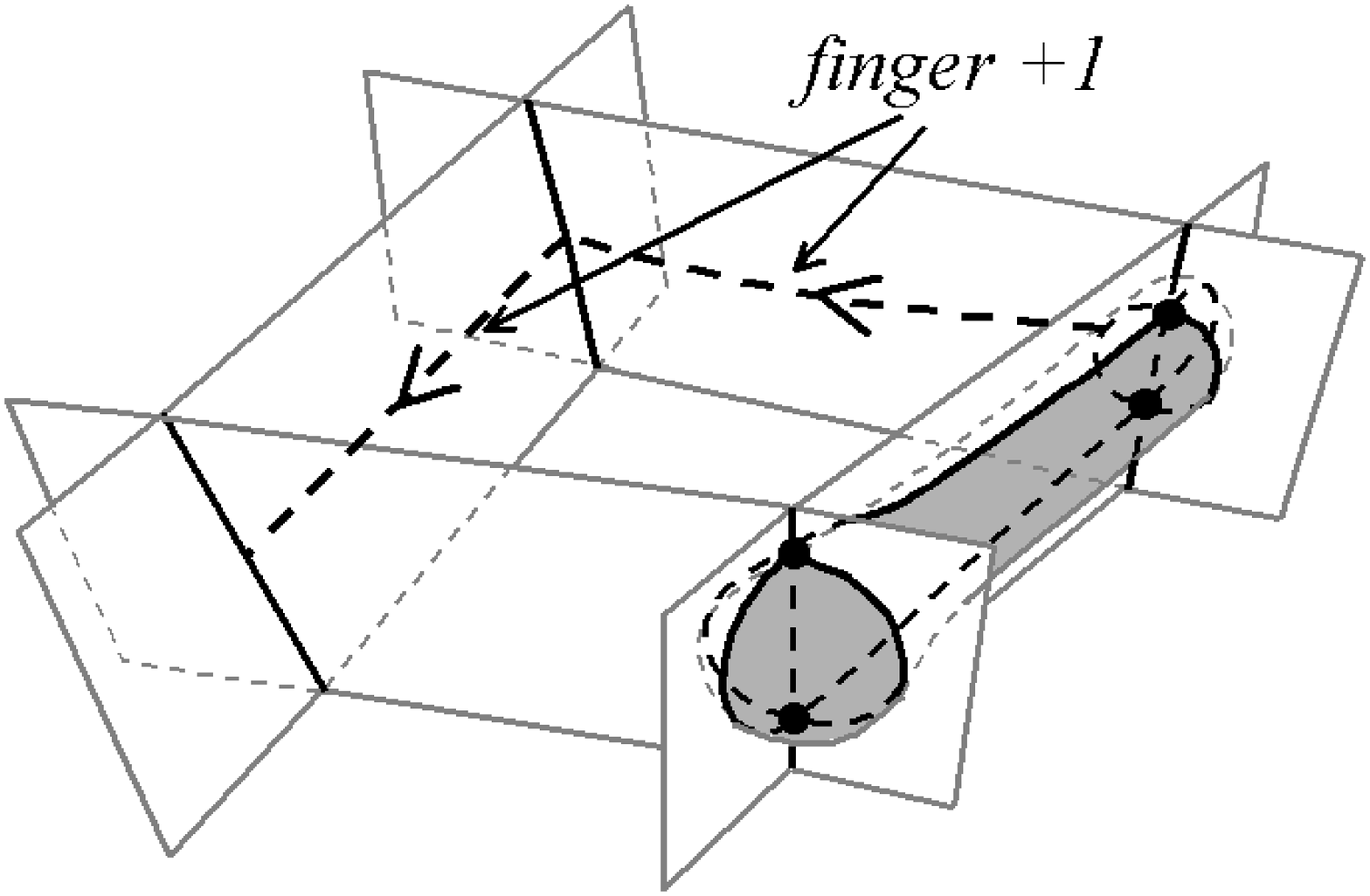}}\\
\subfigure []{ \label{fig10c}
\includegraphics[ width=0.35\textwidth]{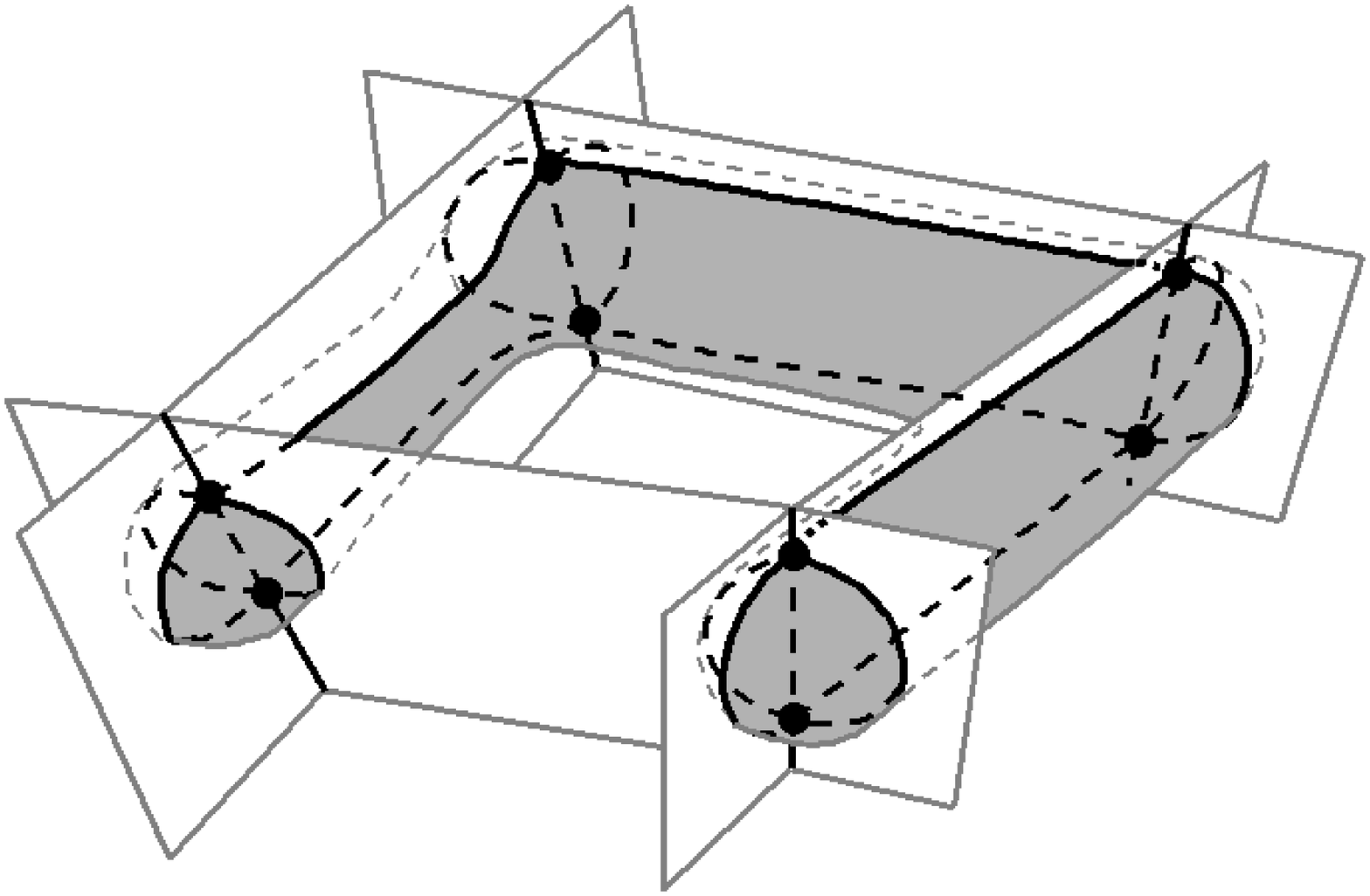}}\hfill
\subfigure[]{ \label{fig10d}
\includegraphics[ width=0.35\textwidth]{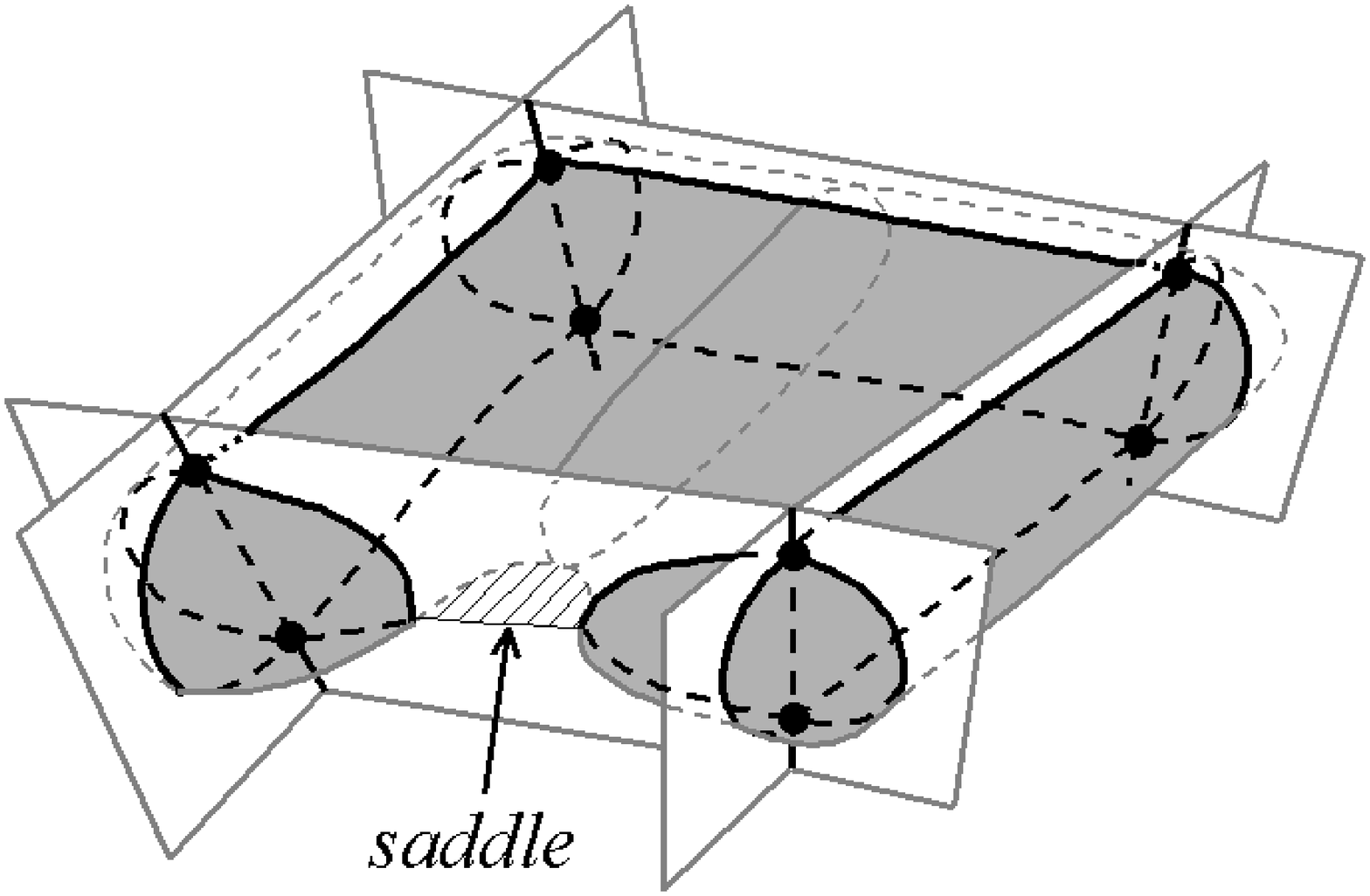}} \caption { }
\end{figure}

\end{proof}

The same operation as the one of inflating of $n$-gons can be done
if we start with a closed disk $w$ in $M$ contained in $\Sigma$.
Consider a closed disk $w$ in $M$ which is the closure of a
regular face of $\Sigma$. In this case we say that $w$ is an
$n$\textit{-gon in} $\Sigma$ if $n>0$ is the number of triple
points of $\Sigma$ lying in $\partial w$. We can inflate it in the
same way as an $n$-gon in a region of $M-\Sigma$, but in this case
we must introduce in Figures 9 and 10 an ''horizontal plane''
corresponding to the sheet of $\Sigma$ that contains $w$. We
inflate a 2-sphere $\Sigma_{P}$ upside the $n$-gon $w$. Then, we
use a finger move +1 to make it cross the $n$-gon $w$. Using
finger moves +2 instead of finger moves +1 and having a little
care in the final saddle move of the previous proof, it can be
seen also that we can connect the 2-sphere $\Sigma_{w}$ with
$\Sigma$ by a spiral piping such that the resulting filling Dehn
surface $\Sigma\#\Sigma_{w}$ is filling homotopic to $\Sigma$. The
final care to which we referred to above, consists on using two
saddle moves (one up and the other down the $n$-gon $w$) and a
final finger move of type -1 (we leave the details to the reader.

In both cases (the $n$-gon $w$ in a region of $M-\Sigma$ or in the Dehn sphere
$\Sigma$), we say that the Dehn surface $\Sigma\#\Sigma_{w}$ is obtained from
$\Sigma$ \textit{by inflating} $w$.

Thus, in the previous two paragraphs, we have seen that starting with a
filling Dehn surface $\Sigma$ and using only filling-preserving moves we can
obtain filling Dehn surface of $M$ more and more complicated, in analogy with
the finer and finer subdivisions of simplicial complexes or cell complexes.

\subsection{Passing over 3-cells.\label{SUBSECTION PassingCells}}

\begin{figure}[htb]
\centering
\includegraphics[ width=0.6\textwidth]{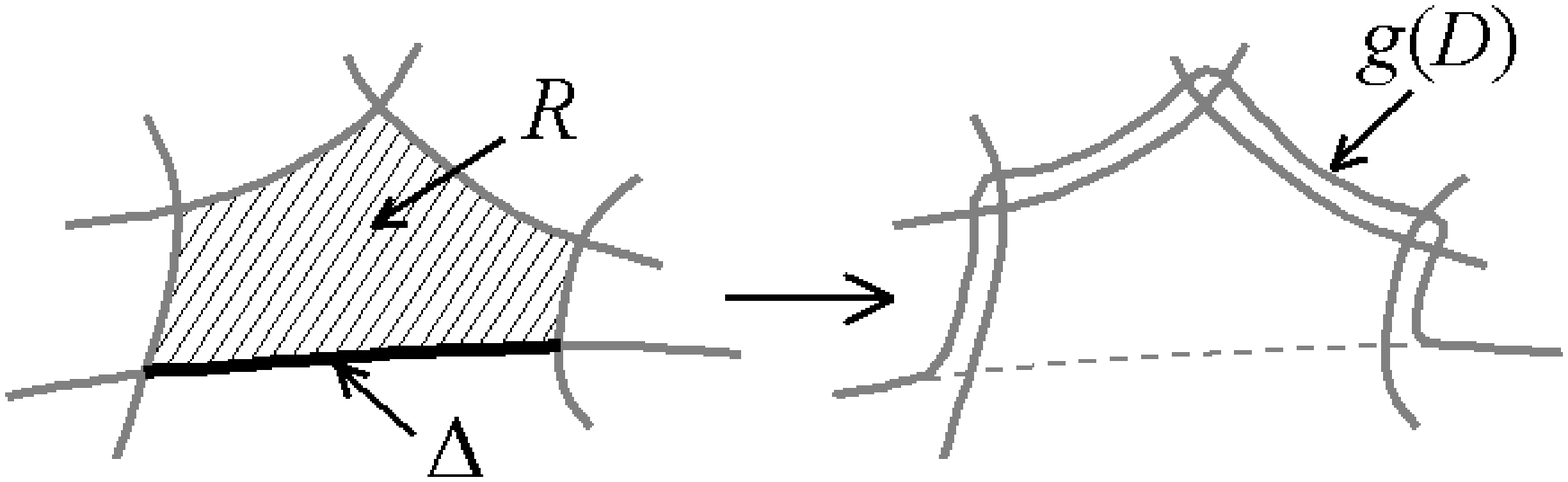} \caption { } \label{fig11}
\end{figure}

Let $\Delta$ be a regular face of $\Sigma$ and $R$ a regular
region of $M-\Sigma$ incident with $\Delta$. The regularity
implies that $cl(\Delta)$ is a closed disk and $cl(R)$ is a closed
3-ball. Take a parametrization $f:S\rightarrow M$ of $\Sigma$, and
take the immersion $g:S\rightarrow M\;$that is obtained from $f$
by a pushing disk $(D,B)$ as it is indicated in Figure
\ref{fig11}. There is an open disk $\tilde{\Delta}$ in $S$ such
that the restriction of $f$ to $\tilde{\Delta}$ is an embedding
and it is $f(\tilde{\Delta})=\Delta$. The pushed disk $D$ in $S$
contains $cl(\tilde {\Delta})$ in its interior and it is as close
to $\tilde{\Delta}$ as necessary, such that $f\mid_{D}$ is an
embedding. The pushing ball $B$ contains the region $R$ and it is
$cl(R)\cap\partial B=cl(\Delta)$. The disk $g(D)$ is a closed disk
outside $cl(R)$ running in parallel to $\partial R-\Delta$. In
\cite{RHomotopies} it is proved the following:

\begin{lemma}
\label{LemmaPassing3-cells}If $g$ is a filling immersion, then it is filling
homotopic to $f$.
\end{lemma}

In paper \cite{RHomotopies} a study is made to ascertain which properties of
the pair $\Delta,R$ are necessary and sufficient for $g$ to be a filling immersion.

\section{Shellability. Smooth triangulations.\label{SECTION Shellability}}

In the proof of Theorem \ref{MAINtheorem} we use an exhaustive use of the
concept of \textit{shelling}.

\begin{definition}
\label{DEFfreen-cell}Let $N$ be an $n$-manifold with boundary, and let
$C\subset N$ be a closed $n$-ball in $N$. We say that $C$ is free in $N$ if
$C\cap\partial N$ is a closed $(n-1)$-ball.
\end{definition}

Let $B$ be a closed $n$-ball, and let $K$ be a regular cell decomposition of
$B$.

\begin{definition}
\label{DEFShellable}$K$ is shellable if there exists an ordering
$C_{1} ,C_{2},...,C_{k}$ of its $n$-cells such that $cl(C_{i})$ is
free in the closure of $(\underset{j\geq i}{\cup}C_{j})$. If this
is the case, we say that the ordering $C_{1},C_{2},...,C_{k}$ is a
shelling of $K$.
\end{definition}

While cell decompositions of 2-disks are always shellable (see Lemma 1 of
\cite{Sanderson}), non-shellable cell decompositions of $n$-balls exist for
$n>2$. In \cite{Bruggeser-Mani} it is proved that every cell decomposition of
an $n$-ball has a shellable subdivision (a cell decomposition $\sigma K$ of
$B$ is a \textit{subdivision} of the cell decomposition $K$ of $B$ if every
cell of $\sigma K$ is contained in a cell of $K$.

For the proof of Theorem \ref{MAINtheorem}, we will work with triangulations
of a 3-manifold $M$ and we will require that a set of non necessary disjoint
3-balls in $M$ (endowed with the induced cell decomposition) are all shellable
at once, and with some special properties. For this, we will use the work of
Whitehead about simplicial collapsings (see \cite{Whitehead} or \cite{Glaser}).

A simplicial complex whose underlying polyhedron is a ball induces in a
natural way a cell decomposition of the ball. Thus, we can consider simplicial
complexes also as particular cell complexes.

If $K$ is a simplicial complex, a simplex of $K$ is \textit{maximal} if it is
not a proper face of another simplex of $K$. If $\epsilon^{i}$ is a maximal
$i$-simplex of $K$, an $(i-1)$-face $\epsilon^{i-1}$ of $\epsilon^{i}$ is
\textit{free} in $K$ if it is not a face of another $i$-simplex of $K$
different from $\epsilon^{i}$. If $\epsilon^{i}$ is maximal in $K$ and
$\epsilon^{i-1}$ is free in $K$, then the result of removing $(\epsilon
^{i},\epsilon^{i-1})$ from $K$ is another simplicial complex $K^{\prime}$ and
it is said that $K^{\prime}$ is obtained from $K$ by a \textit{simplicial
collapsing}. The complex $K$ \textit{collapses simplicially} into a subcomplex
$K^{\prime}$ if $K^{\prime}$ is obtained from $K$ after a finite sequence of
simplicial collapsings. In particular, $K$ is \textit{collapsible} if it
collapses simplicially into a point. If $\sigma K$ is a subdivision of $K$ and
$K^{\prime}$ is a subcomplex of $K$, then $\sigma K^{\prime}$ will denote the
corresponding subcomplex of $\sigma K$.

\begin{theorem}
\cite{Whitehead}\label{THWhiteheadsubdivisions}If $K$ is any finite simplicial
complex, there is a (stellar) subdivision $\sigma K$ of $K$ such that $\sigma
B^{n}$ collapses simplicially into $\sigma B^{n-1}$, where $B^{n}$ is any
subcomplex of $K$ which is a closed $n$-ball and $B^{n-1}$ is any subcomplex
of $\partial B^{n}$ which is a closed $(n-1)$-ball.
\end{theorem}

For triangulations of an $n$-ball, shellability obviously implies
collapsibility. Note that the converse is not obvious because in shellability
we require that the space after each step remains a ball, while in
''collapsings'' it might not even be a manifold. (The example given in
\cite{Rudin} is not shellable, but it is simplicially collapsible, compare
\cite{Chillingworth}). However,\ the converse is almost true in our case
according to the following observation that arises from the proof of Theorem 6
in \cite{Bing}:

\begin{lemma}
\label{LEMABing}If $K$ is a collapsible triangulation of the 3-ball, then the
second derived subdivision of $K$ is shellable.
\end{lemma}

Smooth triangulations of manifolds are introduced in \cite{Whitehead2} for
relating the smooth and PL categories in manifold theory. A triangulation of
an $n$-manifold $N$ is an homeomorphism $h:K\rightarrow N$, where $K$ is a
rectilinear simplicial complex of some euclidean space. If $N$ has a smooth
structure, the triangulation $h$ is \textit{smooth} (with respect to this
structure) if the restriction of $h$ to each simplex of $K$ is a smooth map.
We identify the manifold $N$ with the simplicial complex $K$. In
\cite{Whitehead2} (see also \cite{Munkres}) it is proved that: (i) any
$n$-manifold with a smooth structure admits smooth triangulations; and (ii)
two smooth triangulations of the same smooth manifold have a common smooth
subdivision (!). If $f:S\rightarrow M$ is a transverse immersion from the
surface into the 3-manifold $M$, then there are smooth triangulations $K$ and
$T$ of $S$ and $M$ respectively such that $f$ is simplicial with respect to
them (for more general results of this kind, see \cite{Verona}).

If $f:S\rightarrow M$ is a filling immersion and $K,T$ are triangulations of
$S,M$ respectively such that $f$ is simplicial with respect to them, then the
triangulation $T$ triangulates also the closure of each region of $M-f(S) $.
If $R$ is a regular region of $M-f(S)$, we say that $T$ \textit{shells} $R $
if it induces a shellable triangulation on $cl(R)$. If $R$ is not regular, we
cut first $cl(R)$ along its self-adjacencies to obtain a closed 3-ball
$\widetilde{cl(R)}$. The triangulation $T$ on $cl(R)$ lifts naturally to
$\widetilde{cl(R)}$, and so we say that $T$ \textit{shells} $R$ if the induced
triangulation on $\widetilde{cl(R)}$ is shellable. The triangulation $T$
\textit{shells} the filling immersion $f$ (or the filling Dehn surface $f(S)$)
if $T$ shells each region of $M-f(S)$.

All these results imply the following

\begin{theorem}
\label{THMsmoothtriangulations exist}Let $S_{1},...,S_{k}$ be a
finite collection of surfaces, and for each $i=1,...,k$ let
$f_{i}:S_{i}\rightarrow M $ be a transverse immersion. Then, there
is a smooth triangulation $T$ of $M$ such that:

\begin{enumerate}
\item  the Dehn surfaces $f_{1}(S_{1}),...,f_{k}(S_{k})$ are contained in the
2-skeleton of $T$;

\item  if $f_{i}$ is a filling immersion for some $i=1,...,k$, the
triangulation $T$ shells $f_{i}$;

\item  if $f_{i}$ and $f_{j}$ differ by a pushing disk $(D,B)$ for some
$i,j\in\left\{  1,...,k\right\}  $, then the triangulation $T$
restricted to $B $ collapses simplicially into $f_{j}(D)$.
\end{enumerate}
\end{theorem}

\begin{proof}
First of all, we have seen that for each $i=1,...,k$ there are
smooth triangulations $K_{i},T_{i}$ such that $f_{i}$ is
simplicial with respect to them. According to \cite{Whitehead2},
the smooth triangulations $T_{1} ,...,T_{k}$ have a common smooth
subdivision $T_{0}$. Then, all the Dehn surfaces
$f_{1}(S_{1}),...,f_{k}(S_{k})$ are contained in the 2-skeleton of
$T_{0}$. Take a subdivision $T_{0}^{\prime}$ of $T_{0}$ in the
conditions of Theorem \ref{THWhiteheadsubdivisions}, and take $T$
as the second derived subdivision of $T_{0}^{\prime}$.

If $f_{i}$ and $f_{j}$ differ by a pushing disk $(D,B)$, then $T_{0}$
triangulates $B$ and $f_{j}(D)$, and so because $T_{0}^{\prime}$ has been
chosen following Theorem \ref{THWhiteheadsubdivisions}, the triangulation
$T_{0}^{\prime}$ restricted to $B$ collapses simplicially into $f_{j}(D)$.
Simplicial collapsing is preserved by stellar subdivisions \cite{Whitehead}
and so it is also preserved by derived subdivisions. Thus, $T$ restricted to
$B$ collapses simplicially into $f_{j}(D)$.

If $f_{i}$ is a filling immersion and $R$ is a regular region of
$M-f_{i}(S_{i})$, then $T_{0}$ triangulates $cl(R)$. By the choosing of
$T_{0}^{\prime}$, $T_{0}^{\prime}$ induces a collapsible triangulation of
$cl(R)$, and by Lemma \ref{LEMABing}, $T$ induces a shellable triangulation on
$cl(R)$. If $R$ is not regular, perhaps we need to do more stellar
subdivisions on $cl(R)$ to have the required shelling property on
$\widetilde{cl(R)}$, but this does not alter the previous construction because
stellar subdivisions preserve shellability \cite{Bruggeser-Mani} and collapsibility.
\end{proof}

\begin{definition}
\label{DEF GoodTriangulations}In the hypothesis of the previous Theorem, we
say that $T$ is a good triangulation of $M$ with respect to $f_{1},...,f_{k}$.
\end{definition}

With these results, we have prepared the ground for the following sections.

\section{Inflating triangulations.\label{SECTION InflatingTriangulations}}

Now we will explain how we can associate to any triangulation of the
3-manifold $M$ a filling collection of spheres of $M$.

\begin{figure}[htb]
\centering \subfigure []{ \label{fig12a}
\includegraphics[ width=0.3\textwidth]{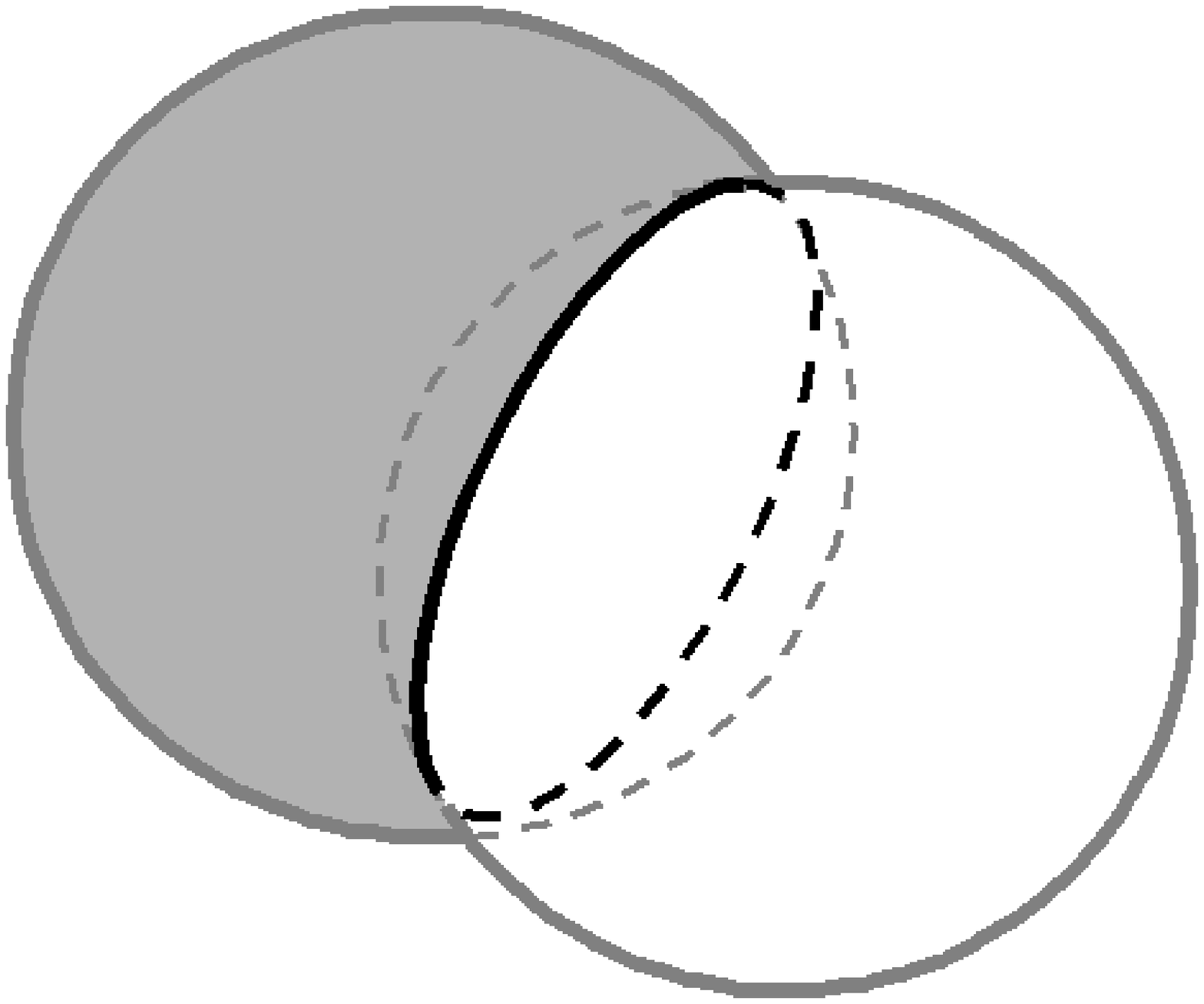}}\hspace{2cm}
\subfigure[]{ \label{fig12b}
\includegraphics[ width=0.3\textwidth]{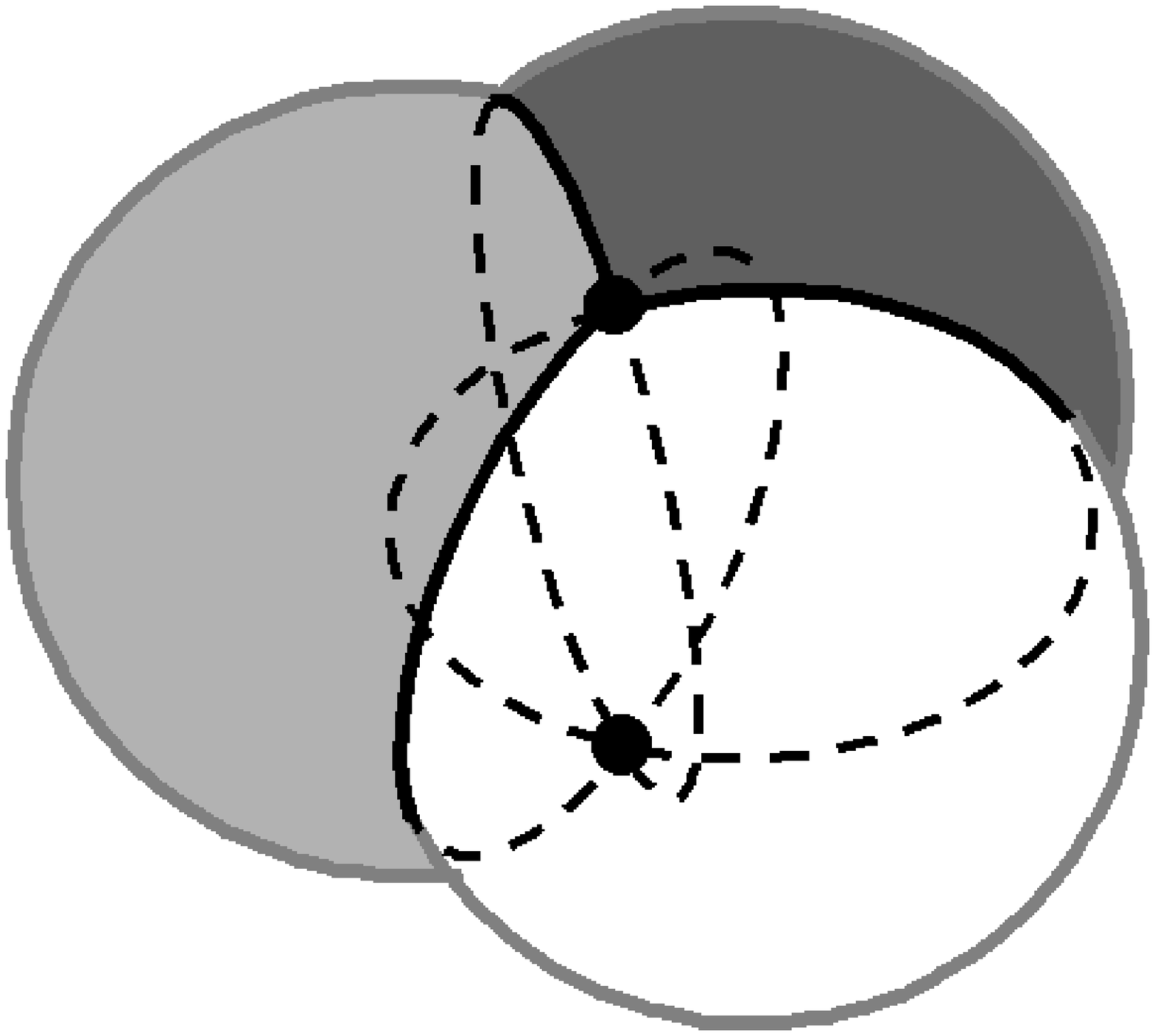}} \caption { }
\end{figure}

Let $B_{1},B_{2}$ be two closed 3-balls in $M$. We say that
$B_{1},B_{2}$ \textit{intersect normally} if they intersect as in
Figure \ref{fig12a}. The 2-spheres $\partial B_{1},\partial B_{2}$
must intersect transversely in a unique simple closed curve. If
$B_{1},B_{2}$ intersect normally, then $B_{1}\cap B_{2}$,
$cl(B_{1}-B_{2})$ and $cl(B_{2}-B_{1})$ are 3-balls. If
$B_{1},B_{2},B_{3}$ are 3-balls in $M$, they \textit{intersect
normally} if they intersect as in Figure \ref{fig12b}. Each pair
$B_{i},B_{j}$ with $i\neq j$ intersect normally and $\partial
B_{1},\partial B_{2}$ and $\partial B_{3}$ intersect transversely
at two triple points.

Let $T$ be a smooth triangulation of $M$. (We refer to the
0-simplexes, 1-simplexes, 2-simplexes and 3-simplexes of $T$ as
vertices, edges, triangles and tetrahedra of $M$, respectively.)
We can construct a filling collection of spheres of $M$ by
''inflating'' $T$ assigning to each simplex $\epsilon$ of the
2-skeleton $T^{2}$ of $T$ a 2-sphere $S\epsilon$ embedded in $M$
in such a way that their union $\mathbb{T}=\underset{\epsilon\in
T^{2}}{\cup}S\epsilon$ fills $M$. We will do this as follows.

\begin{figure}[htb]
\centering \subfigure []{ \label{fig13a}
\includegraphics[ width=0.4\textwidth]{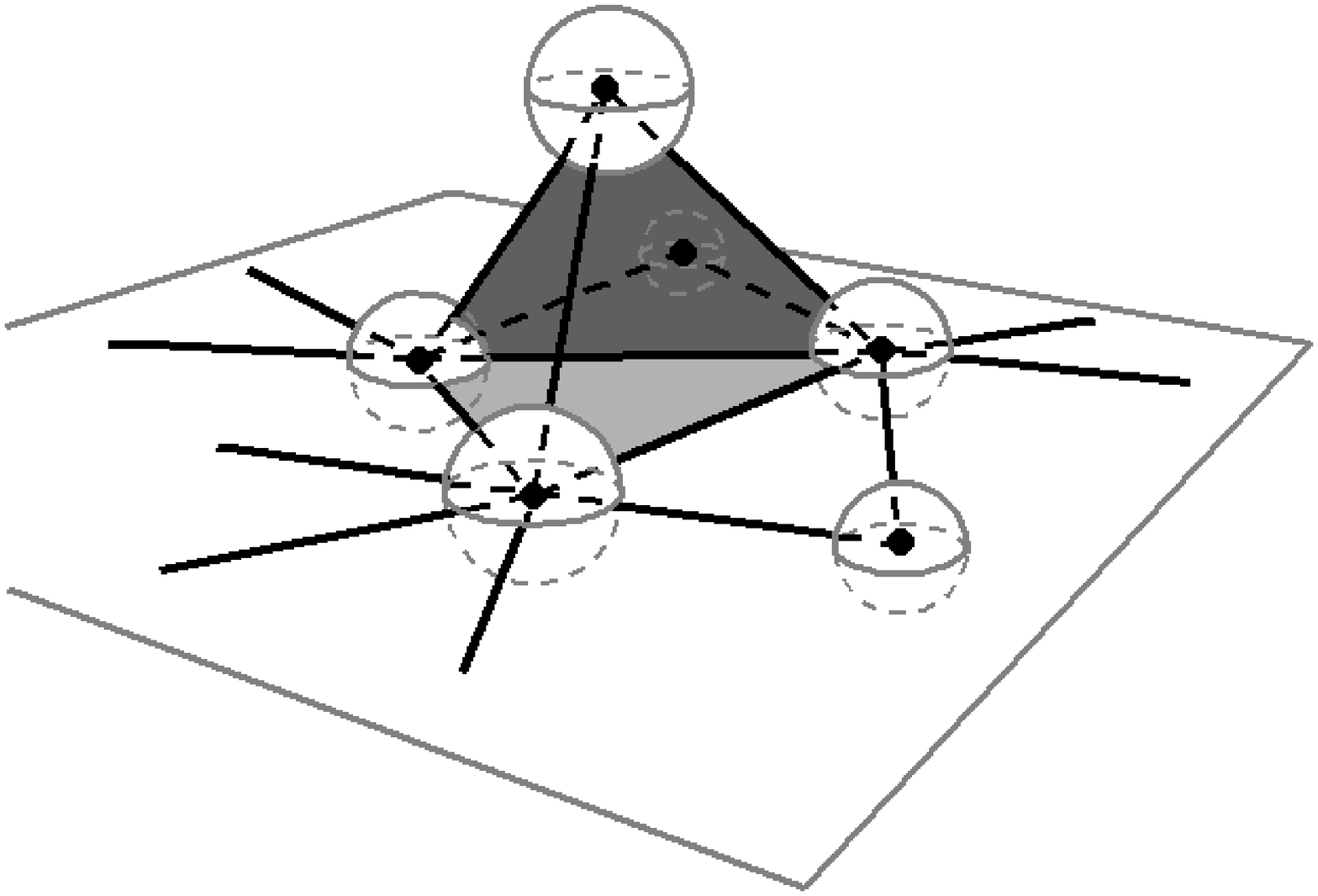}}\hfill \subfigure[]{ \label{fig13b}
\includegraphics[ width=0.4\textwidth]{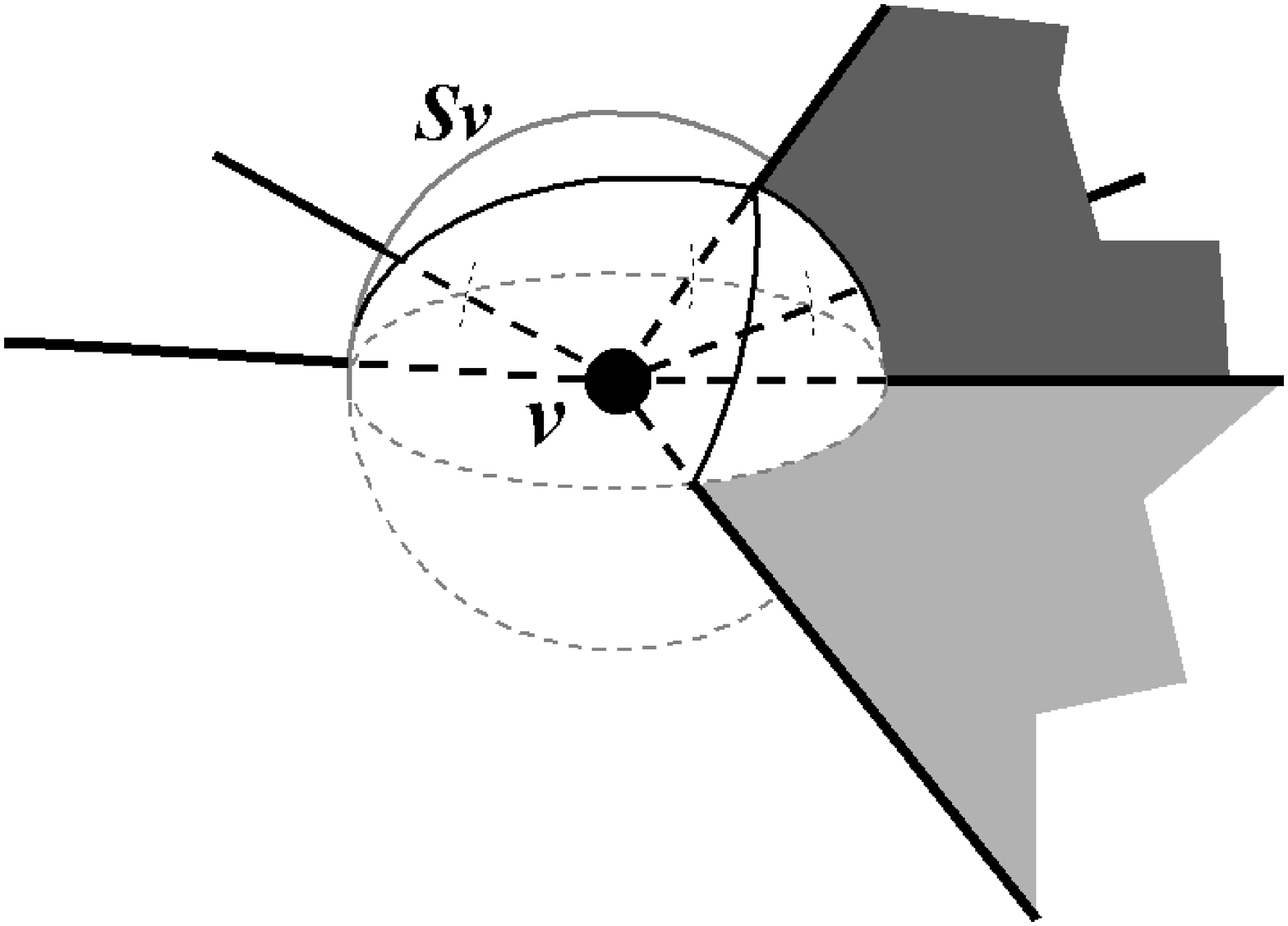}} \caption { }
\end{figure}

First, if $v_{1},...,v_{m_{0}}$ are the vertices of $M$, for
$i=1,...,m_{0}$ the 2-sphere $Sv_{i}\subset M$ bounds a closed
3-ball $Bv_{i}$ in $M$ contained in the open star $star\left(
v_{i}\right)  $ and with $v_{i}$ in its interior. The 2-spheres
$Sv_{1},...,Sv_{m_{0}}$ are pairwise disjoint (Figure
\ref{fig13a}) and the triangulation $T$ of $M$ induces a
triangulation of $Bv_{i}$ as a \textit{cone} from $v_{i}$ over
$Sv_{i}$ for each $i=1,...,m_{0} $ (see Figure \ref{fig13b}). The
2-sphere $Sv_{i}$ intersects transversely each $i$-simplex
$\epsilon^{i}\in star\left(  v_{i}\right)  \subset T$ in a $\left(
i-1\right) $-simplex of this induced triangulation of $Sv_{i}$.

\begin{figure}[htb]
\centering \subfigure []{ \label{fig14a}
\includegraphics[ width=0.4\textwidth]{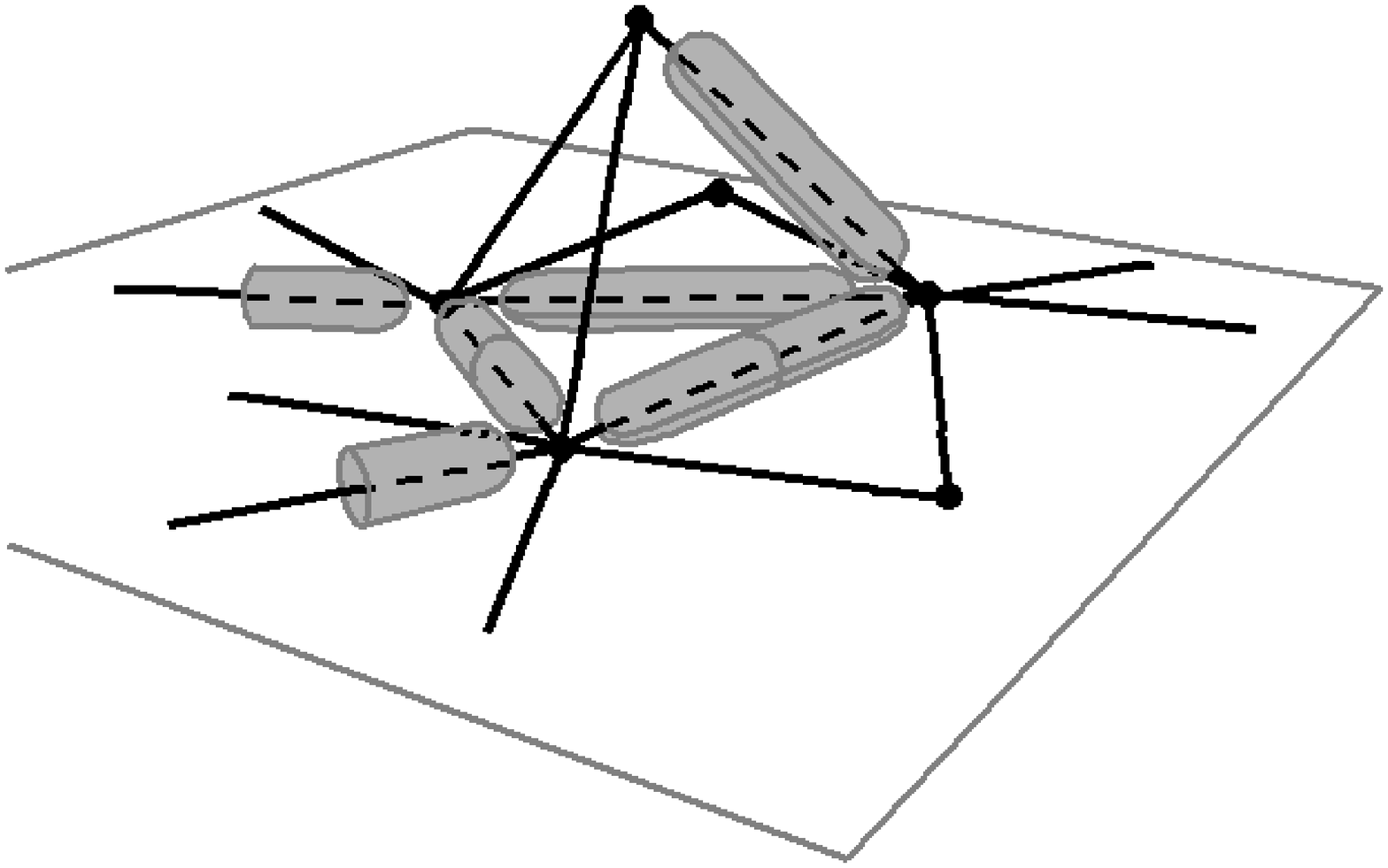}}\hfill \subfigure[]{ \label{fig14b}
\includegraphics[ width=0.4\textwidth]{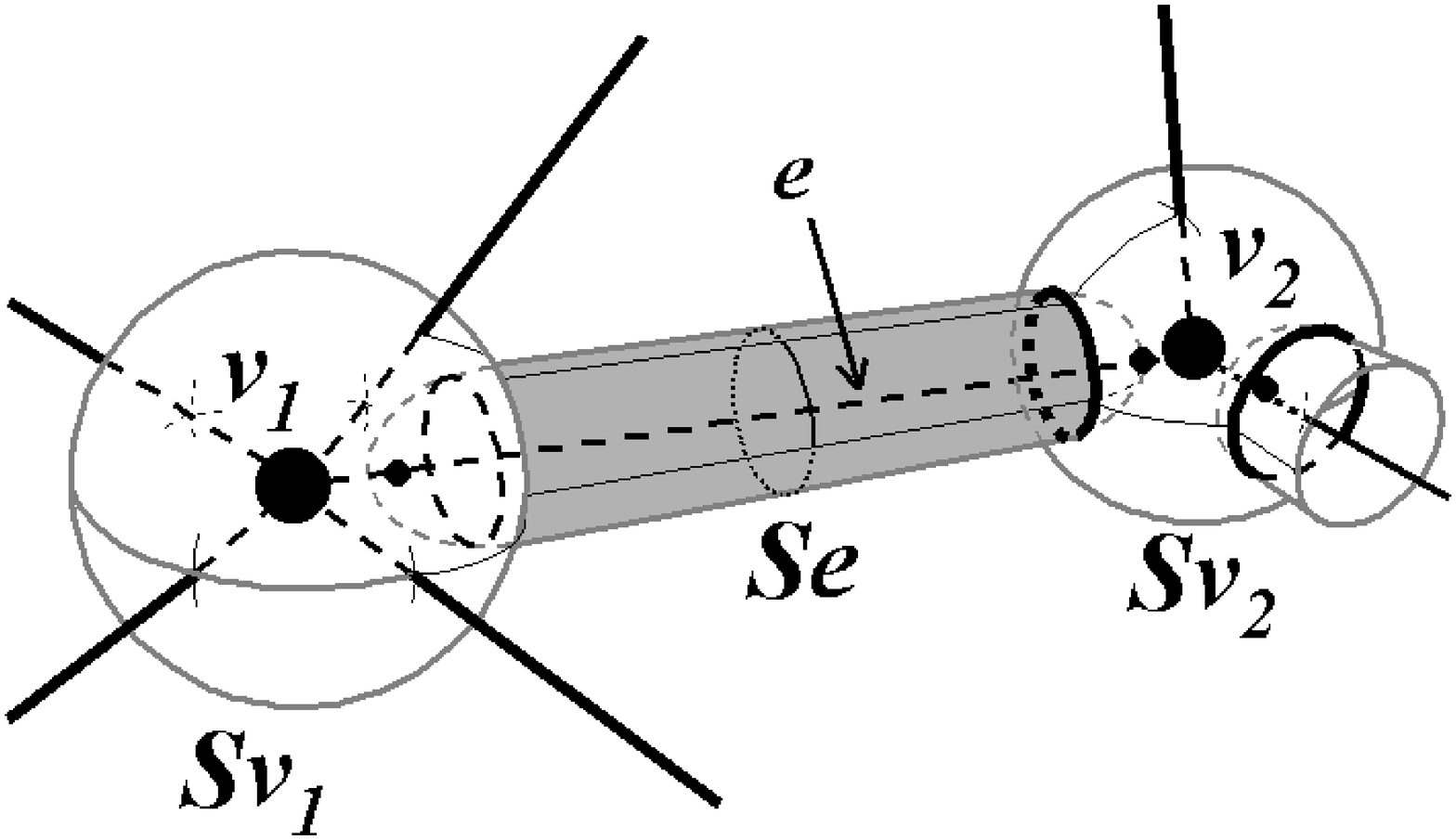}} \caption { }
\end{figure}

If $e_{1},...,e_{m_{1}}$ are the edges of $M$, for $j=1,...,m_{1}$
the 2-sphere $Se_{j}\subset M$ bounds a closed 3-ball $Be_{j}$ in
$M$ as in Figure \ref{fig14a}. The 3-ball $Be_{j}$ is contained in
the open star $star\left( e_{j}\right)  $ and it intersects
$e_{j}$ in a closed sub-arc $\tilde{e} _{j}\subset e_{j}$. The
2-sphere $Se_{j}$ and $e_{j}$ intersect transversely at the
endpoints of the arc $\tilde{e}_{j}$ We take
$Be_{1},...,Be_{m_{1}}$ pairwise disjoint, and for each
$i\in\left\{  1,...,m_{0}\right\}  $ and $j\in\left\{
1,...,m_{1}\right\}  $ the 3-balls $Bv_{i}$ and $Be_{j}$ are also
disjoint unless $v_{i}$ and $e_{j}$ are incident. In this case,
$Bv_{i}$ and $Be_{j}$ intersect normally (Figure \ref{fig14b}) and
$Bv_{i}\cap Be_{j}$ intersects $e_{j}$ in another closed sub-arc
of $e_{j}$. Considering the two points of the intersection
$Se_{j}\cap e_{j}$ as the ''poles'' of $Se_{j}$, each triangle $t$
of $M$ incident with $e_{j}$ intersects $Se_{j}$ transversely in
an open arc which is the interior of a ''meridian'' $a$ with its
endpoints at the poles. The intersection $cl(t)\cap Be_{j}$ is a
closed disk bounded by $a\cup\tilde{e}_{j}$.

\begin{figure}[htb]
\centering \subfigure []{ \label{fig15a}
\includegraphics[ width=0.4\textwidth]{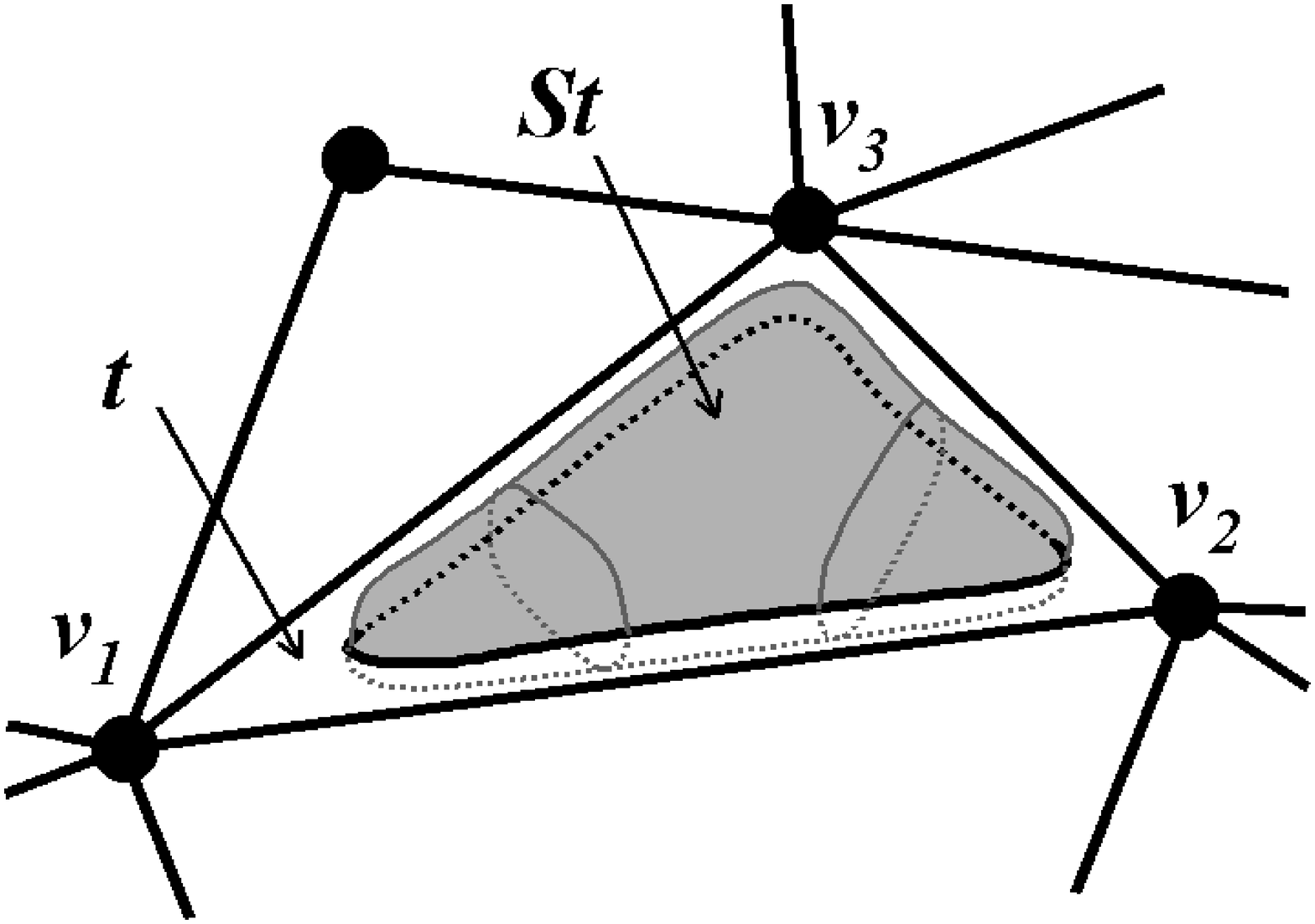}}\hfill \subfigure[]{ \label{fig15b}
\includegraphics[ width=0.4\textwidth]{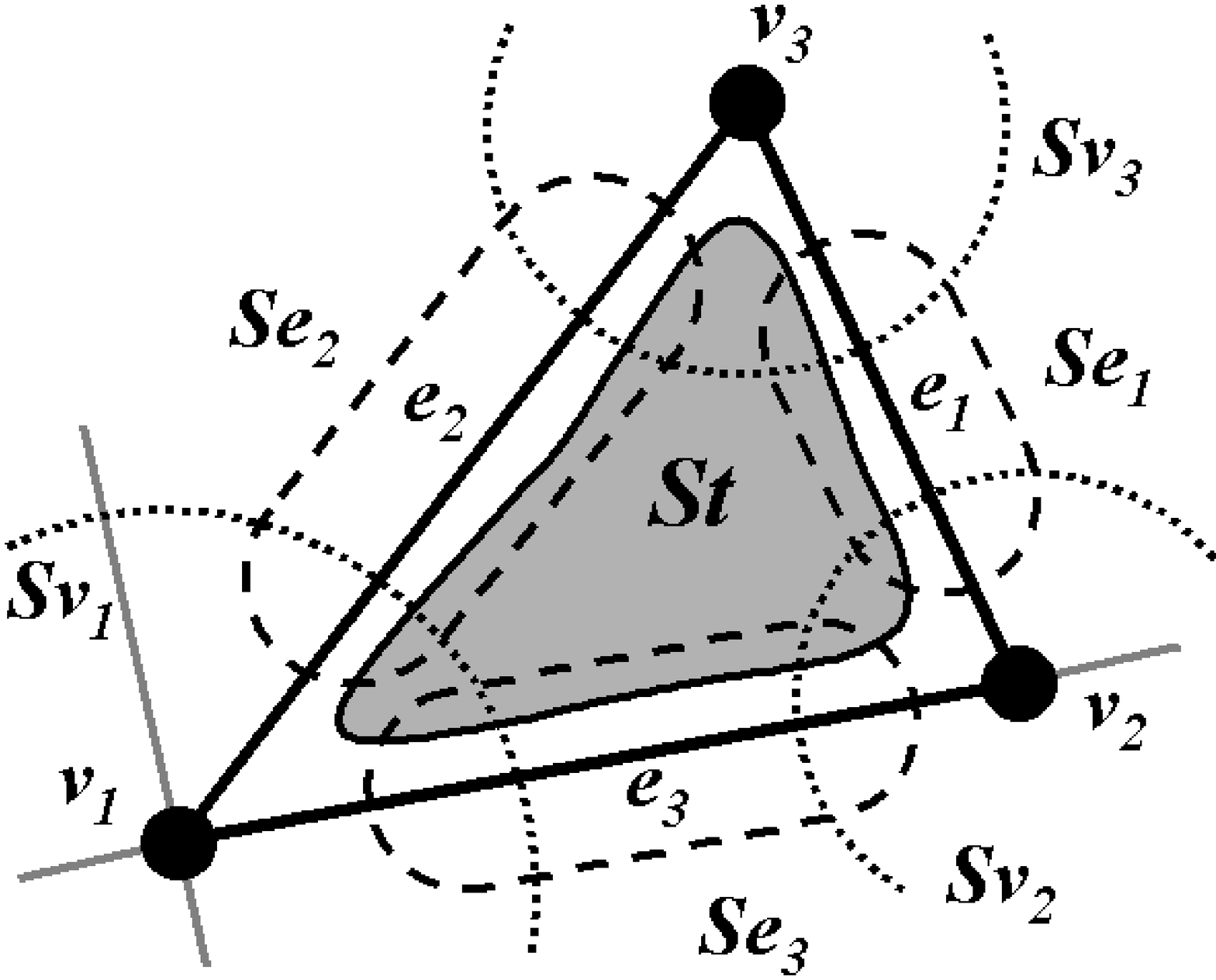}} \caption { }
\end{figure}

Finally, if $t_{1},...,t_{m_{2}}$ are the triangles of $M$, for
$k=1,...,m_{2} $ the 2-sphere $St_{k}$ bounds a 3-ball $Bt_{k}$ as
in Figure \ref{fig15a}. The 3-ball $Bt_{k}$ is contained in the
(open) star $star\left(  t_{k}\right)  $ and it intersects $t_{k}$
in a closed disk $\tilde{t}_{k}\subset t_{k}$, and the
intersection of $St_{k}$ with $t_{k}$ is transverse. The 3-ball
$Bt_{k}$ is disjoint with $B\epsilon$ for every $\epsilon\in
T^{2}$ different from $t_{k}$ unless $\epsilon$ is incident with
$t_{k}$. In this case, $Bt_{k}$ and $B\epsilon$ intersect
normally. Moreover, if we have $v_{i}<e_{j}<t_{k}$, then the
3-balls $Bv_{i},Be_{j},Bt_{k}$ intersect normally (Figure
\ref{fig15b}) and there is one of the two triple points of
$Sv_{i}\cap Se_{j}\cap St_{k}$ in each of the two tetrahedra of
$star\left( t_{k}\right)  $.

If $T$ is a cell decomposition of $M$ instead of a triangulation, the previous
construction can easily be generalized.

It is easy to check that the Dehn surface
$\mathbb{T}=\underset{\epsilon\in T^{2}}{\cup}S\epsilon$ so
constructed is a filling collection of spheres in $M$. Moreover,
$\mathbb{T}$ is regular and it is transverse to the (smooth)
simplexes of the triangulation $T$ of $M$.

In particular, as a Corollary of Theorem \ref{THM fillingcollection implies
fillingSPHERE} this construction implies the main Theorem of \cite{Montesinos}:

\begin{theorem}
$M$ has a nulhomotopic filling Dehn sphere.
\end{theorem}

Note that in this case, in contradistinction of \cite{Montesinos} or
\cite{Anewproof}, we have not made any assumption about the orientability of
$M$.

The following result follows directly from the construction.

\begin{proposition}
\label{PROPinflatetriangulationsFILLS}Let $S$ be a surface and $f:S\rightarrow
M$ a transverse immersion. Let $K,T$ be triangulations of $S,M $ respectively
such that $f$ is simplicial with respect to them. Then $f\left(  S\right)
\cup\mathbb{T}$ is a regular filling surface of $M$.
\end{proposition}

In the previous Proposition the immersion $f$ can be any transverse immersion,
filling or not. Assume now that $f:S\rightarrow M$ is a filling immersion and
put $\Sigma:=f(S)$. Let $K,T$ be triangulations of $S$ and $M$ respectively
such that $f$ is simplicial with respect to them. By the previous Proposition,
$\Sigma\cup\mathbb{T}$ fills $M$, and by the same methods of the proof of
Theorem \ref{THM fillingcollection implies fillingSPHERE}, we can obtain from
$\Sigma\cup\mathbb{T}$ a unique filling Dehn surface of $M$. If we look at the
proof of Theorem \ref{THM fillingcollection implies fillingSPHERE}, this can
be done in many different ways because there are many possibilities for
performing the spiral pipings. We say that each filling Dehn sphere
$\Sigma^{\prime}$ of $M$ that is obtained from $\Sigma$ in this way is a
$T$\textit{-inflating} of $\Sigma$. Let $\Sigma^{\prime}$ be a $T$-inflating
of $\Sigma$. By Proposition \ref{PROPinflatetriangulationsFILLS} and
Proposition \ref{PROPSpiralPipingPreserveFillingness}, $\Sigma^{\prime}$ is
regular because spiral pipings preserve regularity. There is an immersion
$f^{\prime}:S\rightarrow M$ parametrizing $\Sigma^{\prime}$ that \textit{comes
from} $f$ in a natural way, that is, $f^{\prime}$ agrees with $f$ in most of
$S$ except in the small disks where we perform the pipings. We say also that
$f^{\prime}$ is a $T$\textit{-inflating} of $f$. The first application of
shellability is the next result proved in \cite{RHomotopies}.

\begin{proposition}
\label{PROPshellableimpliesfillinghomotopictoinflated}If $T$ shells $f$, then
there is a $T$-inflating $f^{\prime}$ of $f$ filling homotopic to $f$.
Moreover, we can choose $f^{\prime}$ such that there are only two spiral
pipings connecting $\Sigma$ with components of $\mathbb{T}$ and such that the
rest of spiral pipings are performed around triple points of $\mathbb{T}$.
\end{proposition}

By Theorem \ref{THMsmoothtriangulations exist}, passing to suitable
subdivisions we can assume that $T$ shells $f$, and thus we have:

\begin{corollary}
\label{CorolarYfillinghomotopictoREGULAR}If $f:S\rightarrow M$ is a filling
immersion, then $f$ is filling homotopic to a regular filling immersion.
\end{corollary}

The proof of Proposition
\ref{PROPshellableimpliesfillinghomotopictoinflated} is made by
repeatedly applications of the constructions of sections
\ref{SUBSECTION InflatingDoublePoint} and \ref{SUBSECTION
InflatingDisks}, using that each region of $M-\Sigma$ has a
shellable triangulation and that each triangulation of a 2-disk is
shellable \cite{Sanderson}. As an example, we will illustrate the
starting point of this construction in which we ''inflate'' a
tetrahedron of $T$.

\begin{example}
\label{EXAMPLE Inflating a Tetrahedron}Inflating a tetrahedron.

\vspace {10pt}
\begin{figure}[htb]
\centering
\includegraphics[ width=0.45\textwidth ]{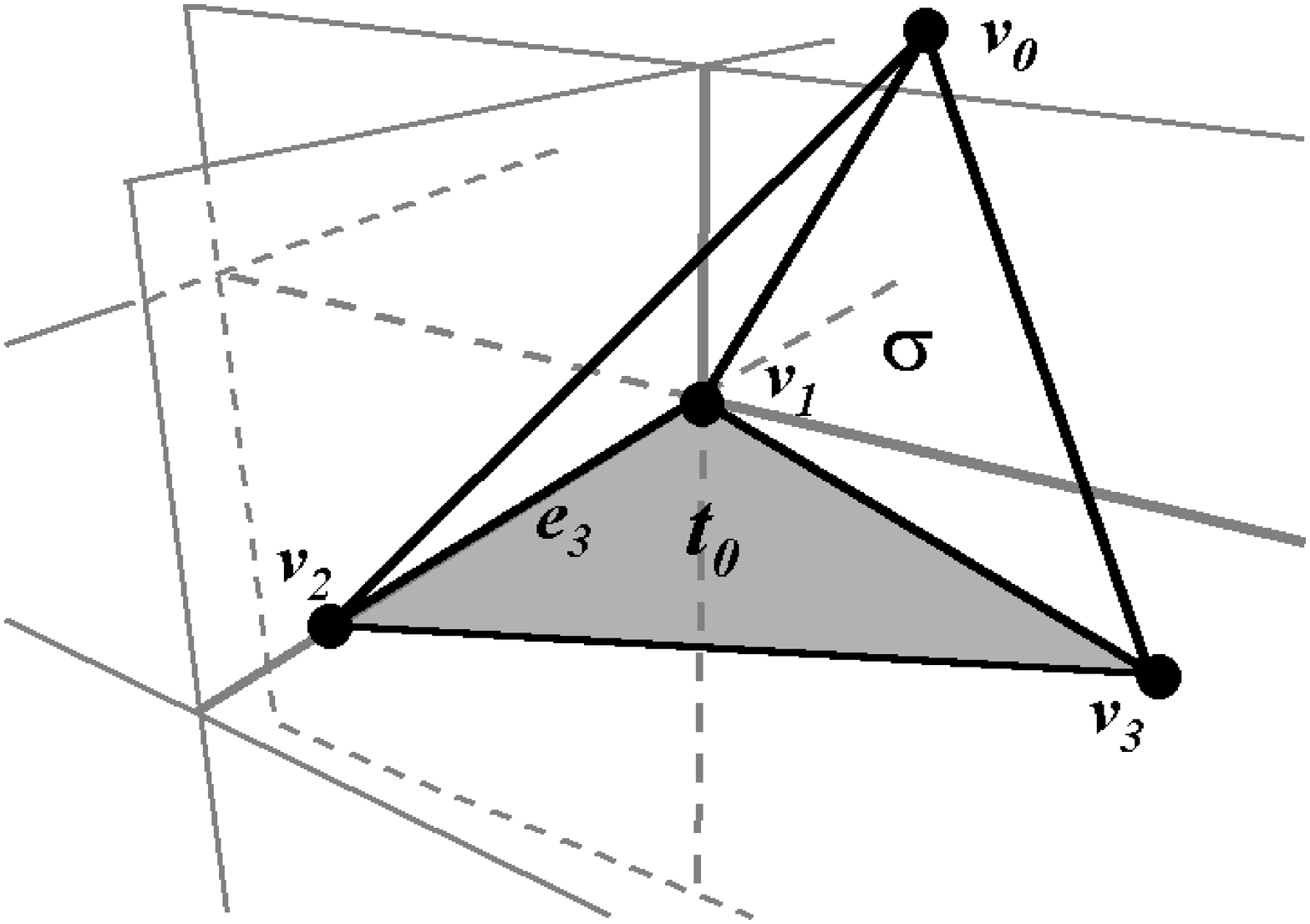} \caption { } \label{fig16}
\end{figure}

Let $f$, $\Sigma$ and $T$ be as in Proposition
\ref{PROPshellableimpliesfillinghomotopictoinflated}. Let $R$ be a
region of $M-\Sigma$. Because $T$ induces a shellable
triangulation on $R$, there is a tetrahedron $\sigma$ of $T$ such
that $cl(\sigma)$ is free in $cl(R)$. Assume that $\sigma$ is as
in Figure \ref{fig16}. In this case, the intersection
$cl(\sigma)\cap cl(R)$ is the closure of a triangle $t_{0}$ of
$\sigma$, and $t_{0}$ has exactly one vertex $v_{1}$ which is a
triple point of $\Sigma$ and exactly one edge $e_{3}$ (incident
with $v_{1}$) contained in a double curve of $\Sigma$.

We can think that the triangle $t_{0}$ is the triangle $t$ of
Figure \ref{fig15b} before.

\begin{enumerate}
\item  Consider the intersection point $Q_{0}$ of the 2-sphere $Sv_{1}$ with
the edge $e_{3}$. We inflate $Q_{0}$ to obtain a small 2-sphere
$\Sigma _{Q_{0}}$ piped with $\Sigma$ (Fig. \ref{fig17}(b)). After
a finger move +2 through $v_{1}$ (Fig. \ref{fig17}(c)) and ambient
isotopy of $M$ (Fig. \ref{fig17}(d)), $\Sigma_{Q_{0}} $ is
transformed into $Sv_{1}$ and $\Sigma\#\Sigma_{Q_{0}}$ is
transformed into $\Sigma\#Sv_{1}$.

\vspace {10pt}

\begin{figure}[htb]
\centering
\includegraphics[ width=0.8\textwidth
]{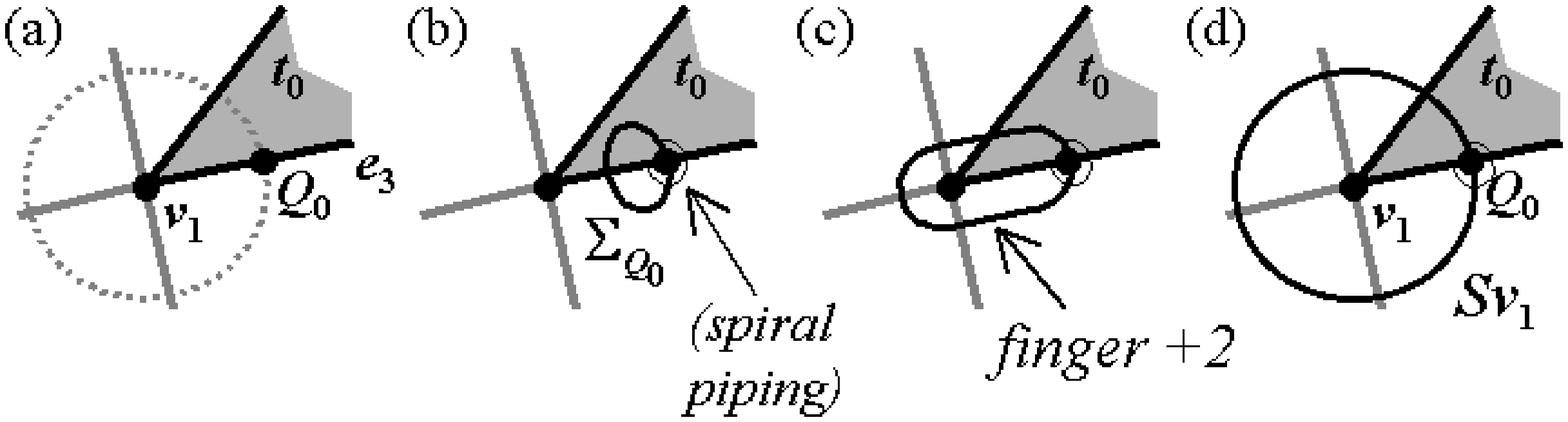} \caption { } \label{fig17}
\end{figure}

\item  Take the point $Q_{1}$ of intersection of $Se_{3}$ with $e_{3}$ that
lies inside $Bv_{1}$ (Fig. \ref{fig18}(a)). We inflate $Q_{1}$
(Fig. \ref{fig18}(b)) and apply a piping passing move through
$Q_{0}$ (Fig. \ref{fig18}(c)) to obtain $\Sigma \#Sv_{1}\#Se_{3}$
(Fig. \ref{fig18}(d)).

\item  Take the triple point $P_{0}$ of $\mathbb{T}$ in the intersection
$Sv_{1}\cap Se_{3}\cap St_{0}$ that lies in $\sigma$ (Fig.
\ref{fig19a}). This point $P_{0}$ is now a double point of
$\Sigma\#Sv_{1}\#Se_{3}$. We inflate $P_{0}$. If $A$ is the
intersection point of the double curve $Sv_{1}\cap Se_{3}$ of
$\mathbb{T}$ with the triangle $t_{0}$, after a finger move +2
through $A$ (Fig. \ref{fig19b}) and an ambient isotopy of $M$
(Figs. \ref{fig19c} and \ref{fig19d}), $\Sigma_{P_{0}}$ is
transformed into $St_{0}$ and $\Sigma\#Sv_{1}\#Se_{3}$ into
$\Sigma\#Sv_{1}\#Se_{3}\#St_{0}$.

\vspace {10pt}

\begin{figure}[htb]
\centering
\includegraphics[ width=0.8\textwidth
]{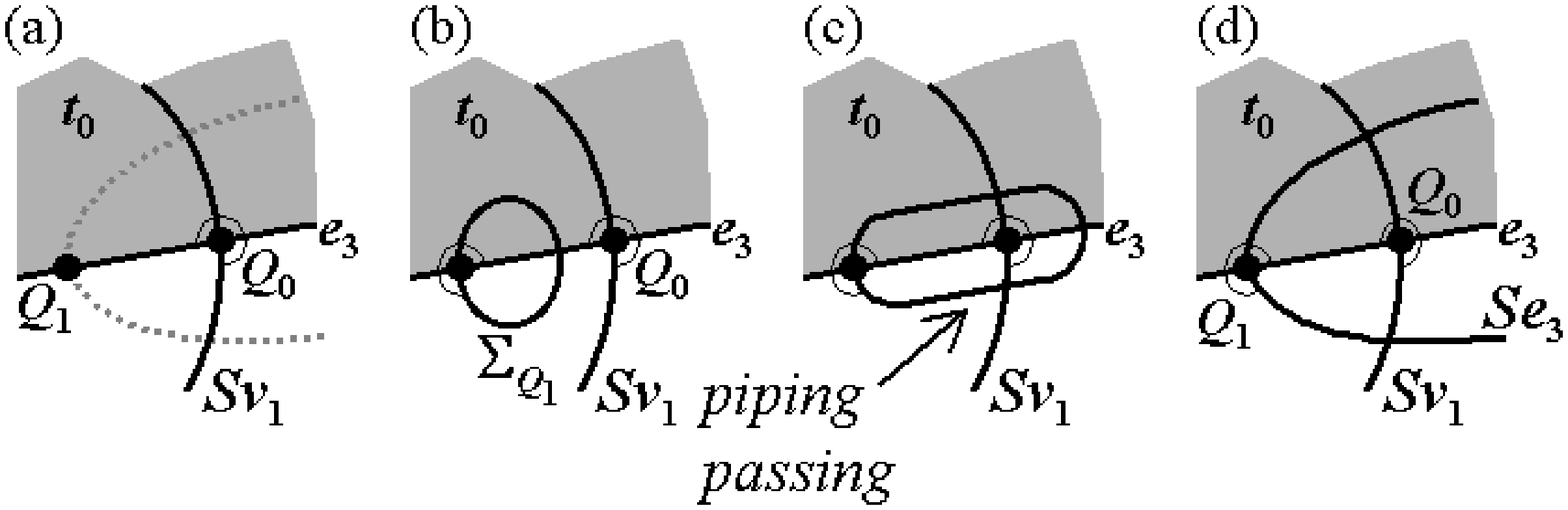} \caption { } \label{fig18}
\end{figure}

\item  Consider now the other vertex $v_{2}$ of $e_{3}$ different from $v_{1}
$. Take the triple point $P_{1}$ of $\mathbb{T}$ in the
intersection $Sv_{2}\cap Se_{3}\cap St_{0}$ that lies in $\sigma$.
We inflate $P_{1}$ to obtain the (piped) 2-sphere $\Sigma_{P_{1}}$
contained in $Bv_{2}$. We need now two consecutive finger moves +2
and an ambient isotopy of $M$ for transforming $\Sigma_{P_{1}}$
into $Sv_{1}$ (Fig. \ref{fig20}).

\item  In a similar way as in 4, if $e_{1}$ is the edge of $t_{0}$ incident
with $v_{2}$ and different from $e_{3}$, we inflate $Se_{1}$ from
the intersection $Sv_{2}\cap St_{0}$ (Fig. \ref{fig21}(a)). After
this, we inflate $Sv_{3} $ from $Se_{1}\cap St_{0}$ and $Se_{2}$
from $Sv_{3}\cap St_{0}$ (Fig. \ref{fig21}(b)). After a final
finger move +2 (Fig. \ref{fig21}(c)) we have just inflated all the
2-spheres $S\epsilon$ for $\epsilon<t_{0}$.

\begin{figure}[htb]
\centering \subfigure []{ \label{fig19a}
\includegraphics[ width=0.35\textwidth ] {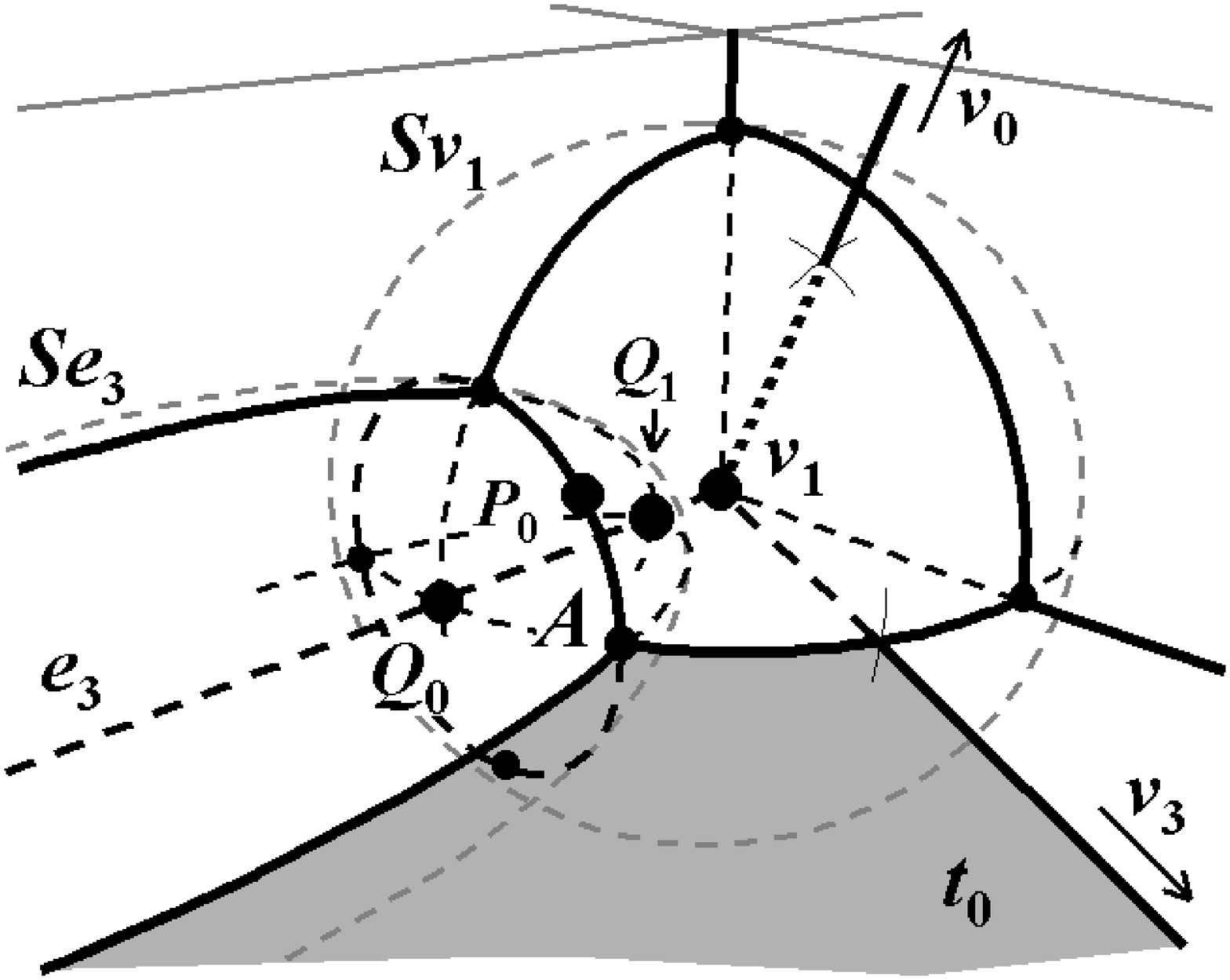}}\hspace{1cm} \subfigure[]{ \label{fig19b}
\includegraphics[ width=0.35\textwidth ]{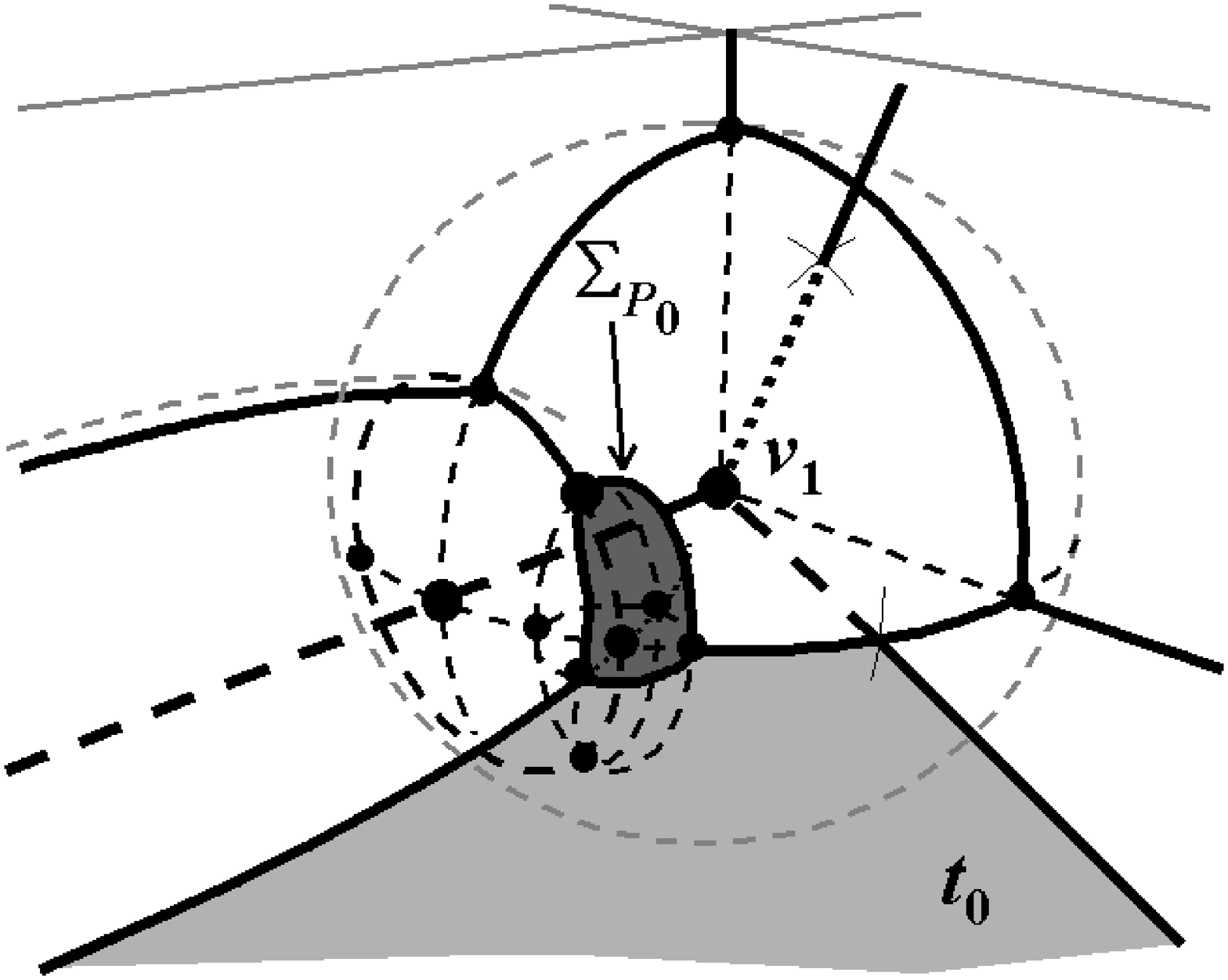}}\\
\vspace {10pt} \subfigure []{ \label{fig19c}
\includegraphics [
height=0.15\textheight ] {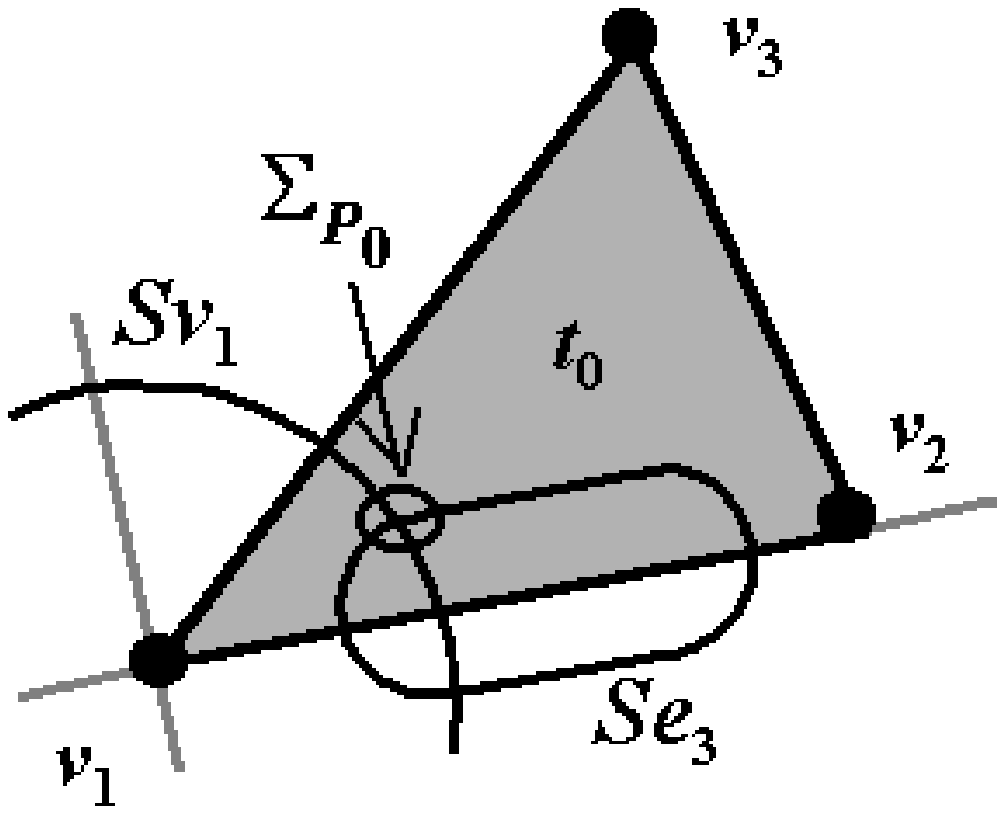}}\hspace{0.5cm} \subfigure{
\includegraphics [
height=0.15\textheight ] {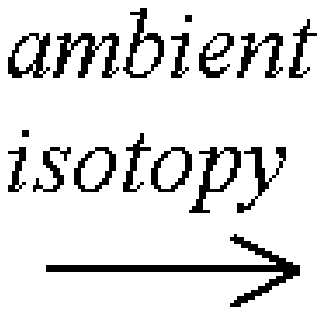}}\hspace {0.5cm}
\addtocounter{subfigure}{-1} \subfigure [] {\label{fig19d}
\includegraphics[
height=0.15\textheight ]{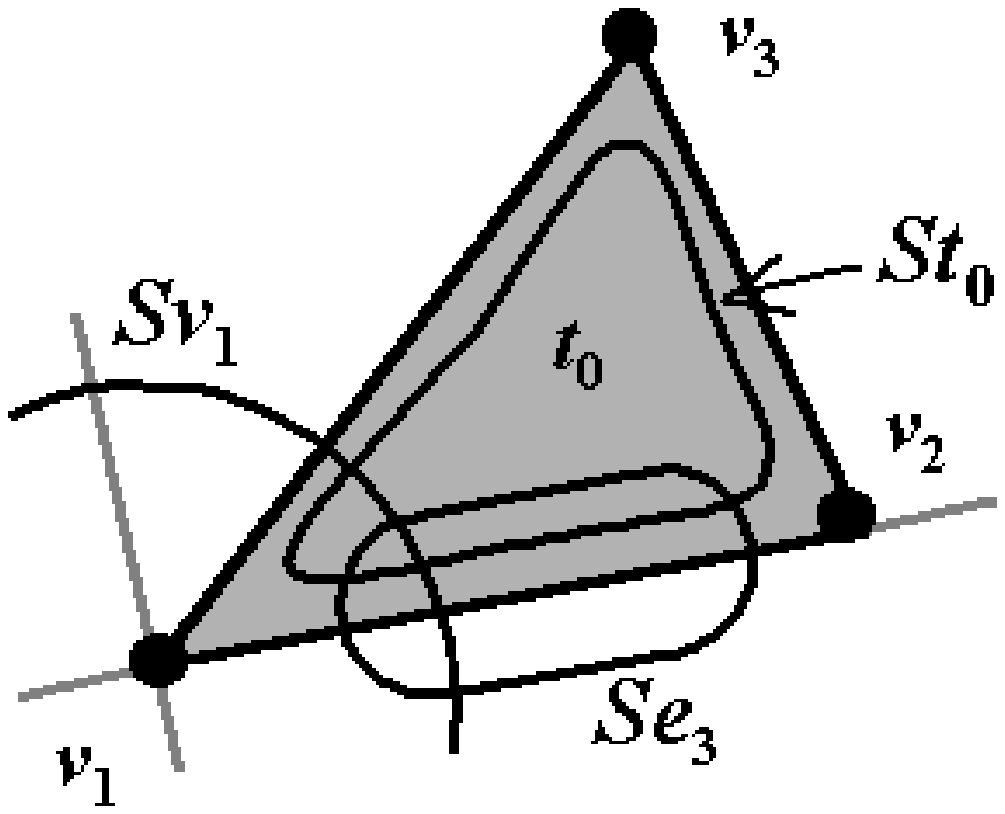}} \caption { }
\end{figure}

\item  Now, consider the triangle $t_{1}$ of $\sigma$ different from $t_{0}$
which is incident with $e_{3}$. Let $P_{2}$ be the triple point of $Sv_{1}\cap
Se_{3}\cap St_{1}$ that lies in $\sigma$. We inflate $P_{2}$ to obtain the
2-sphere $\Sigma_{P_{2}}$ contained in $Bt_{1}$. The triangle $t_{1}$ is not
contained in the Dehn surface $\Sigma$, and because of this we need only a
finger move +1 to make $\Sigma_{P_{2}}$ cross $Sv_{2}$, and then an ambient
isotopy of $M$ for transforming $\Sigma_{P_{2}}$ into $St_{1}$.

\vspace {10pt}

\begin{figure}[htb]
\centering
\includegraphics[ width=0.8\textwidth
]{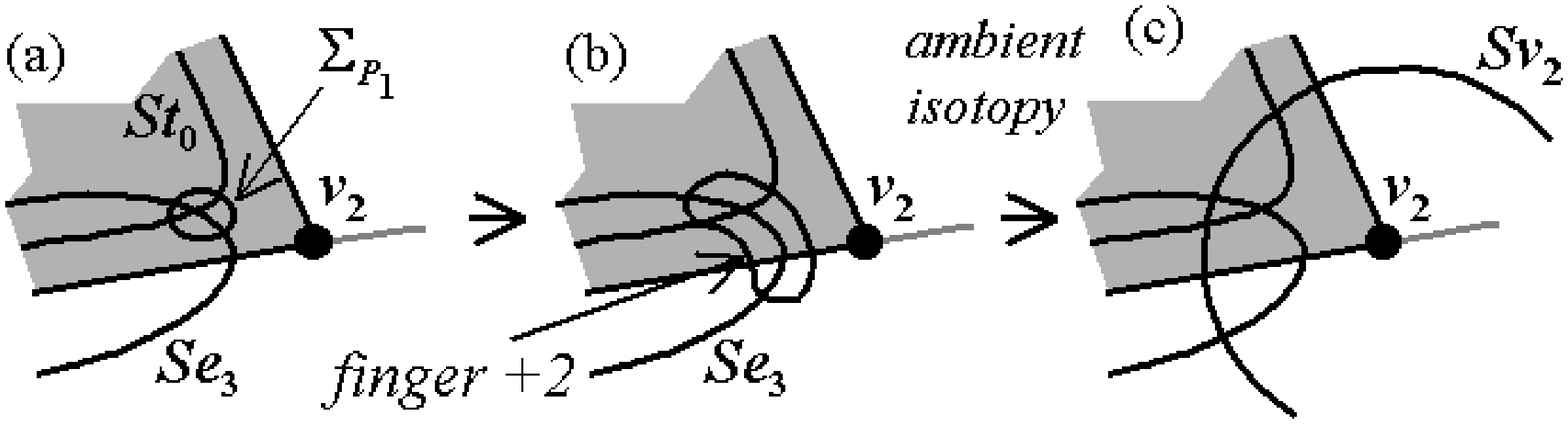} \caption { } \label{fig20}
\end{figure}

\item  Now, in a similar (but simpler because $t_{1}\not \subset\Sigma$) way
as in 4 and 5, we inflate inductively all the remaining 2-spheres $S\epsilon$
for $\epsilon<t_{1}$.

\item  If $t_{2}$ is another triangle of $\sigma$ different from $t_{1}$, in a
similar way as in 6 and 7 we inflate the 2-sphere $St_{2}$ and the remaining
2-spheres $S\epsilon$ for $\epsilon<t_{2}$. When inflating $St_{2}$ there will
be a slight difference with the case of $t_{1}$ because we have just inflated
the 2-spheres corresponding to two edges of $t_{2}$.

\vspace {10pt}

\begin{figure}[htb]
\centering
\includegraphics[ width=0.8\textwidth
]{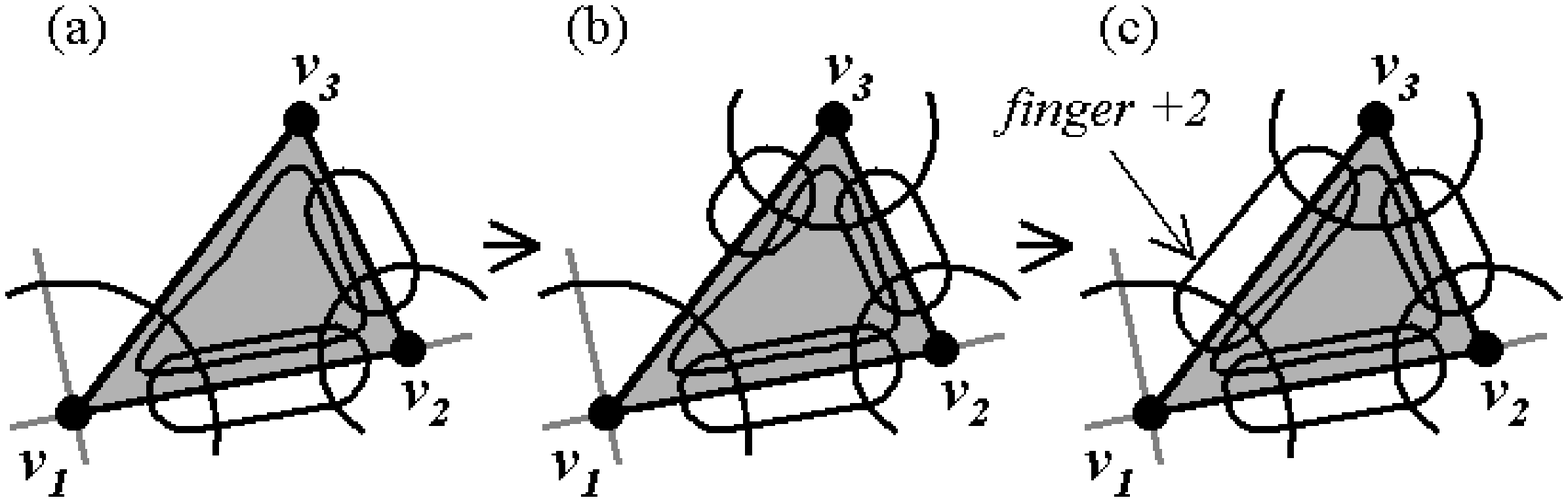} \caption { } \label{fig21}
\end{figure}

\item  For the final triangle $t_{3}$, we have just inflated all the 2-spheres
$S\epsilon$ for $\epsilon<t_{3}$ with the exception of $St_{3}$. If
$\Sigma^{\prime}$ is the filling Dehn surface that we have at this moment, we
can see that the closure of the part of $t_{3}$ lying outside the 3-balls
$B\epsilon$'s with $\epsilon<t_{3}$ is a 6-gon $w$ in $M-\Sigma^{\prime}$.
Inflating $w$ as in section \ref{SUBSECTION InflatingDisks} and after an
ambient isotopy of $M$ we get also $St_{3}$.
\end{enumerate}
\end{example}

In the previous example, we see that the way of constructing the filling Dehn
surface $\Sigma^{\prime}$ from the filling Dehn surface $\Sigma$ in
Proposition \ref{PROPshellableimpliesfillinghomotopictoinflated} is in some
sense to make $\Sigma$ \textit{grow} inductively following a path given by the
triangulation $T$. The growing path we will follow for proving Proposition
\ref{PROPshellableimpliesfillinghomotopictoinflated} in \cite{RHomotopies}
will not be exactly as in Example \ref{EXAMPLE Inflating a Tetrahedron}. There
(in \cite{RHomotopies}) we will inflate first from $\Sigma$ all the 2-spheres
of $\mathbb{T}$ corresponding to the simplexes of $T$ contained in $\Sigma$
starting in a similar way as in Steps 1 to 5 of the Example \ref{EXAMPLE
Inflating a Tetrahedron}. Then, the shellability conditions imposed to $T$
will give us the growing path of $\Sigma$ on the regions of $M-\Sigma$ using
similar methods to that of Steps 6 to 9 of Example \ref{EXAMPLE Inflating a
Tetrahedron}.

We will say that the $T$-inflating $\Sigma^{\prime}$ of $\Sigma$ ($f^{\prime}$
of $f$) as in previous Proposition is a $T$\textit{-growth} of the filling
Dehn sphere $\Sigma$ (of the filling immersion $f$). Note that to be a
$T$-growth is stronger than to be a $T$-inflating.

The next (but not the last) application of the constructions of section
\ref{SECTION Shellability} is the following.

Let $f,g:S\rightarrow M$ be transverse immersions that differ by
the pushing disk $(D,B)$, and assume that $f$ is a filling
immersion. In this situation, the immersion $g$ will not
necessarily be a filling immersion. Consider triangulations $K,T$
of $S$ and $M$ respectively such that $f$ and $g$ are simplicial
with respect to them. Take a $T$-inflating $f^{\prime}$ of $f$
such that $f^{\prime}$ agree with $f$ over $D$ (because
$f^{\prime}$ agrees with $f$ in most of $S$ except in some small
disks of $S$, we require that these small disks do not intersect
$D$), and consider the immersion $g^{\prime}$ that is obtained
from $f^{\prime}$ after the pushing disk $(D,B)$. We can assume
that $g^{\prime}$ agree with $g$ except in the disks of $S$ where
$f^{\prime}$ ''disagrees'' with $f^{\prime}$. The Dehn surface
$g^{\prime}(S)$ is obtained from $g(S)\cup\mathbb{T}$ by spiral
pipings, and because $g$ is also simplicial (with respect to
$K,T$), $g^{\prime}$ is a $T$-inflating of $g$.

\begin{figure}[htb]
\centering \subfigure []{ \label{fig22a}
\includegraphics[
width=0.4\textwidth ] {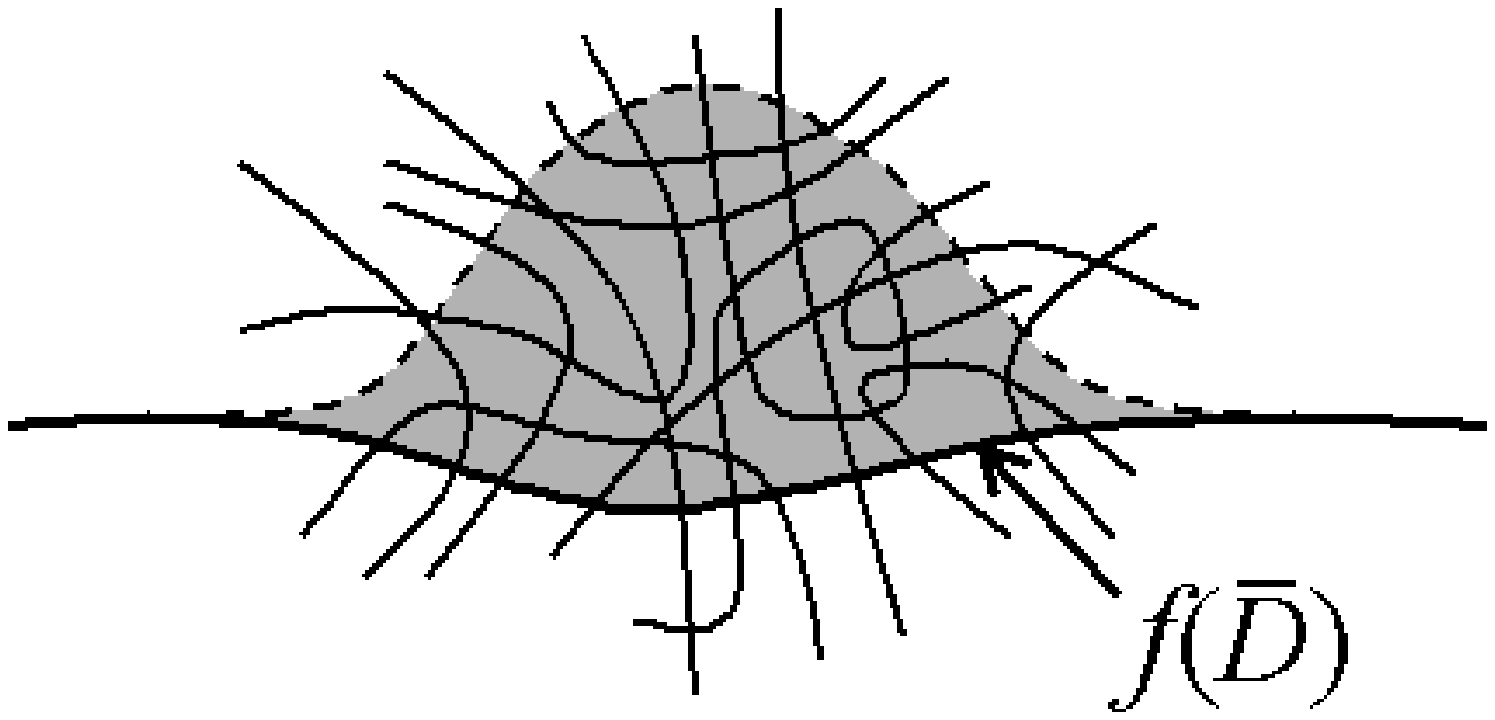}}\hfill \subfigure[]{ \label{fig22b}
\includegraphics[ width=0.4\textwidth
]{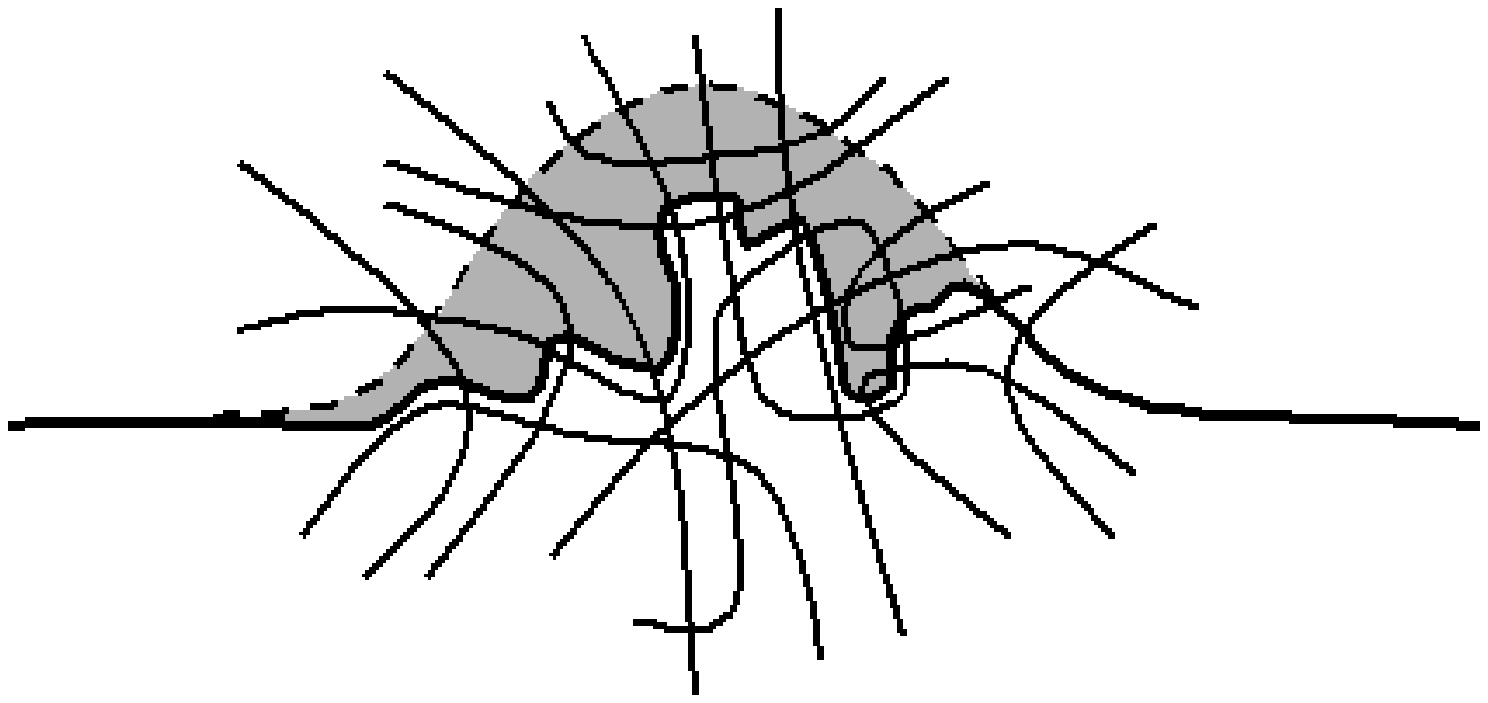}}\\
\vspace {10pt} \subfigure []{ \label{fig22c}
\includegraphics [
width=0.4\textwidth ] {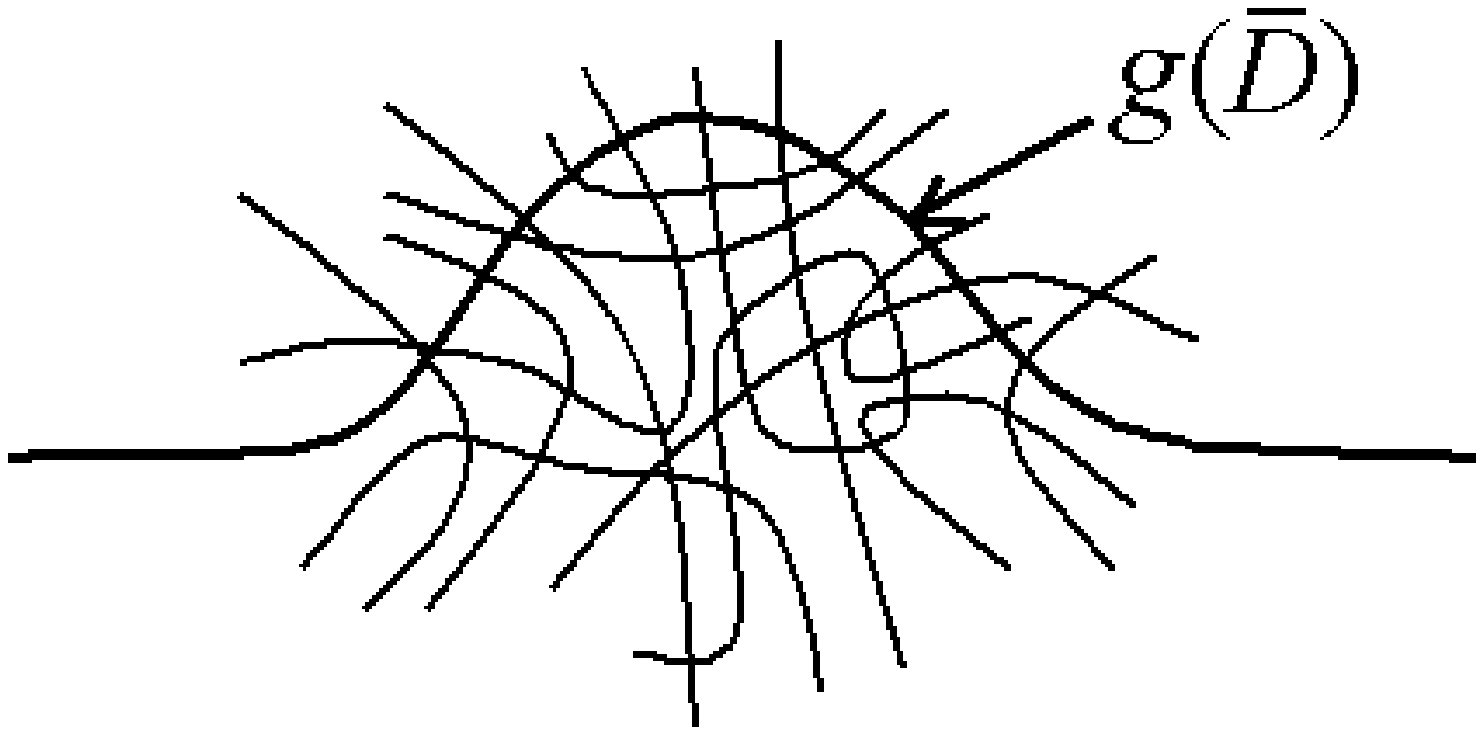}} \caption { }
\end{figure}

Because both $f,g$ are simplicial (with respect to $K,T$), the triangulation
$T$ induces a triangulation of the pushing ball $B$.

\begin{lemma}
[Key Lemma 2]\label{KeyLemma2}If the (induced) triangulation of $B$ collapses
simplicially into $g(D)$, then $f^{\prime}$ is filling homotopic to
$g^{\prime}$.
\end{lemma}

\begin{proof}
[Sketch of the Proof]Because $B$ is triangulated by $T$, it can be
easily shown that $f(S)\cup g(S)\cup\mathbb{T}$ induces a cell
decomposition of $B$. If $B$ collapses simplicially into $g(D)$,
then we can define a special shelling of this cell decomposition
of $B$ induced by $f(S)\cup g(S)\cup \mathbb{T}$. This special
shelling has the property that it will allow us to apply Lemma
\ref{LemmaPassing3-cells} repeatedly (Figure \ref{fig22b}) to the
filling Dehn sphere $f(S)\cup\mathbb{T}$ until we get
$g(S)\cup\mathbb{T}$. Substituting finger moves 2 by piping
passing moves where it is required, this deformation of
$f(S)\cup\mathbb{T}$ into $g(S)\cup\mathbb{T}$ defines also a
deformation of $f^{\prime}$ into $g^{\prime}$ by filling
preserving moves.
\end{proof}

\section{Filling pairs.\label{SECTION FillingPairs}}

Let $\Sigma_{1}$ and $\Sigma_{2}$ be two filling Dehn surfaces of $M$. Assume
by simplicity that both are regular.

If we are not given more information about $\Sigma_{1}$ and $\Sigma_{2}$, we
don't know how they are related to each other. The only thing we can say, if
$M$ is not $S^{3}$, is that they must have nonempty intersection.

\begin{definition}
\label{DEF Fillingpair}We say that $\Sigma_{1}$ and $\Sigma_{2}$ form a
filling pair in $M$ if their union $\Sigma_{1}\cup\Sigma_{2}$ is also a
regular filling Dehn surface $M$.
\end{definition}

In particular, if $\Sigma_{1}$ and $\Sigma_{2}$ form a filling pair in $M$,
they intersect transversely.

If $\Sigma_{1}$ and $\Sigma_{2}$ are a filling pair in $M$, then $\Sigma_{2}$
induces a cell decomposition on the closure of each region of $M-\Sigma_{1}$
and viceversa. Because both $\Sigma_{1}$ and $\Sigma_{2}$ are regular, all
these induced cell decompositions are also regular. If $R_{1}$ is a region of
$M-\Sigma_{1}$, then we say that $\Sigma_{2}$ \textit{shells} $R_{1}$ if
$\Sigma_{2}$ induces a shellable cell decomposition of the 3-ball $cl(R_{1})$.
We say that $\Sigma_{2}$ \textit{shells} $\Sigma_{1}$ if $\Sigma_{2}$ shells
each region of $M-\Sigma_{1}$.

\begin{definition}
\label{DEF MutuallyShellable}Let $\Sigma_{1}$ and $\Sigma_{2}$ form a filling
pair in $M$. We say that $\Sigma_{1}$ and $\Sigma_{2}$ are mutually shellable
if $\Sigma_{1}$ shells $\Sigma_{2}$ and $\Sigma_{2}$ shells $\Sigma_{1}$.
\end{definition}

The following result is proved in detail in \cite{RHomotopies}.

\begin{proposition}
\label{PROPfillinghomotopictoMutuallyShellable}Let $\Sigma_{1},\Sigma_{2}$ be
regular filling Dehn surfaces of $M$ such that they intersect transversely. If
$f_{1}:S_{1}\rightarrow M$ parametrizes $\Sigma_{1}$, then $f_{1}$ is filling
homotopic to an immersion $f_{1}^{\prime}:S_{1}\rightarrow M$ such that
$\Sigma_{1}^{\prime}:=f^{\prime}(S_{1})$ and $\Sigma_{2}$ form a mutually
shellable filling pair in $M$.
\end{proposition}

\begin{proof}
[Sketch of the Proof]Let $f_{2}:S_{2}\rightarrow M$ be a
parametrization of $\Sigma_{2}$, and let $T$ be a good
triangulation of $M$ with respect to $f_{1}$ and $f_{2}$
(definition \ref{DEF GoodTriangulations}). Then, $T$ shells every
region of $M-\Sigma_{1}$ and every region of $M-\Sigma_{2}$. The
union $\Sigma_{1}\cup\Sigma_{2}\cup\mathbb{T}$ is a regular
filling Dehn surface of $M$ by Proposition
\ref{PROPinflatetriangulationsFILLS}. Take a $T$-growth
$f_{1}^{\prime}$ of $f_{1}$, and put $\Sigma_{1}^{\prime}
=f_{1}^{\prime}(S_{1})$. We make the spiral pipings that transform
$\Sigma _{1}\cup\mathbb{T}$ into $\Sigma_{1}^{\prime}$ small
enough, such that they do not intersect $\Sigma_{2}$. Because
regularity is preserved by spiral pipings, it is not difficult to
see that $\Sigma_{1}^{\prime}\cup\Sigma_{2}$ is a regular filling
Dehn surface of $M$. It is also easy to see that $\Sigma_{2} $
induces a shellable cell decomposition on every region of
$M-\Sigma _{1}^{\prime}$ using that $\Sigma_{2}$ is a subcomplex
of $T$ and the construction of $\mathbb{T}$. The non-trivial part
is to check that $\Sigma_{1}^{\prime}$ induces a shellable cell
decomposition on every component of $M-\Sigma_{2}$. This is made
in detail in \cite{RHomotopies}, and it is parallel to the proof
of Key Lemma 2 above. If $R_{2}$ is a region of $M-\Sigma_{2}$,
the first thing that it is checked is that $\mathbb{T} \;$induces
a shellable cell decomposition on $cl(R_{2})$. This is done
following the proof of Lemma \ref{LEMABing} in \cite{Bing}, using
that the restriction to $cl(R_{2})$ of the triangulation $T$ is
simplicially collapsible. After this, it is seen that the presence
of $\Sigma_{1}$ do not alter this shellability property. Finally,
the presence of the spiral pipings might affect the shelling in
some cases. In \cite{RHomotopies} it is explained how this can
occur, and how we can modify locally $f_{1}^{\prime}$ by
filling-preserving moves until $f_{1}^{\prime}$ verifies the
statement of Proposition
\ref{PROPfillinghomotopictoMutuallyShellable}.
\end{proof}

\section{Simultaneous growings.\label{SECTION Growing}}

The last application of shellability will be the following one.

\begin{figure}[htb]
\centering \subfigure []{ \label{fig23a}
\includegraphics[ width=0.4\textwidth ] {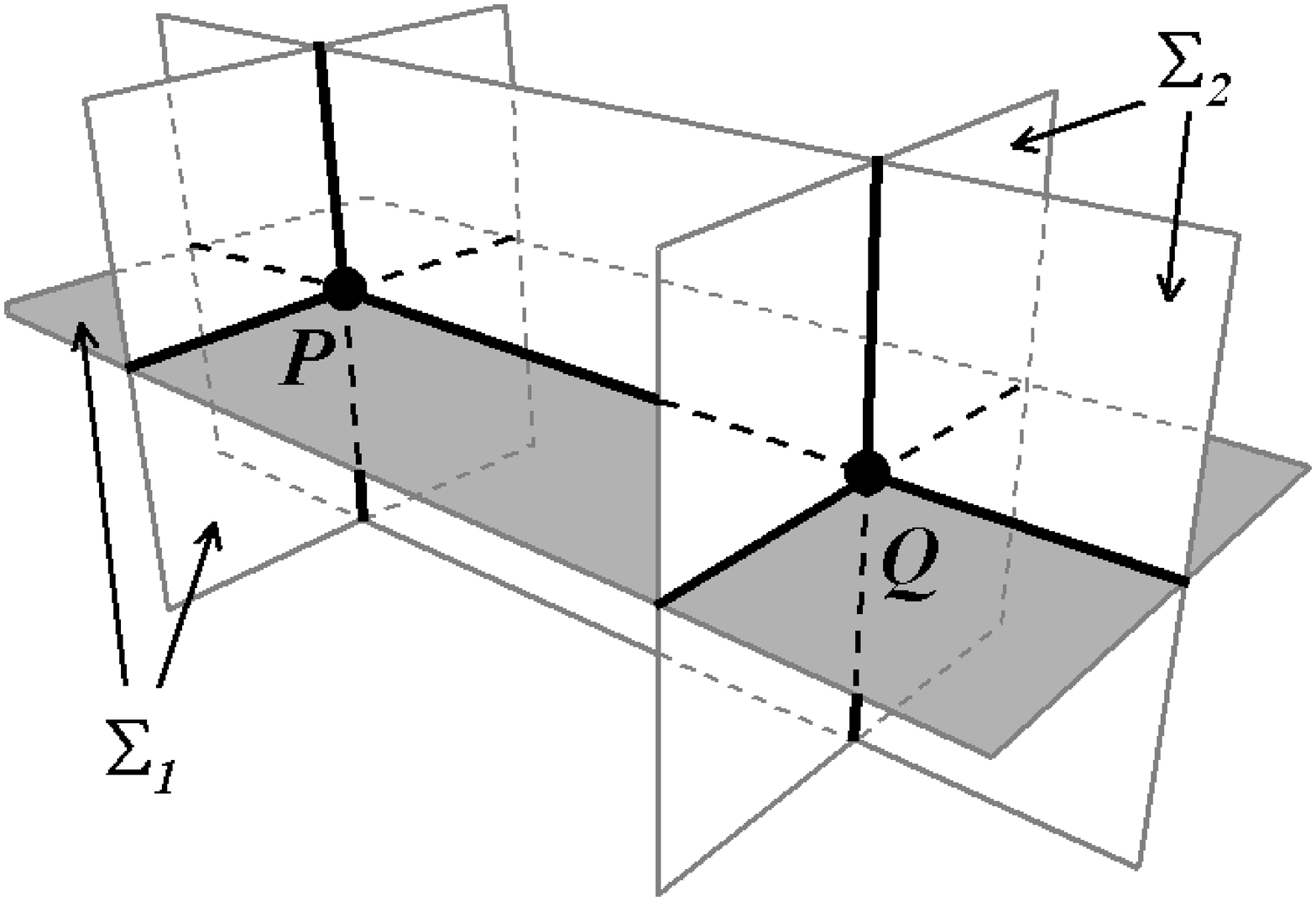}}\hfill \subfigure {
\includegraphics[ height=0.15\textheight, width=0.05\textwidth ] {fig35ab} }\hfill
\addtocounter {subfigure}{-1} \subfigure[]{ \label{fig23b}
\includegraphics[ width=0.4\textwidth ] {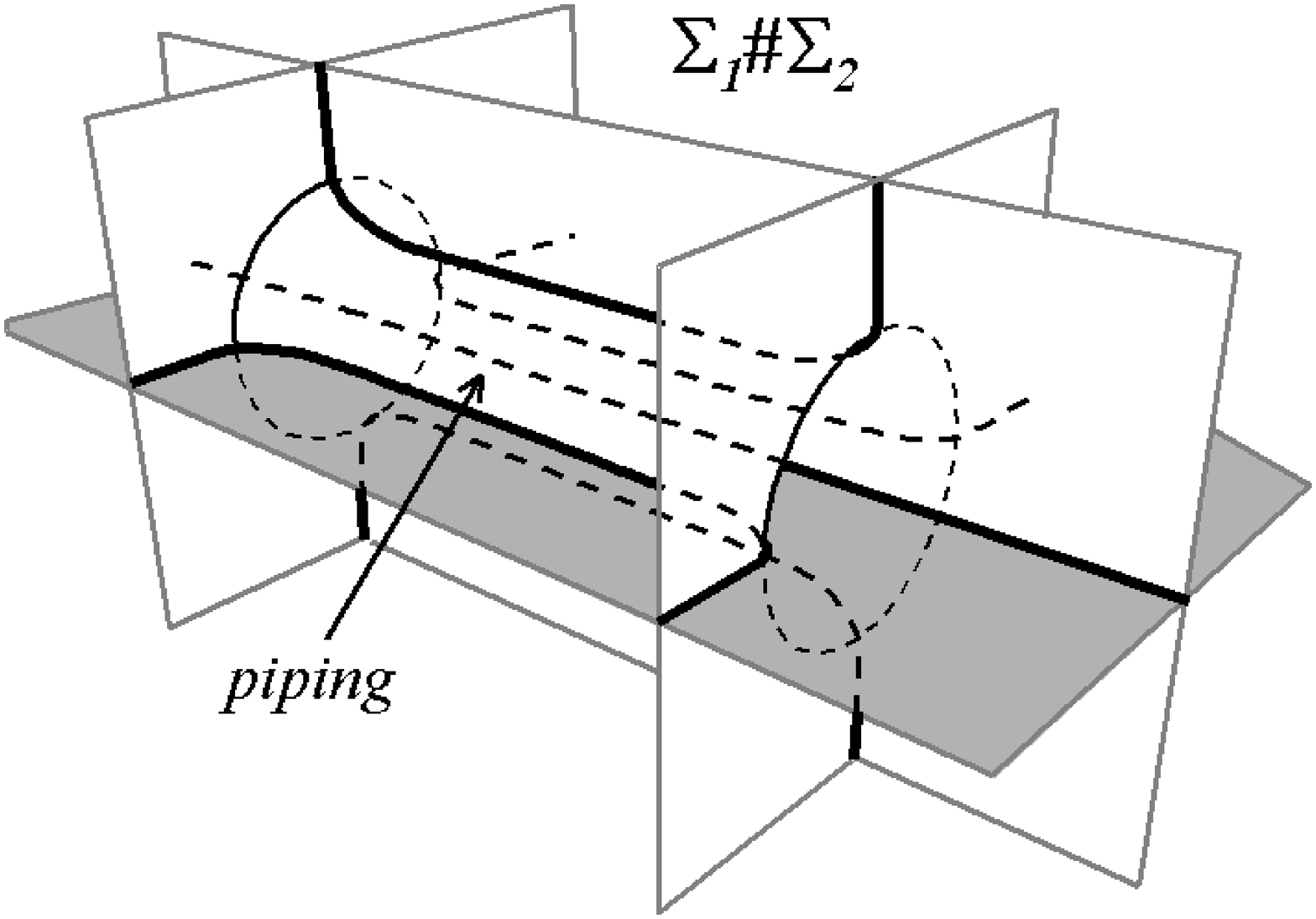}} \caption { }
\end{figure}

Assume that $\Sigma_{1},\Sigma_{2}$ are two regular filling Dehn
surfaces of $M$ and that there are two points $P,Q$ where
$\Sigma_{1}$ and $\Sigma_{2}$ intersect as in Figure \ref{fig23a}.
We call $\Sigma_{1}\#\Sigma_{2}$ to the Dehn surface of $M$ that
arises piping $\Sigma_{1}$ with $\Sigma_{2}$ near $P,Q$ as in
Figure \ref{fig23b}. We assume also that the points $P,Q$ have the
property that $\Sigma_{1}\#\Sigma_{2}$ is another filling Dehn
surface of $M$. Let $f:S_{1}\rightarrow M,g:S_{2}\rightarrow M$ be
parametrizations of $\Sigma _{1},\Sigma_{2}$ respectively.
Consider the two small disks $\delta_{1} ,\delta_{2}$ of $S_{1}$
and $S_{2}$ respectively whose respective images by $f$ and $g$
disappear after the piping. In this situation we can construct a
parametrization $f\#g:S_{1}\#S_{2}\rightarrow M$ of
$\Sigma_{1}\#\Sigma_{2}$ ''coming'' from $f,g$, where the surface
$S_{1}\#S_{2}$ is the result of identify $S_{1}-\delta_{1}$ and
$S_{2}-\delta_{2}$ along the boundary of $\delta_{1}$ and
$\delta_{2}$. We can assume also that the immersion $f\#g$ agrees
with $f$ over $S_{1}-\delta_{1}$ and that $f\#g$ agrees with $g$
over $S_{2}-\delta_{2}$.

Let $K_{1},K_{2},T$ be triangulations of $S_{1},S_{2},M$ that make
$f,g$ and $f\#g$ simplicial, and assume that $T$ shells every
region of $M-\Sigma_{1}$. Consider a $T$-growth $f^{\prime}$ of
$f$ as in Proposition
\ref{PROPshellableimpliesfillinghomotopictoinflated}, such that
$f^{\prime}$ agrees with $f$ in all of $S_{1}$ except in two small
disks. We can take these two disks of $S_{1}$ where $f^{\prime}$
''disagrees'' with $f$ such that their respective images by $f$
are far away from Figure 23 (that is, they do not affect the
piping between $\Sigma_{1}$ and $\Sigma_{2}$). In this situation,
we can consider also the ''piped immersion'' $f^{\prime}\#g:S_{1}
\#S_{2}\rightarrow M$ that agrees with $f^{\prime}$ in
$S_{1}-\delta_{1}$ and with $g$ in $S_{2}-\delta_{2}$, as the
result of pasting $f^{\prime}$ and $g$ by means of the piping in
exactly the same way as $f$ was pasted with $g$ in $f\#g$.

\begin{remark}
\label{RemarkexistPQ}if $\Sigma_{1},\Sigma_{2}$ are the surfaces $\Sigma
_{1}^{\prime},\Sigma_{2}$ that result from the proof of the previous
Proposition, then there always exists such pair of points $P,Q$ as considered here.
\end{remark}

We know that $f^{\prime}$ is filling homotopic to $f$ because it is a
$T$-growth of $f$, but we also have:

\begin{lemma}
[Key Lemma 3]\label{Key Lemma3}If $\Sigma_{2}$ shells $\Sigma_{1}$, then we
can choose $f^{\prime}$ such that $f^{\prime}\#g$ is a $T$-growth of $f\#g$.
\end{lemma}

In particular, if $\Sigma_{2}$ shells $\Sigma_{1}$, then we can choose
$f^{\prime}$ such that $f^{\prime}\#g$ is filling homotopic to $f^{\prime}$.

This Lemma is also proved in detail in \cite{RHomotopies}. The
required property that $\Sigma_{2}$ shells $\Sigma_{1}$ implies
that the growing of $f$ into $f^{\prime}$ in the proof of
Proposition \ref{PROPshellableimpliesfillinghomotopictoinflated}
can be adapted to $\Sigma_{2}$ in such a way that \textit{the
growing from }$f$\textit{\ into }$f^{\prime}$\textit{\ defines
simultaneously a growing from }$f\#g$ \textit{\ into}
$f^{\prime}\#g$ when we introduce $\Sigma_{2}$.

\section{Proof of Theorem \ref{MAINtheorem}.\label{SECTION ProofofMainTheorem}}

With these tools we can make a sketch of the proof of Theorem
\ref{MAINtheorem}.

\begin{proof}
[Sketch of the Proof of Theorem \ref{MAINtheorem}]We are given a pair
$\Sigma_{1},\Sigma_{2}$ of nulhomotopic filling Dehn spheres of $M$ and two
parametrizations $f,g$ of them respectively. We will introduce the following
notation: we take two different 2-spheres $S_{1},S_{2}$ and we will consider
that $S_{i}$ is the domain of $\Sigma_{1}$ for $i=1,2$. In particular, it is
$\Sigma_{1}=f(S_{1})$ and $\Sigma_{2}=g(S_{2})$.

Modifying $f$ if necessary by an ambient isotopy of $M$ we can assume that
$\Sigma_{1}$ and $\Sigma_{2}$ have nonempty transverse intersection.

By Proposition \ref{PROPfillinghomotopictoMutuallyShellable} and
Remark \ref{RemarkexistPQ}, we can assume that
$\Sigma_{1},\Sigma_{2}$ form a mutually shellable filling pair of
spheres of $M$ and that there are two points $P,Q$ of
$\Sigma_{1}\cup\Sigma_{2}$ where $\Sigma_{1}$ and $\Sigma_{2}$
intersect as in Figure \ref{fig23a} in section \ref{SECTION
Growing}.

\begin{figure}[htb]
\centering \subfigure []{ \label{fig24a}
\includegraphics[ width=0.45\textwidth ] {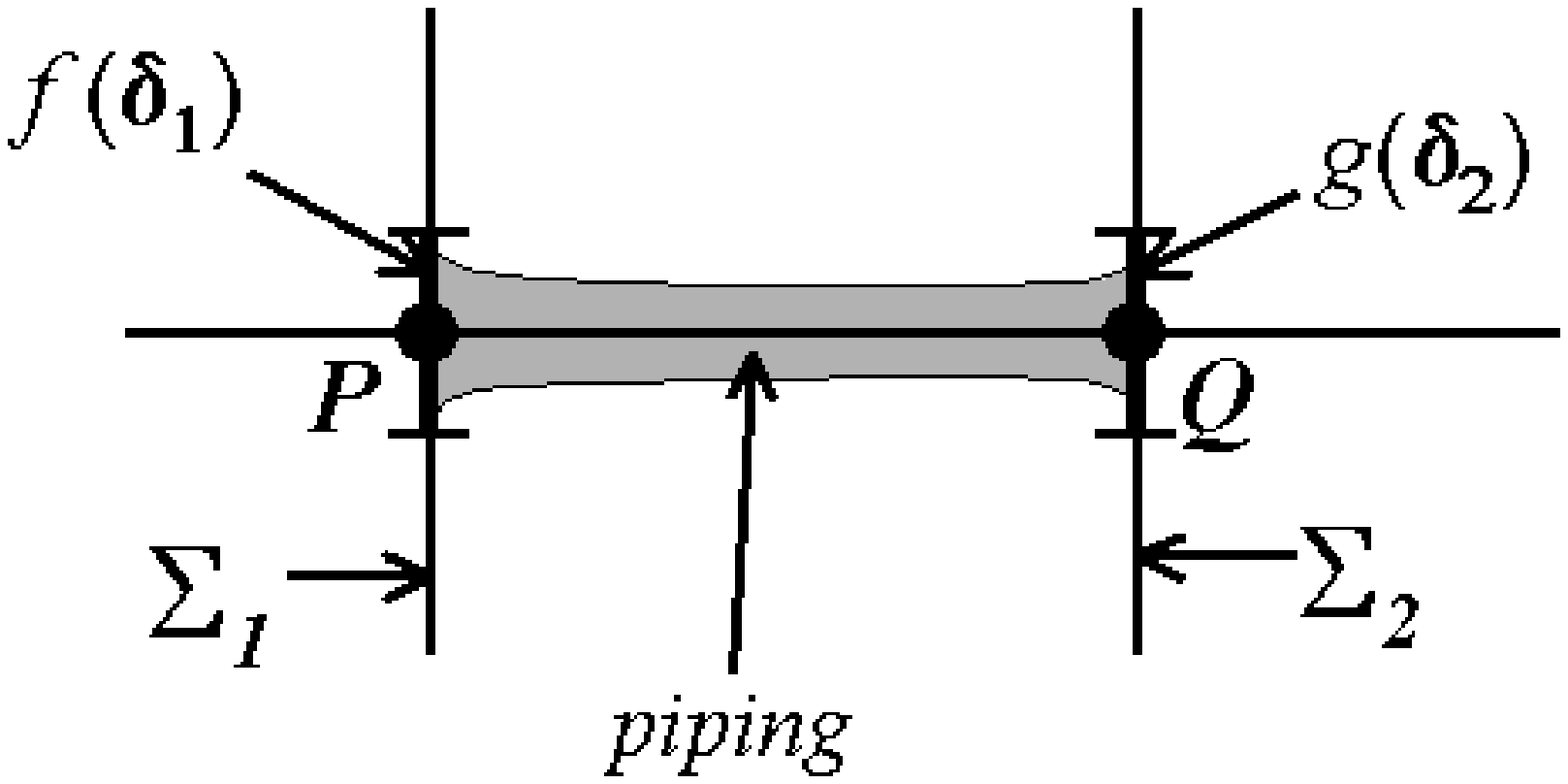}}\hfill
\subfigure {
\includegraphics[ height=0.15\textheight, width=0.05\textwidth ] {fig35ab} }\hfill
\addtocounter {subfigure}{-1} \subfigure[]{ \label{fig24b}
\includegraphics[ width=0.45\textwidth ] {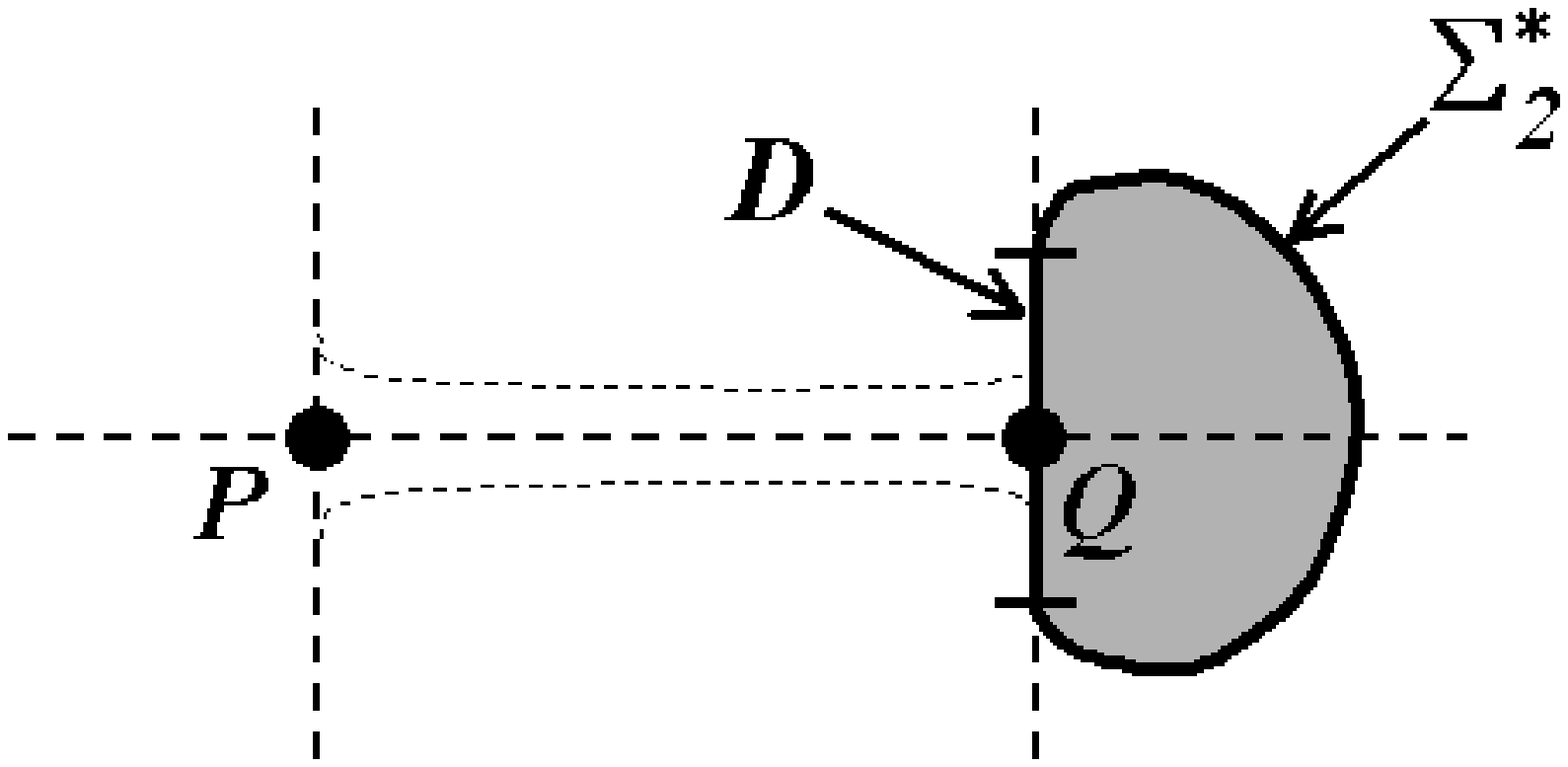}} \caption { }
\end{figure}

Consider the filling Dehn surface $\Sigma_{1}\#\Sigma_{2}$ and the
parametrization $f\#g$ as in section \ref{SECTION Growing}. We
consider also the disks $\delta_{1}\subset S_{1}$ and
$\delta_{2}\subset S_{2}$ as in section \ref{SECTION Growing}. We
denote $f\#g$ by $h$ for simplicity. $\medskip$

$\medskip h$\textsl{\ is filling homotopic to }$f$\textsl{.}

$\medskip$Consider a small standardly embedded 2-sphere
$\Sigma_{2}^{\ast} $ and a parametrization
$g_{\ast}:S_{2}\rightarrow M$ of $\Sigma_{2}^{\ast}$ as in Figure
\ref{fig24b}. This sphere shares with $\Sigma_{2}$ a 2-disk
$\bar{D}$ containing $g(\delta_{2})$ in its interior, and the
immersions $g$ and $g_{\ast}$ agree over
$\tilde{D}:=g_{\ast}^{-1}(\bar{D})$.

By the Key Lemma 1, we can deform $g$ into $g_{\ast}$ by a finite
sequence of transverse pushing disks leaving $\tilde{D}$ fixed.
Let $(D_{1},B_{1} ),...,(D_{k},B_{k})$ be this sequence of pushing
disks, and let $g=g_{0} ,g_{1},...,g_{k}=g_{\ast}:S_{2}\rightarrow
M$ be the sequence of transverse immersions such that $g_{i}$ is
obtained from $g_{i-1}$ by the pushing disk $(D_{i},B_{i})$.

\begin{figure}[htb]
\centering \subfigure []{ \label{fig25a}
\includegraphics[ width=0.45\textwidth ] {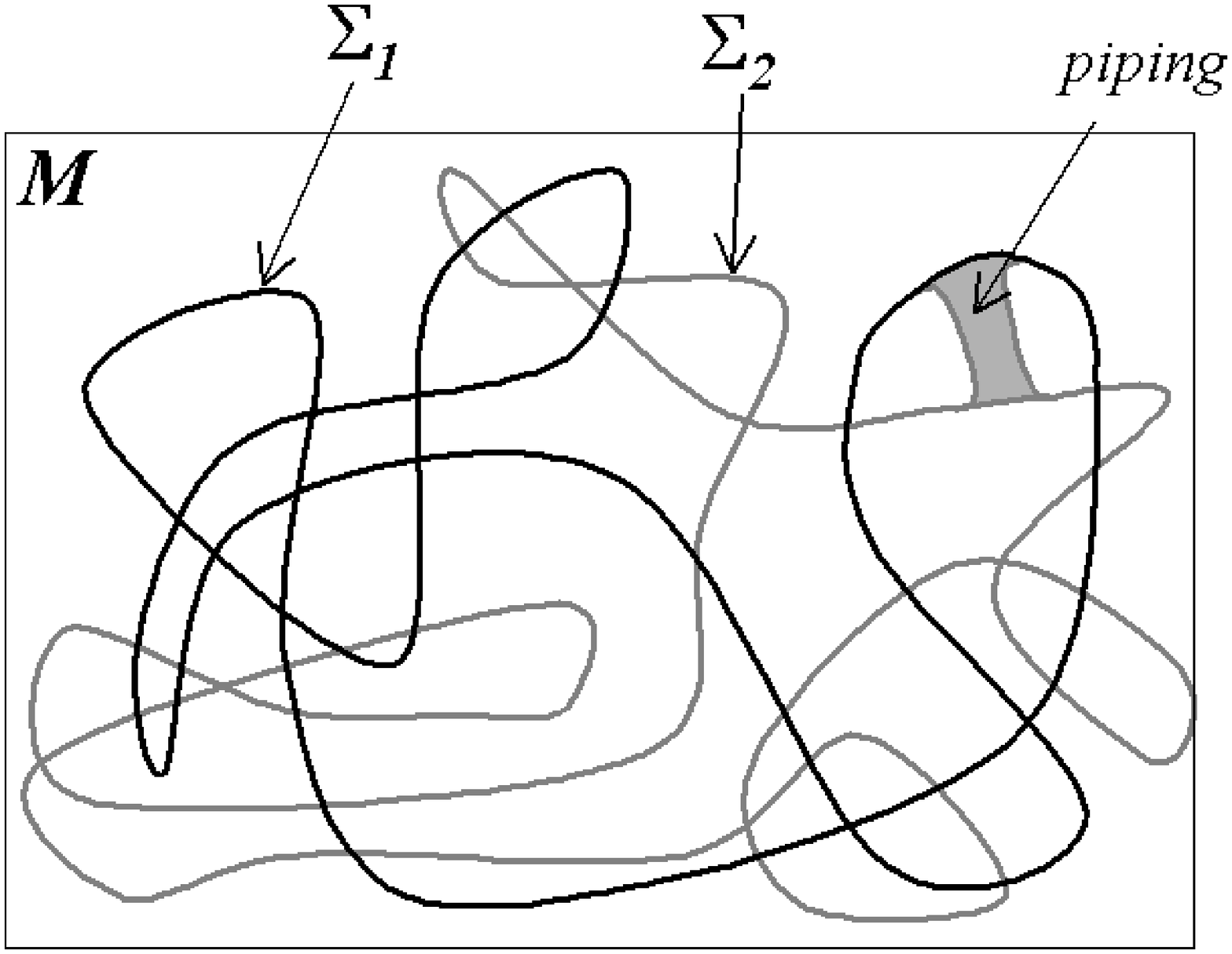}}\hfill \subfigure[]{
\label{fig25b}
\includegraphics[ width=0.45\textwidth ]{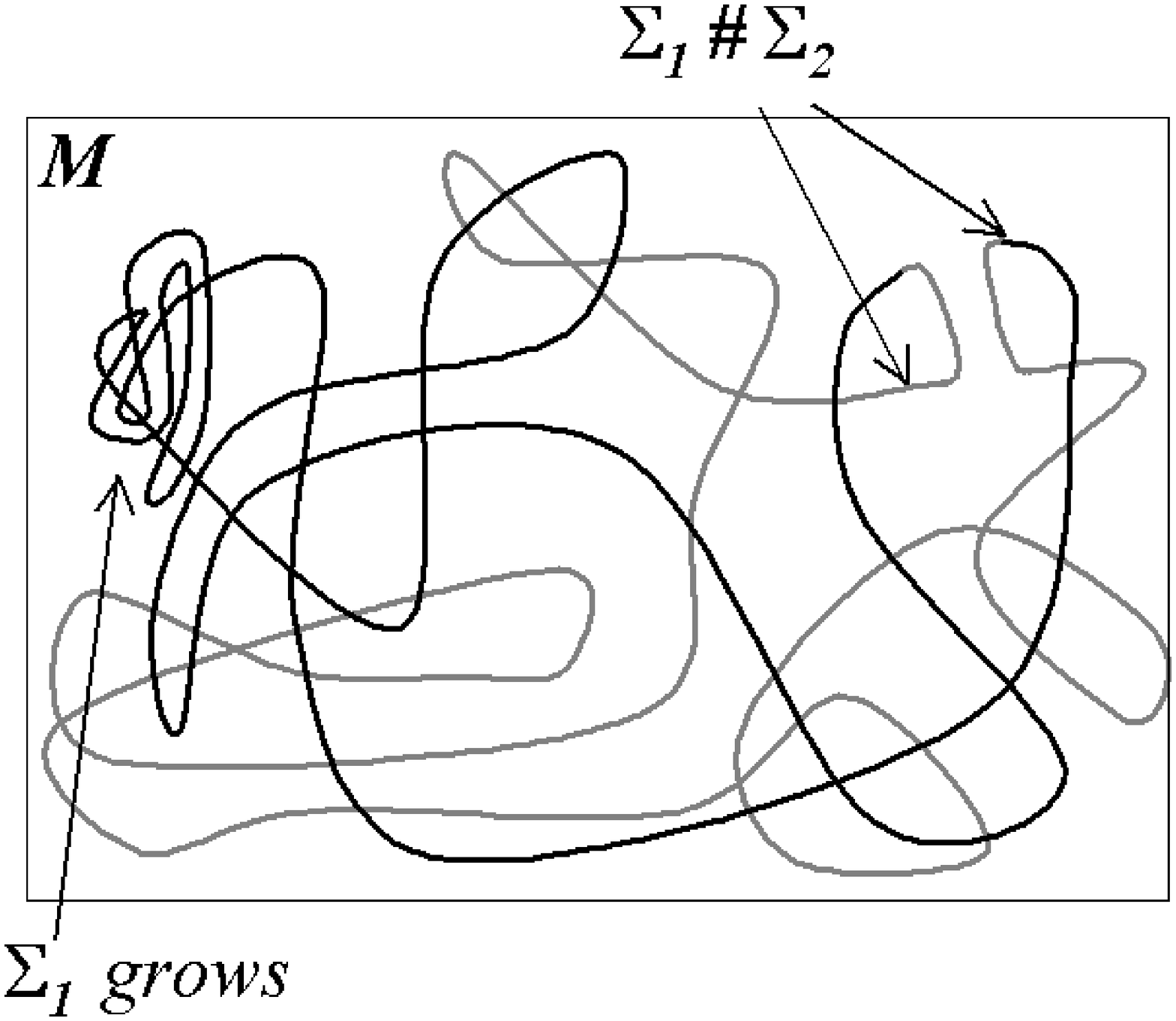}} \caption { }
\end{figure}

Modifying slightly $\Sigma_{1}$ and the piping between $\Sigma_{1}$ and
$\Sigma_{2}$ by an ambient isotopy of $M$ if necessary, we can assume that all
these pushing disks are transverse also with $\Sigma_{1}$ and $\Sigma
_{1}\#\Sigma_{2}$. Because the pushing disks $(D_{i},B_{i})$ leave $\tilde{D}$
fixed, we can think of these pushing disks as acting over the immersion
$h=f\#g$ instead of over $g$, and we can consider the sequence of transverse
immersions $h=h_{0},h_{1},...,h_{k}:S_{1}\#S_{2}\rightarrow M$ such that
$h_{i}$ is obtained from $h_{i-1}$ by the pushing disk $(D_{i},B_{i})$. Note
that $h_{k}(S_{1}\#S_{2})=\Sigma_{1}\#\Sigma_{2}^{\ast}$, where $\Sigma
_{1}\#\Sigma_{2}^{\ast}$ is obtained piping $\Sigma_{1}$ with $\Sigma
_{2}^{\ast}$ exactly in the same way as $\Sigma_{1}$ is piped with $\Sigma
_{2}$, and then a final transverse pushing disk $(D_{k+1},B_{k+1})$ transforms
$h_{k}$ into $f$. Thus, we can assume that \textit{there is a finite sequence
of transverse pushing disks leaving }$S_{1}-\delta_{1}$\textit{\ fixed\ that
transform }$h=f\#g$\textit{\ into }$f$.

\begin{figure}[htb]
\centering \subfigure []{ \label{fig26a}
\includegraphics[ width=0.45\textwidth ] {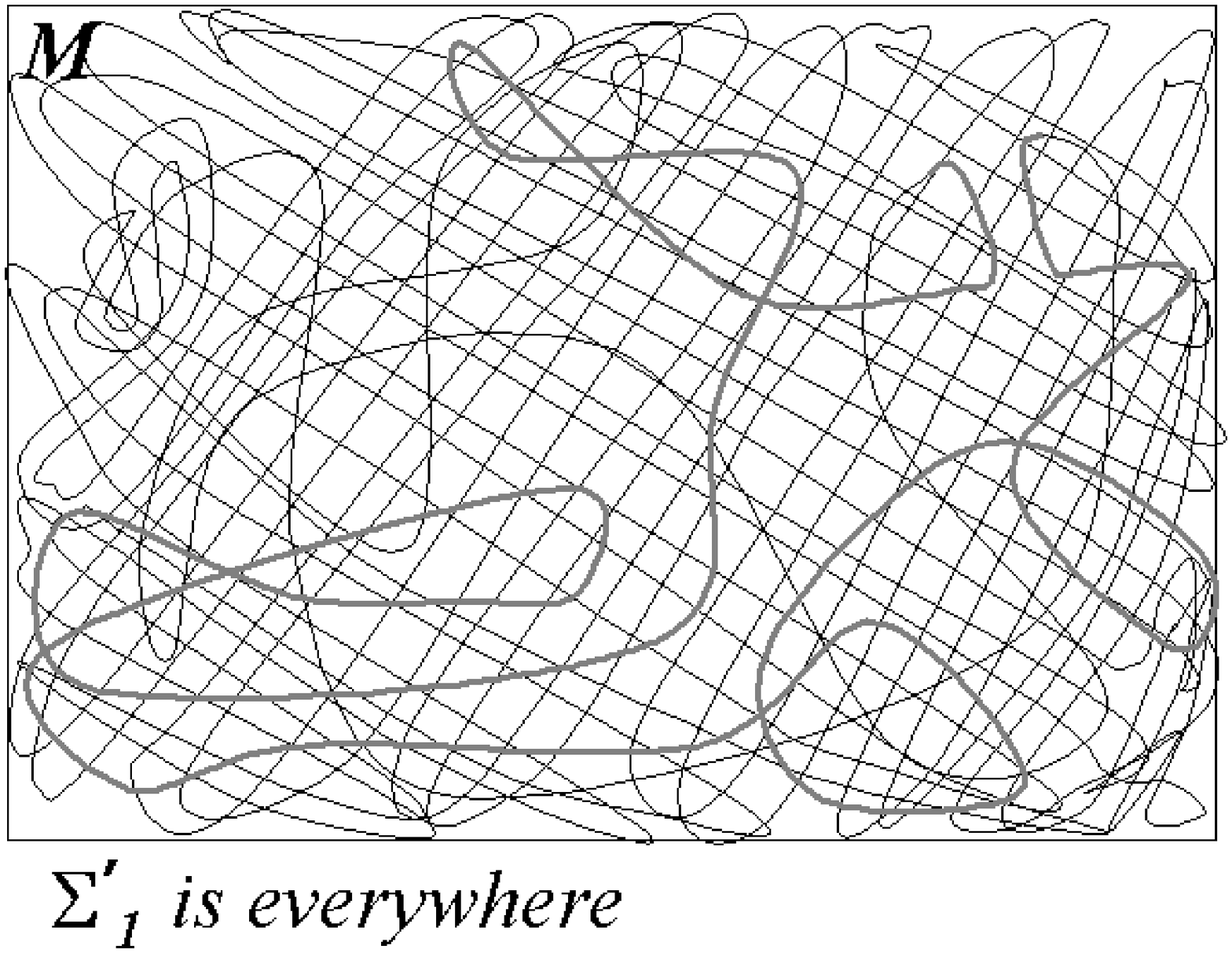}}\hfill \subfigure[]{
\label{fig26b}
\includegraphics[ width=0.45\textwidth ]{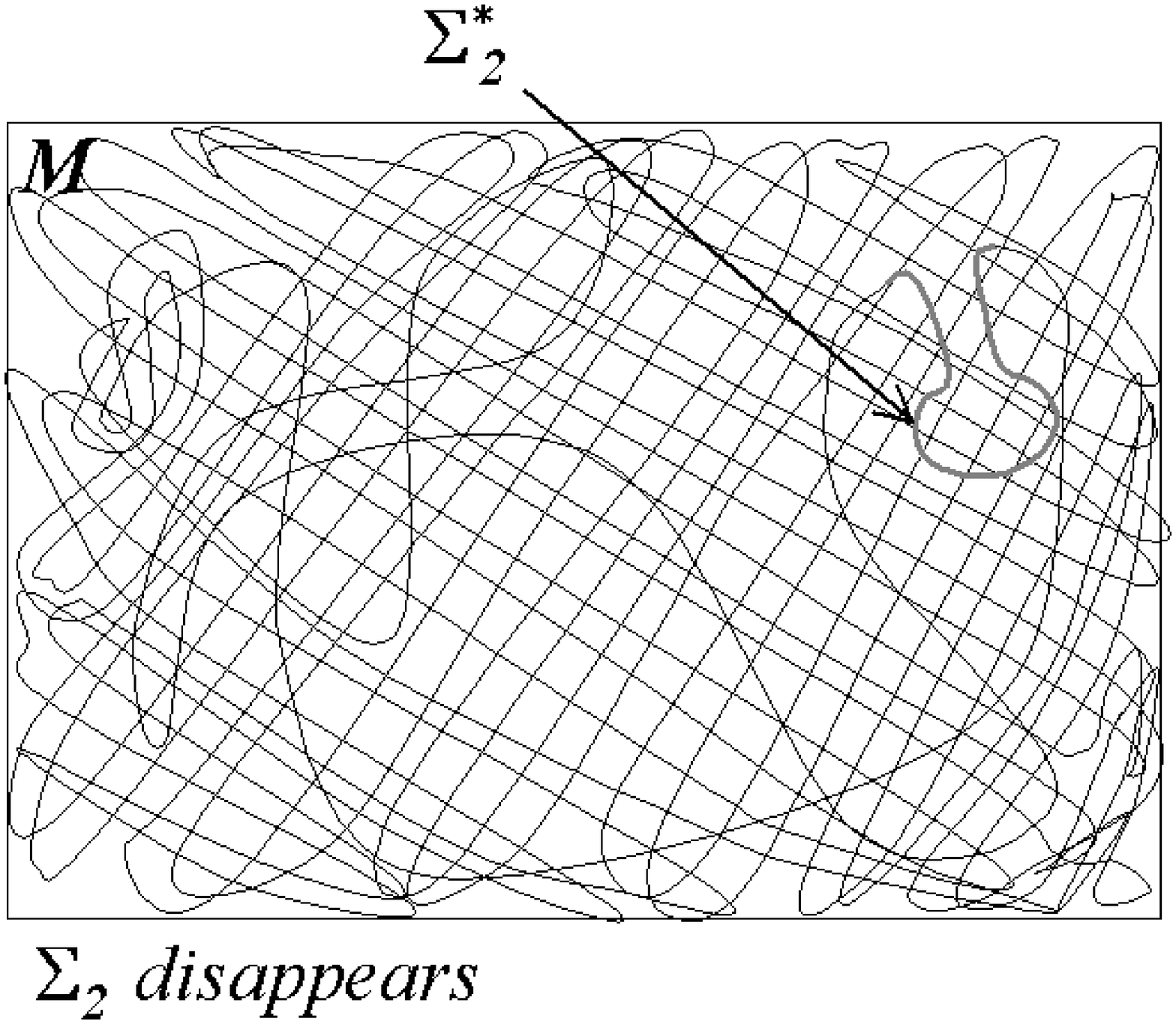}} \caption { }
\end{figure}

Take a good triangulation $T$ of $M$ with respect to
$f,g,g_{1},...,g_{k} ,h,h_{1},...,h_{k}$ (see Theorem
\ref{THMsmoothtriangulations exist} and Definition \ref{DEF
GoodTriangulations}).

The triangulation $T$ shells $f$ because $f$ is a filling immersion, so
consider a $T$-growth $f^{\prime}$ of $f$ such that the pipings of $\Sigma
_{1}$ with the components of $T$ do not affect $\Sigma_{2}$ neither the piping
between $\Sigma_{1}$ and $\Sigma_{2}$. Because $\Sigma_{1}$ and $\Sigma_{2}$
form a filling pair, in particular we have that $\Sigma_{2}$ shells
$\Sigma_{1}$. By Key Lemma 3, we can take $f^{\prime}$ such that it defines
also a $T$-growth $h^{\prime}:=f^{\prime}\#g$ of $h=f\#g$.

Consider now the sequence of immersions $h^{\prime}=h_{0}^{\prime}
,h_{1}^{\prime},...,h_{k}^{\prime},h_{k+1}^{\prime}=f^{\prime}$
such that $h_{i}^{\prime}$ is obtained from $h_{i-1}^{\prime}$ by
the pushing disk $(D_{i},B_{i})$, for $i=1,...,k+1$.

Note that by construction, each $h_{i}^{\prime}$ is a
$T$-inflating of $h_{i} $. Because of the choosing of $T$, for
each $i=1,...,k+1$ the triangulation $T$ restricted to the pushing
ball $B_{i}$ collapses simplicially into
$g(D_{i})=h^{\prime}(D_{i})$. By Key Lemma 2 this implies that
$h_{i}^{\prime }$ is filling homotopic to $h_{i-1}^{\prime}$ for
each $i=1,...,k+1$.

\begin{figure}[htb]
\centering \subfigure []{ \label{fig27a}
\includegraphics[ width=0.45\textwidth ] {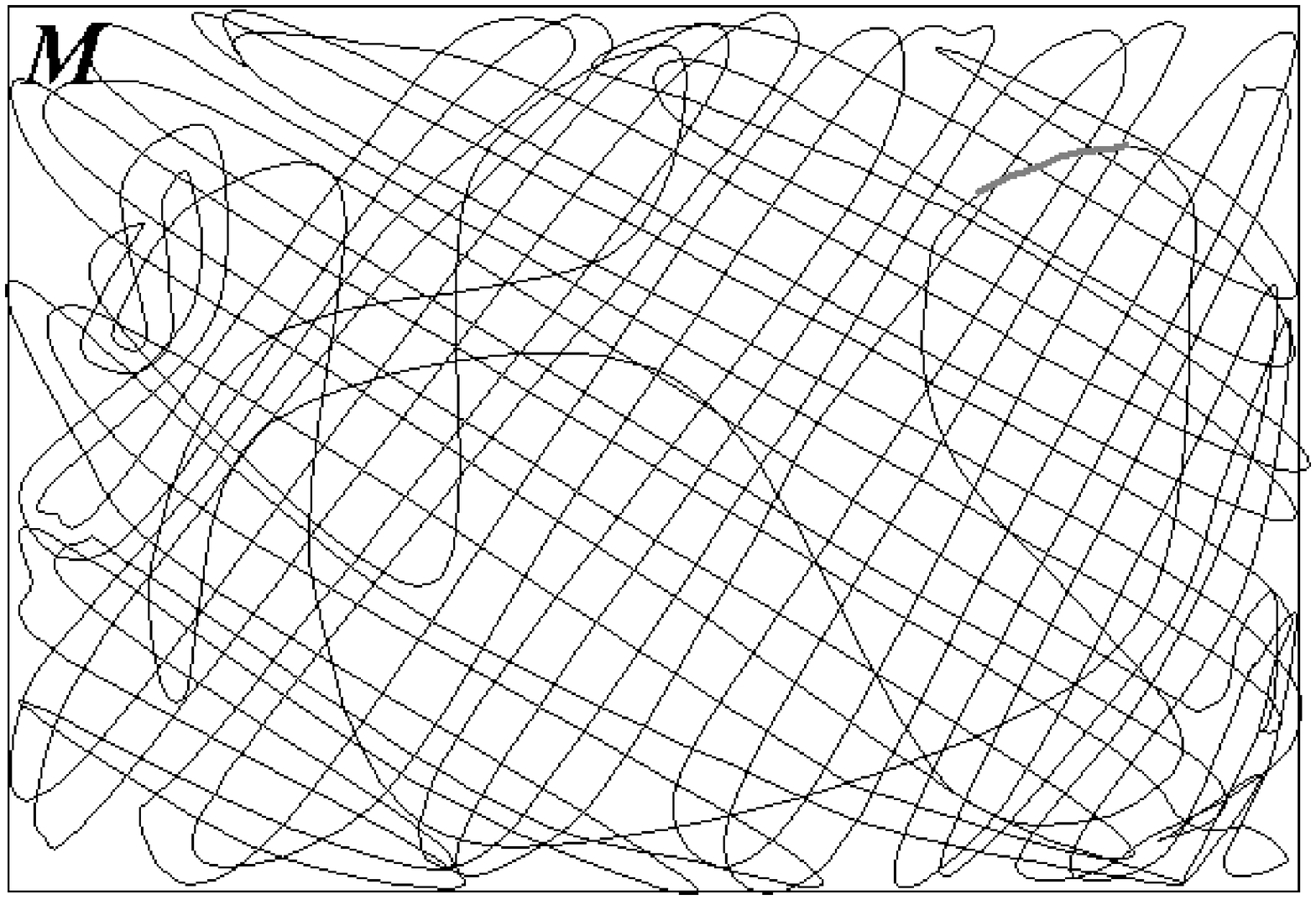}}\hfill \subfigure[]{
\label{fig27b}
\includegraphics[ width=0.45\textwidth ]{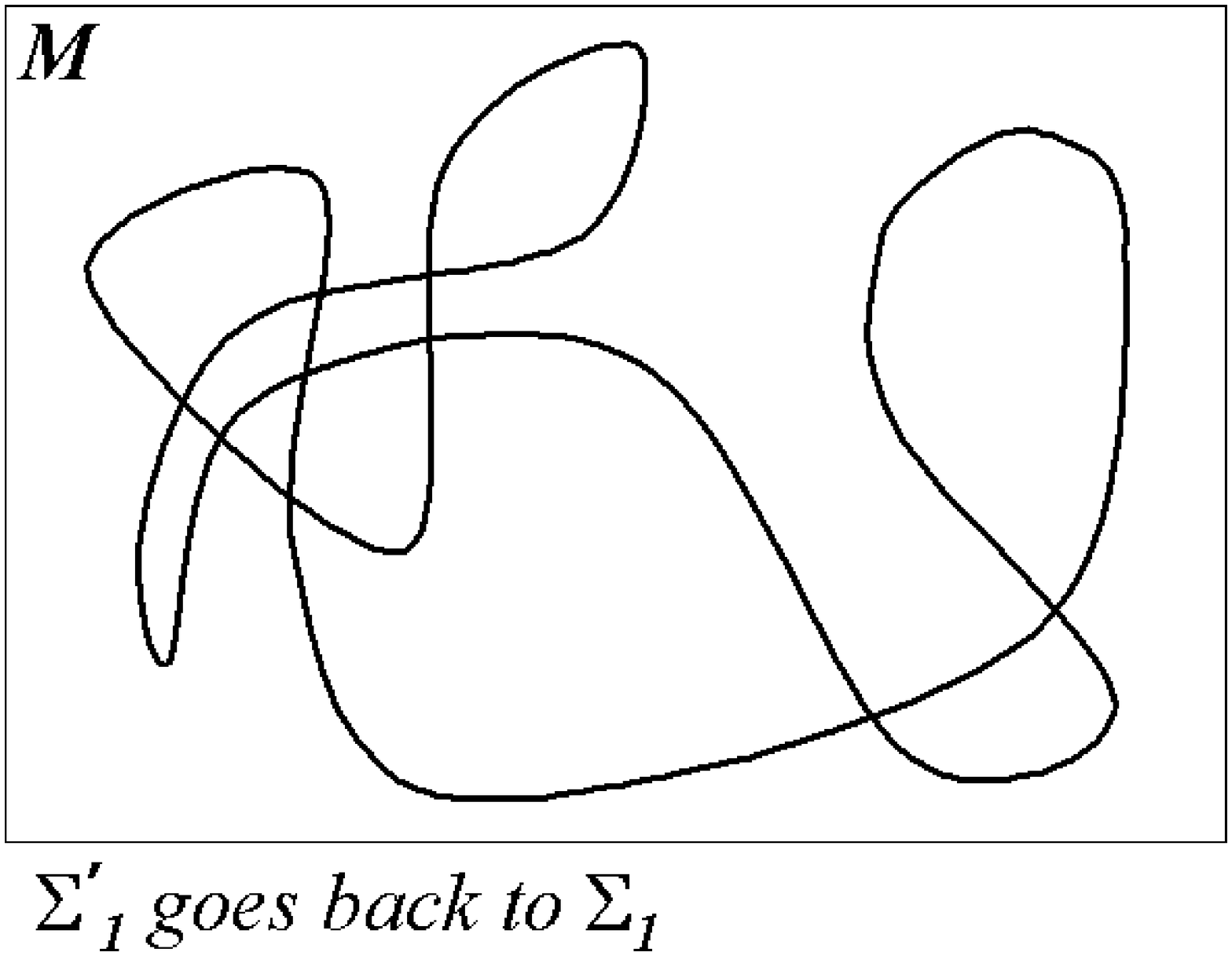}} \caption { }
\end{figure}

Resuming, we have that $h$ is filling homotopic to $h^{\prime}$ because
$h^{\prime}$ is a $T$-growth of $h$, $h^{\prime}$ is filling homotopic to
$f^{\prime}=h_{k+1}^{\prime}$ by repeatedly applications of Key Lemma 2 to the
pushing disks $(D_{i},B_{i})$ for $i=1,...,k+1$, and $f^{\prime}$ is filling
homotopic to $f$ because $f^{\prime}$ is a $T$-growth of $f$. Therefore, $h$
is filling homotopic to $f$.

By the same arguments, we have that $h$ is filling homotopic to
$g$ and thus $f$ and $g$ are filling homotopic.

\end{proof}

\section{Diagrams.\label{SECTION Diagrams}}

Let $\Sigma\subset M$ be a Dehn surface in $M$, and let $f:S\rightarrow M$ be
a parametrization of $\Sigma$. As we have pointed out in section \ref{SECTION
Introduction}, the singular set $S(f)$, together with the information about
how its points are identified by $f$ is what we call the Johansson diagram of
$\Sigma$. We will give a more detailed definition of Johansson diagram. This
new definition is equivalent to the definition given in \cite{Papa}. We assume
now by simplicity that both $S$ and $M$ are orientable.

Let $\bar{\gamma}:S^{1}\rightarrow M$ be a parametrization of a
double curve of $\Sigma$. Because both $S$ and $M$ are orientable,
the inverse image by $f$ of $\bar{\gamma}(S^{1})$ is the union of
two different closed curves in $S(f)$. There are exactly two
different immersions $\gamma_{1},\gamma _{2}:S^{1}\rightarrow S$
such that $f\circ\gamma_{1}=f\circ\gamma_{2} =\bar{\gamma}$. In
this situation, we say that $\gamma_{1}$ and $\gamma_{2}$ are
\textit{lifted curves} of $\bar{\gamma}$ under $f$ and that they
are \textit{sisters} under $f$.

A \textit{complete parametrization} of the singularity set
$\bar{S}(\Sigma)$ of $\Sigma$ is a set $\mathcal{\bar{D}}=\left\{
\bar{\alpha}_{1},\bar{\alpha }_{2},...,\bar{\alpha}_{m}\right\}  $
of immersions from $S^{1}$ into $M$ such that: (i) each
$\bar{\alpha}_{i}$ parametrizes a double curve of $\Sigma$; (ii)
$\bar{\alpha}_{i}(S^{1})\neq\bar{\alpha}_{j}(S^{1})$ if $i\neq j$;
and (iii)
$\bar{S}(\Sigma)=\underset{i=1}{\overset{m}{\cup}}\bar{\alpha}_{i}
(S^{1})$. If $\mathcal{\bar{D}}$ is a complete parametrization of
$\bar {S}(\Sigma)$ and we denote by $\mathcal{D}$ the set of all
lifted curves of the curves of $\mathcal{\bar{D}}$, the map
$\tau:\mathcal{D}\rightarrow \mathcal{D}$ that assigns to each
curve of $\mathcal{D}$ its sister curve under $f$ defines a free
involution of $\mathcal{D}$. The pair $(\mathcal{D} ,\tau)$
contains all the information about the singular set $S(f)$ and
about how the points of $S(f)$ are identified by the map $f$: two
different points $A,B\in S$ verify $f(A)=f(B)$ if and only if
there is a parametrized curve $\alpha\in\mathcal{D}$ and a $z\in
S^{1}$ with $A=\alpha(z)$ and $B=\tau \alpha(z)$.

The pair $(\mathcal{D},\tau)$ of the previous paragraph is the model we will
presently use for defining (\textit{abstract}) \textit{diagram}. We have seen
that every Dehn surface has an associated Johansson diagram. Thus, we can
define an \textit{abstract diagram} in a surface $S$ as a collection
$\mathcal{D}$ of closed curves in $S$ together with a free involution
$\tau:\mathcal{D}\rightarrow\mathcal{D}$, such that the curves of
$\mathcal{D}$ can be coherently identified by $\tau$. The natural question now
is if this diagram $(\mathcal{D},\tau)$ is the Johansson diagram coming from a
transverse immersion $f:S\rightarrow M$ of $S$ into some orientable 3-manifold
$M$. If this occurs, we say that the (abstract) diagram $(\mathcal{D},\tau)$
in the surface $S$ is \textit{realizable }(compare \cite{Papa}) and that the
immersion $f$ \textit{realizes} the diagram $(\mathcal{D},\tau)$.

\begin{figure}[htb]
\centering \subfigure []{ \label{fig28a}
\includegraphics[ height=0.2\textheight ] {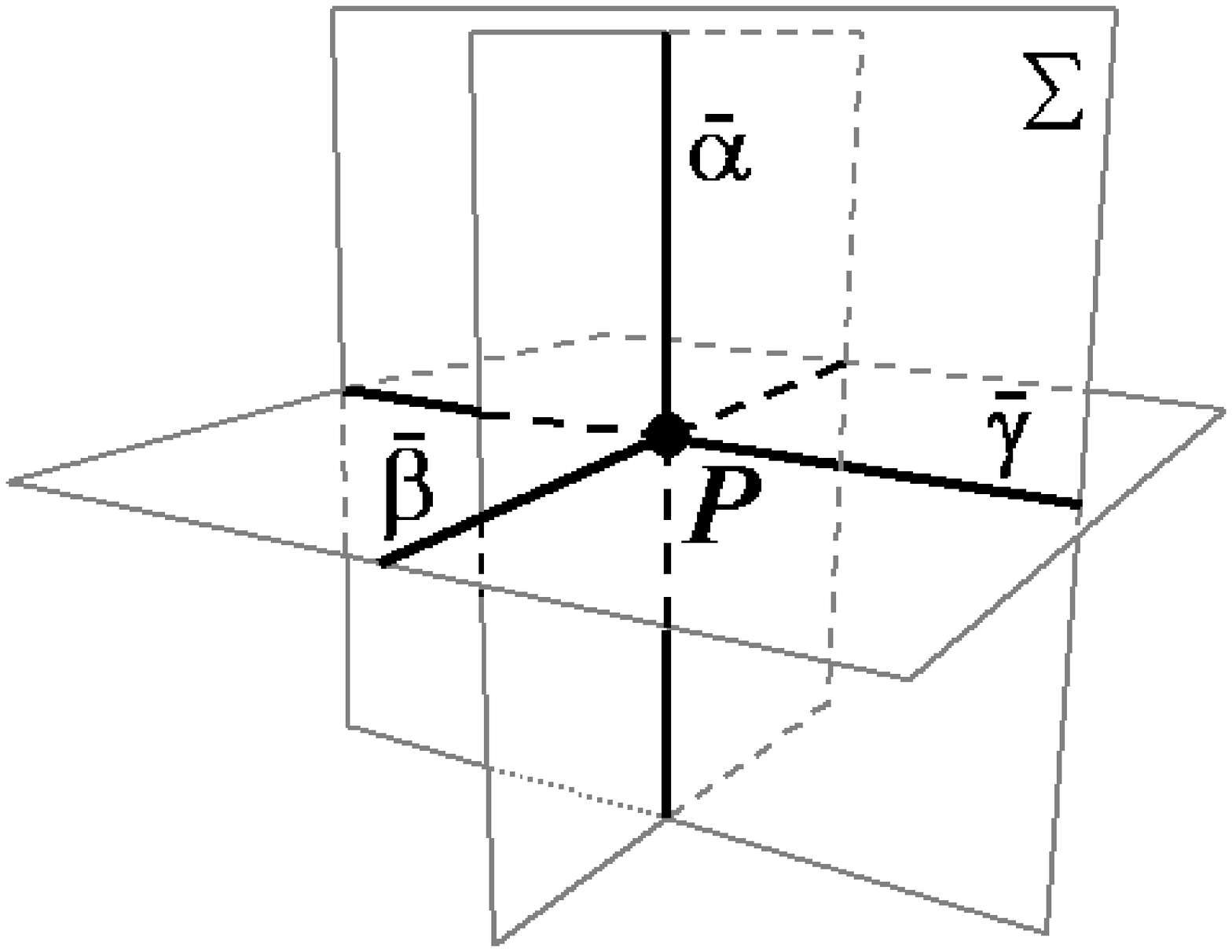}}\hspace {0.6cm}
\subfigure {
\includegraphics[ height=0.2\textheight ] {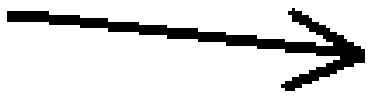} }\hspace {0.6cm} \addtocounter {subfigure}{-1}
\subfigure[]{ \label{fig28b}
\includegraphics[ height=0.2\textheight ] {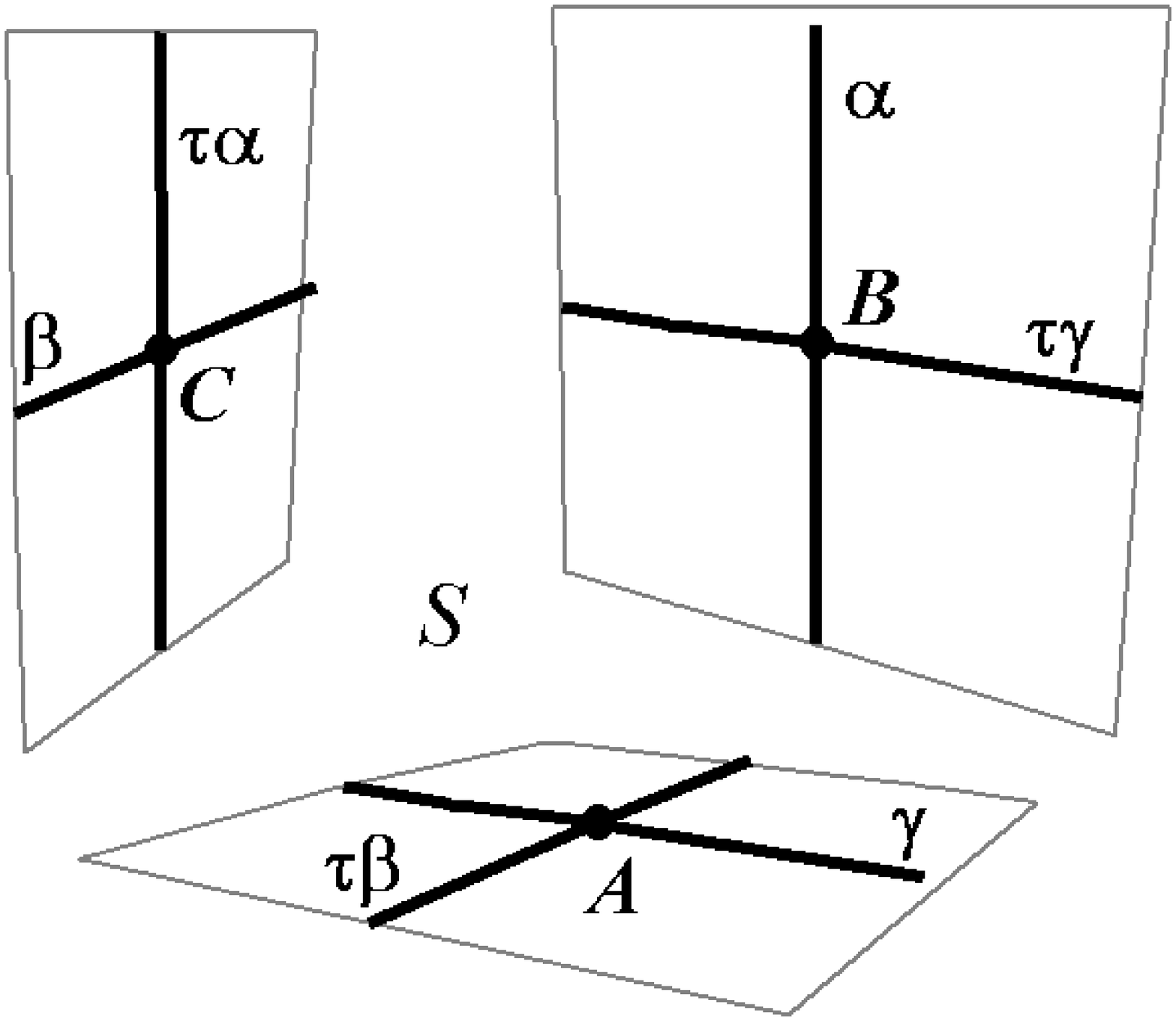}} \caption { }
\end{figure}

The first condition that must verify the curves of $\mathcal{D}$
is that they intersect transversely as in Figure \ref{fig28b} at
some points of $S$ which are the \textit{double points} of the
diagram $(\mathcal{D},\tau)$. We define that two different points
$A,B\in S$ are \textit{related} by the diagram $(\mathcal{D}
,\tau)$ if there is a curve $\alpha\in\mathcal{D}$ and a $z\in
S^{1}$ with $A=\alpha(z)$ and $B=\tau\alpha(z)$. With this
notation, each double point $A$ of the diagram will be related
with two points $B,C$ of the surface $S$. If $(\mathcal{D},\tau)$
is realizable, these two points $B,C$ must be different and they
must be also related by the diagram (see Figure \ref{fig28b}).
Thus, the double points of the diagram must be arranged in
\textit{triplets} of pairwise related points (the diagram is
\textit{riveted} in the notation of \cite{Carter}). If
$f:S\rightarrow M$ realizes the diagram, each of these triplets
compose the inverse image by $f$ of a triple point of $f$. In
Figure 29 we have given names to the double points of the diagram
in such a way that related points are denoted with the same name.
We consider two diagrams on $S$ as equivalent if they are related
by an homeomorphism of $S$ or by reparametrization of the curves
of the diagram.

The following is a survey of the main result of \cite{Johansson2}
about the realizability of diagrams. We will denote the diagram
$(\mathcal{D},\tau)$ simply by $\mathcal{D}$.

\begin{figure}[htb]
\centering \subfigure []{ \label{fig29a}
\includegraphics[ height=0.2\textheight ] {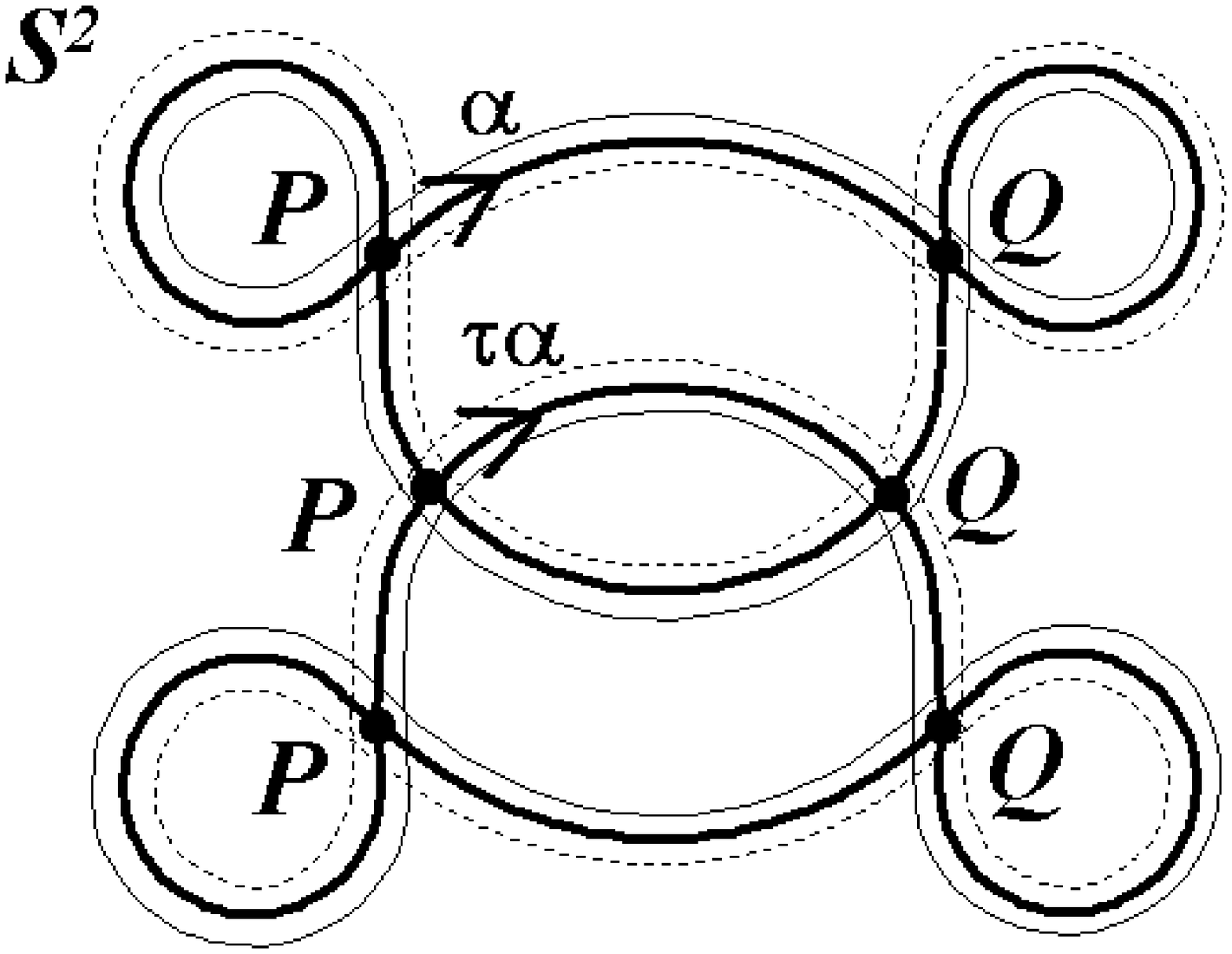}}\hfill \subfigure[]{ \label{fig29b}
\includegraphics[ height=0.2\textheight ] {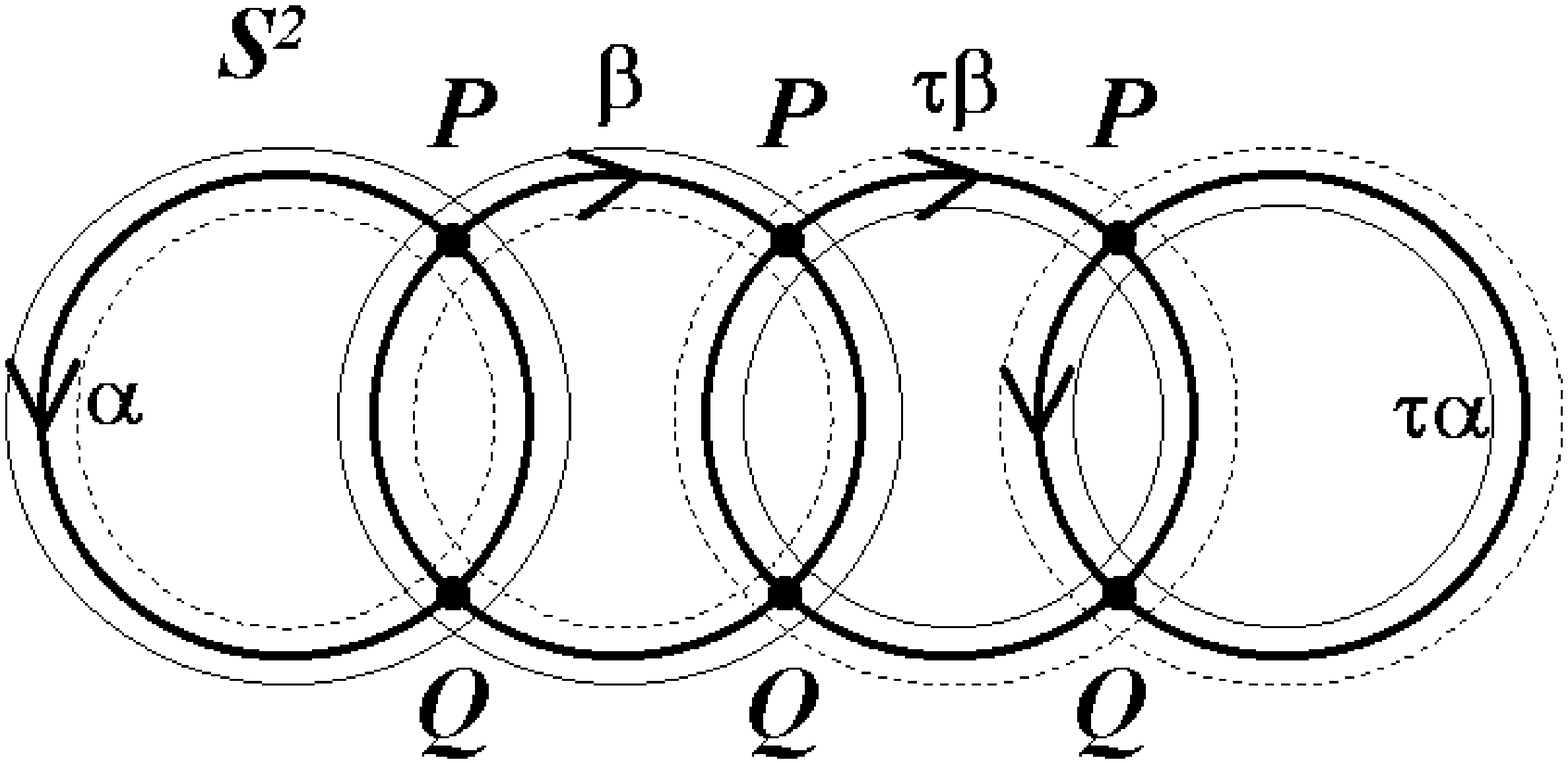}} \caption { }
\end{figure}

\begin{figure}[htb]
\centering
\includegraphics[ width=0.5\textwidth ]{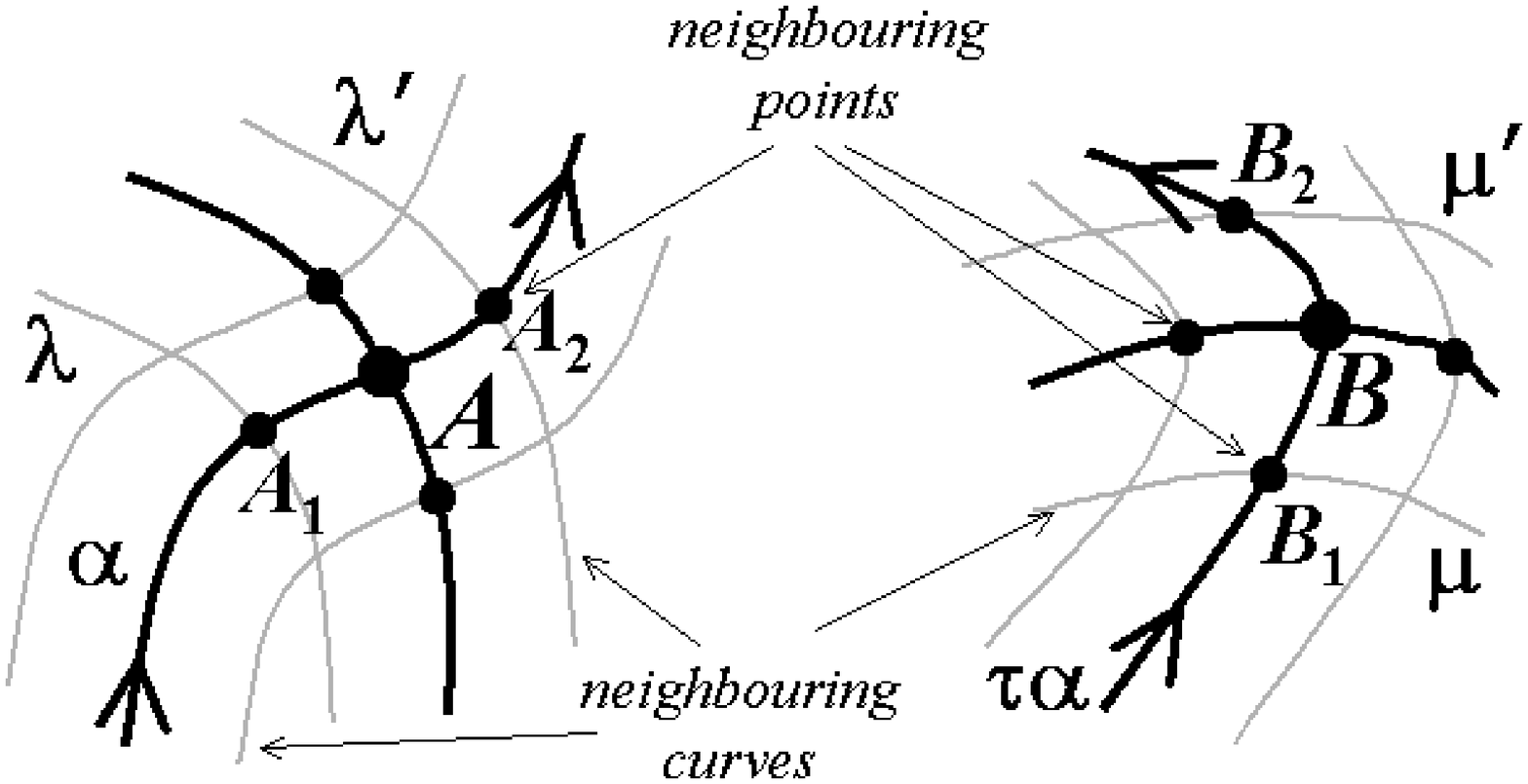} \caption { } \label{fig30}
\end{figure}

Assume that we are given an abstract diagram $\mathcal{D}$ on the
surface $S$ (Figure 29). For each $\alpha\in\mathcal{D}$, we
consider two \textit{neighbouring curves}
$\lambda,\lambda^{\prime}$ that run parallel to $\alpha$ and such
that $\lambda$ and $\lambda^{\prime}$ lie on different sides of
$\alpha$ (Figure 29). We say that the two neighbouring curves of
the same curve of the diagram are \textit{opposite} neighbouring
curves of the diagram. The neighbouring curves of the diagram can
be taken such that they only intersect near the double points of
the diagram and exactly as depiced in Figures 29 and \ref{fig30}.
We call \textit{neighbouring points} of the diagram to the
intersection points of the neighbouring curves with the curves of
the diagram. With these assumptions, there appear four
neighbouring points around each double point of the diagram
(Figure \ref{fig30}). Consider two related double points $A,B$ of
the diagram. Because they are related, there is a curve $\alpha
\in\mathcal{D}$ and a $z\in S^{1}$ with $A=\alpha(z)$ and
$B=\tau\alpha(z)$. If we orient the curves $\alpha,\tau\alpha$
using the standard orientation of $S^{1}$, then near $A$ the curve
$\alpha$ pass through the points $A_{1},A,A_{2}$ in this order,
where $A_{1},A_{2}$ are neighbouring points of the diagram. In the
same way, near $B$ the curve $\tau\alpha$ pass trough the points
$B_{1},B,B_{2}$ in this order, where $B_{1},B_{2}$ are
neighbouring points of the diagram. We assume that the
neighbouring curves are so chosen that in this situation $A_{i}$
is related by the diagram with $B_{i}$ for $i=1,2$ (see Figure
\ref{fig30}).

Once we have drawn the neighbouring curves of the diagram as in
the previous paragraph, we give some definitions. If two
neighbouring curves $\lambda,\mu$ pass through related
neighbouring points, as the curves $\lambda$ and $\mu$ of Figure
\ref{fig30}, we say that $\lambda,\mu$ are \textit{elementary
related}. If we orient all the curves of the diagram using the
standard orientation of $S^{1}$ and if we consider the surface $S$
oriented, for a curve $\alpha\in \mathcal{D}$ we say that the
neighbouring curve of $\alpha$ lying on the left-hand side of
$\alpha$ is \textit{elementary G-related} with the neighbouring
curve of $\tau\alpha$ lying on the right-hand side of $\tau
\alpha$, and equivalently, that the neighbouring curve of $\alpha$
lying on the right-hand side of $\alpha$ is elementary G-related
with the neighbouring curve of $\tau\alpha$ lying on the left-hand
side of $\tau\alpha$.

\begin{definition}
\label{DEF-G-Classes}\cite{Johansson2}Two neighbouring curves $\lambda,\mu$ of
the diagram $\mathcal{D}$ are in the same G-class if there exists a finite
sequence $\lambda=\lambda_{0},\lambda_{1},...,\lambda_{k}=\mu$ of neighbouring
curves of the diagram such that $\lambda_{i-1}$ is elementary related or
elementary G-related with $\lambda_{i}$ for $i=1,...,k $.
\end{definition}

In the diagrams of Figures \ref{fig29a} and \ref{fig29b} we have
drawn in the same way the neighbouring curves in the same G-class.
From the construction of G-classes can be checked whithout
difficulty the following

\begin{lemma}
\label{LeMMA G-Classes-singularity set}If $\mathcal{D}$ is a diagram in $S$
and $f:S\rightarrow M$ realizes $\mathcal{D}$, then the number of G-classes of
$\mathcal{D}$ is twice the number of connected components of the singularity
set $\bar{S}(f)$.
\end{lemma}

The following Theorem appears in \cite{Johansson2}.

\begin{theorem}
\label{THM Johansson G-Clases}\cite{Johansson2}A diagram $\mathcal{D}$ in the
orientable surface $S$ is realizable by a transverse immersion $f:S\rightarrow
M$ of $S$ into an orientable 3-manifold $M$ if and only if there are no
opposite neighbouring curves of the diagram in the same G-class.
\end{theorem}

This Theorem gives an easy method for checking realizability on a wide class
of diagrams. An analogue result is given in \cite{Carter} for diagrams with no
closed components in surfaces with boundary.

Though Theorem \ref{THM Johansson G-Clases} was stated in
\cite{Johansson2} for diagrams in the 2-disk without singular
boundary points, as it is pointed out in \cite{Papa} the proof can
be extended directly to the case stated here. More exactly,
Theorem \ref{THM Johansson G-Clases} is also true if we remove
from $S$ a finite number of open disks not touching the diagram.

\begin{figure}[htb]
\centering
\includegraphics[ width=0.45\textwidth ]{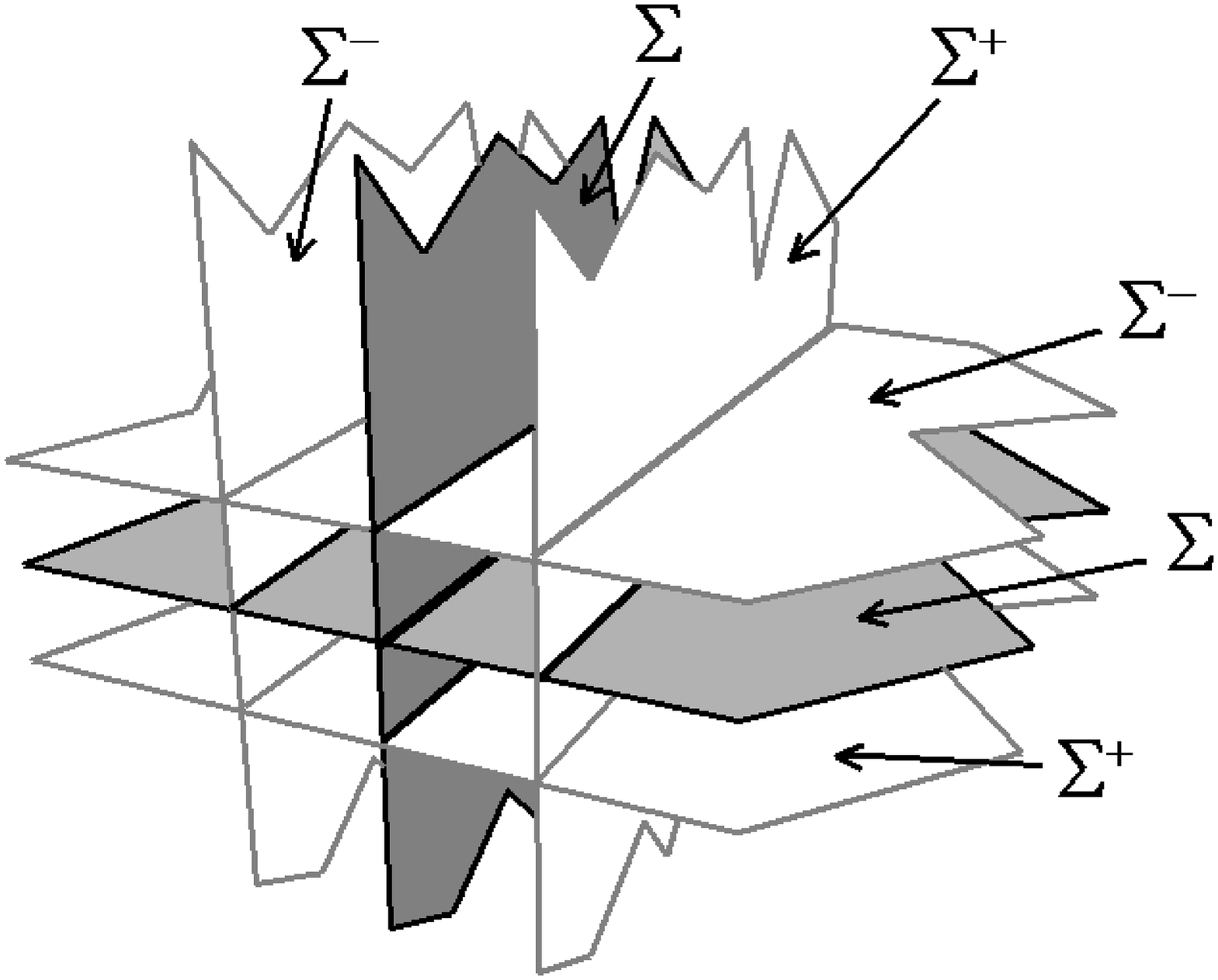} \caption { } \label{fig31}
\end{figure}

The key for proving Theorem \ref{THM Johansson G-Clases} is
\textit{2-sidedness}. Every immersion $f:S\rightarrow M$ with both
$S$ and $M$ orientable is \textit{2-sided}, and this 2-sidedness
is reflected in the neighbouring curves of the diagram. An
immersion $f:S\rightarrow M$ is 2-sided if there exists an
immersion $F:S\times\left[  -1,1\right]  \rightarrow M$ with
$F(X,0)=f(X)$ for every $X\in S$. Put $\Sigma=f(S)$. If $f$ is
2-sided and transverse, we can choose $F$ as close to $f$ as we
want, such that in a neighbourhood of a double curve of $\Sigma$
the image of $F$ looks like Figure \ref{fig31}. In that Figure, we
can see that the \textit{lower sheet} $\Sigma
^{-}=F(S\times\left\{ -1\right\}  )$ and the \textit{upper sheet}
$\Sigma ^{+}=F(S\times\left\{  1\right\}  )$ intersect $\Sigma$ in
some curves that behave exactly as the images by $f$ of the
neighbouring curves of the Johansson diagram $\mathcal{D}$ of $f$.
We can assume without loss of generality that in this case, the
neighbouring curves of the diagram compose exactly the inverse
image $f^{-1}(\Sigma^{-}\cup\Sigma^{+})$. In \cite{Johansson2} it
is proved the following Proposition.

\vspace {10pt}

\begin{proposition}
\label{PROPsameG-class-sameSheet}If two neighbouring curves of $\mathcal{D}$
are in the same G-class, then their images by $f$ must be contained in the
same sheet $\Sigma^{-}$ or $\Sigma^{+}$.
\end{proposition}

This implies that if $\mathcal{D}$ is realizable by a 2-sided immersion, there
cannot be two opposite neighbouring curves in the same G-class.

On the other hand, if there are not two opposite neighbouring
curves of the abstract diagram $\mathcal{D}$ on $S$ in the same
G-class, we can make an identification $\sim$ on the thickened
surface $S\times\left[  -1,1\right]  $ compatible with the diagram
such that neighbourhoods of sister curves are identified as in
Figure \ref{fig31}. The quotient
$\hat{M}(\mathcal{D})=S\times\left[ -1,1\right]  /\sim$ is a
3-manifold with boundary and it verifies: (i) the canonical
projection $\pi:S\times\left[  -1,1\right] \rightarrow\hat
{M}(\mathcal{D})$ is an immersion; (ii) if we take the inclusion
$j:S\rightarrow S\times\left[  -1,1\right]  $ given by
$j(X)=(X,0)$, then $\pi\circ j$ is a transverse immersion
realizing $\mathcal{D}$; and (iii) $\hat{M}(\mathcal{D}) $ is
orientable. See \cite{Johansson1} and \cite{Johansson2} for more
details.

Going back to the immersion $f$, if it is a filling immersion, the singularity
set $\bar{S}(f)$ must be connected and by Lemma \ref{LeMMA
G-Classes-singularity set} this implies that the Johansson diagram
$\mathcal{D}$ of $f$ has only two G-classes of neighbouring curves. In this
case, the manifold with boundary $\hat{M}(\mathcal{D})$ constructed from
$\mathcal{D}$ as in the previous paragraph \textit{is uniquely determined by
}$\mathcal{D}$ and it is homeomorphic to a regular neighbourhood of the
filling Dehn surface $\Sigma\subset M$. Because $f$ is a filling immersion the
boundary of $\hat{M}(\mathcal{D})$ must be composed by a union of 2-spheres.
Pasting a 3-ball to $\hat{M}(\mathcal{D})$ along each boundary component we
obtain a closed 3-manifold $M(\mathcal{D})$ homeomorphic to $M$. This is the
way for reconstructing a 3-manifold $M$ from a Johansson representation of
$M$.

\begin{figure}[htb]
\centering \subfigure []{ \label{fig32a}
\includegraphics[ width=0.4\textwidth ] {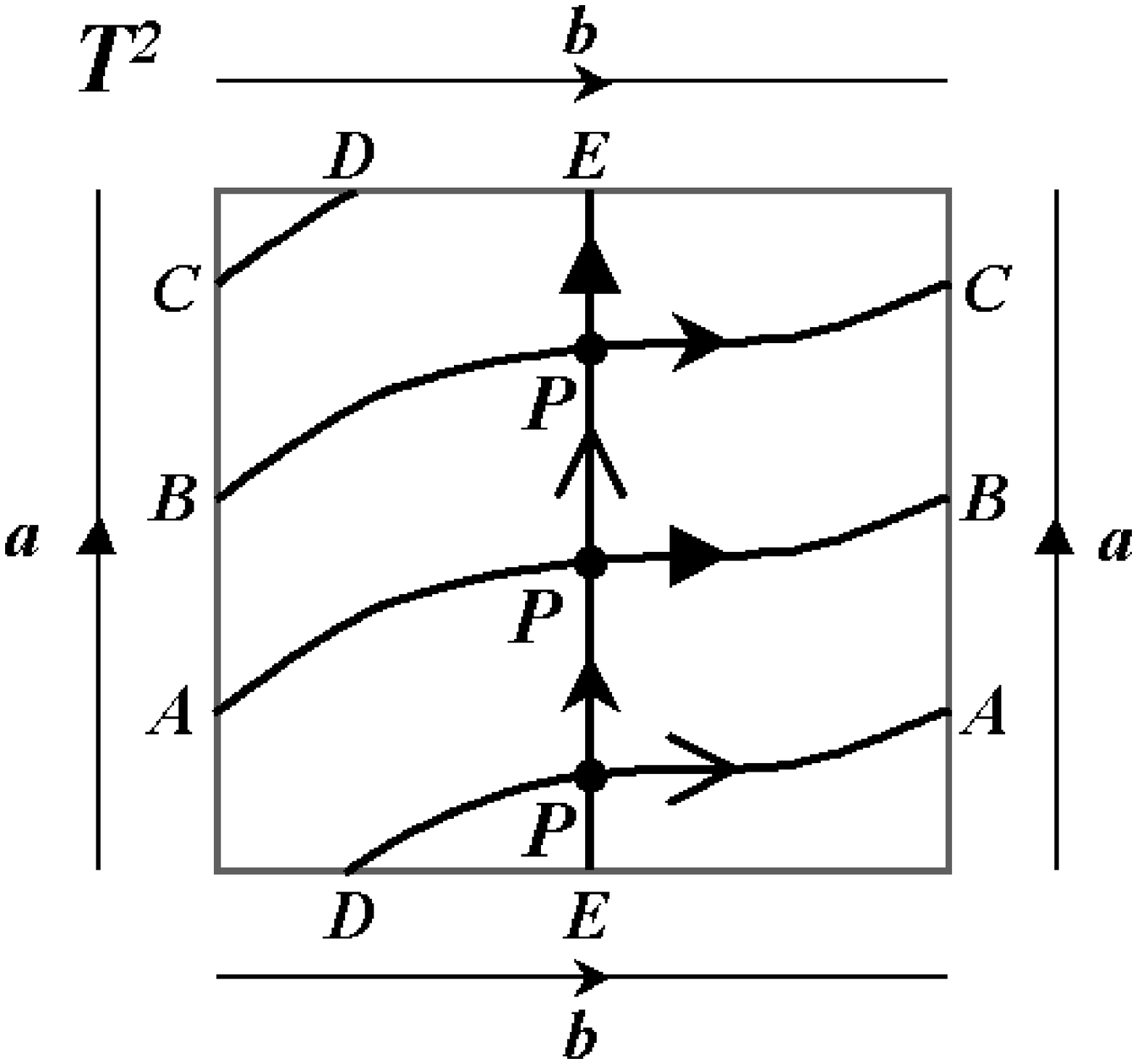}}\hfill \subfigure[]{ \label{fig32b}
\includegraphics[ width=0.4\textwidth ] {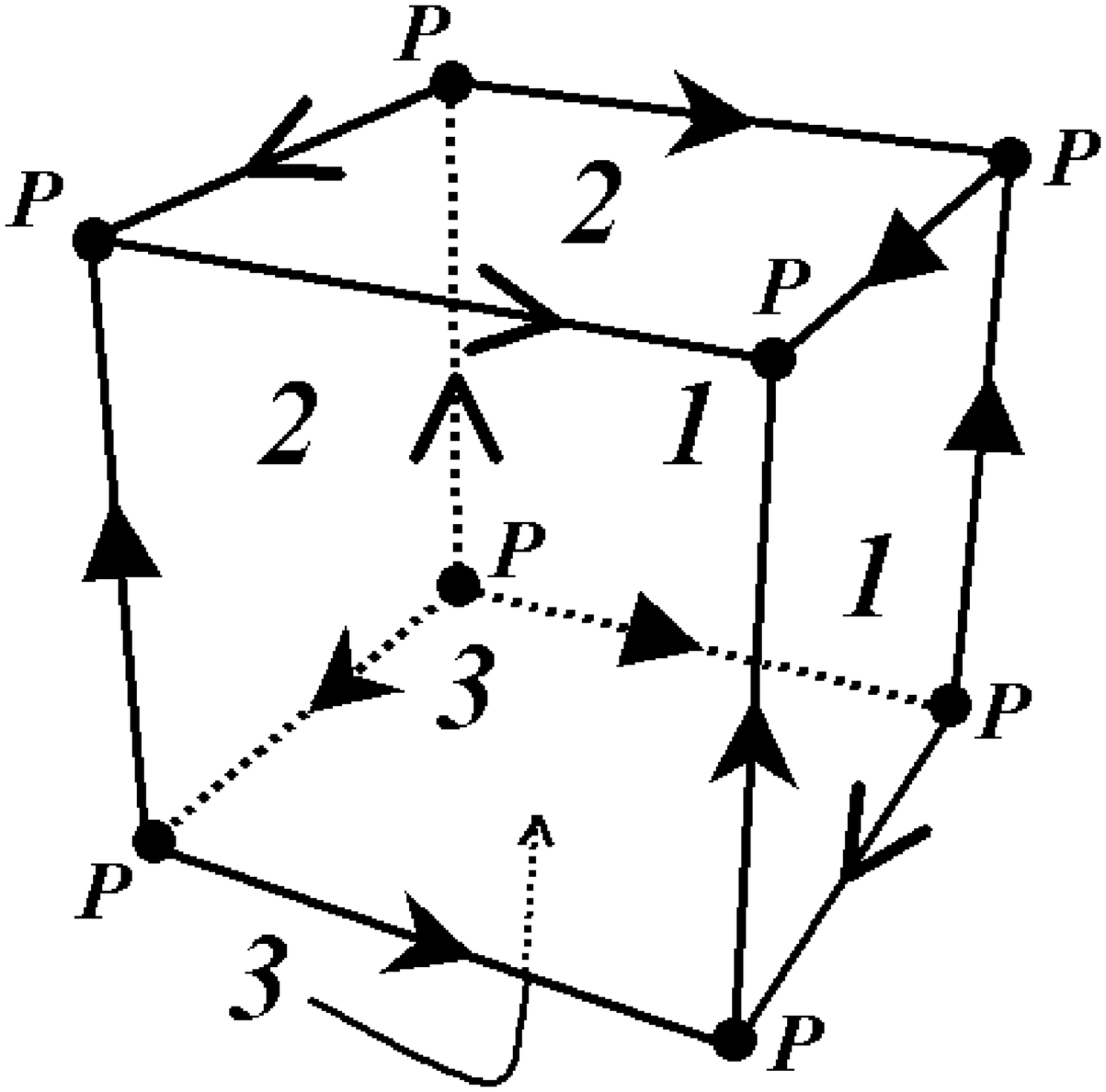}} \caption { }
\end{figure}

Now assume that we are given the realizable diagram $\mathcal{D}$
in the orientable surface $S$ and that we want to know if it is
the Johansson diagram of a filling Dehn sphere of some 3-manifold.
First of all, the diagram $\mathcal{D}$ must \textit{fill} the
surface $S$, (that is, $S-\mathcal{D}$ must be a disjoint union of
open 2-disks, and $\mathcal{D}-\left\{ \text{double points of
}\mathcal{D}\right\}  $ must be a disjoint union of open
intervals), and in particular the curves of $\mathcal{D}$ must
compose a connected graph on $S$. By Lemma \ref{LeMMA
G-Classes-singularity set} this implies that $\mathcal{D}$ has
exactly two opposite G-classes of neighbouring curves (because
$\mathcal{D}$ is realizable). As in the previous paragraph, the
construction of $\hat{M}(\mathcal{D})$ is uniquely determined by
the diagram and we have that $\mathcal{D}$ is the Johansson
diagram of a filling Dehn sphere of some 3-manifold $M\;$if and
only if $\partial\hat {M}(\mathcal{D})$ is a collection of
2-spheres. If this occurs, we say that $\mathcal{D}$ is a
\textit{filling diagram} and pasting a 3-ball to $\hat
{M}(\mathcal{D})$ along each boundary component of
$\hat{M}(\mathcal{D})$ we obtain the required closed 3-manifold
$M(\mathcal{D})$ that is also uniquely determined by the diagram
$\mathcal{D}$. The construction of $\hat {M}(\mathcal{D})$ can be
made in an algorithmic way. The diagrams of Figures \ref{fig29a},
\ref{fig29b} and \ref{fig32a} are all examples of (realizable)
filling diagrams. The diagram of Figure \ref{fig29a} appears in
the original paper \cite{Johansson1} and it is a Johansson
representation of the 3-sphere. Its corresponding filling Dehn
sphere is \textit{Johansson's sphere} (see Fig. 8 of
\cite{A.Shima2}). The diagram of Figure \ref{fig29b} represents
$S^{2}\times S^{1}$. The diagram of Figure \ref{fig32a} is a
diagram of a filling Dehn torus $\Sigma_{0}$ with only one triple
point in an Euclidean 3-manifold $M$. This euclidean manifold
coincides with the Seifert manifold
$M(S_{333})=(Oo0\mid-1;(3,1),(3,1),(3,1))$ (see \cite[p.
155]{MONTLIbro}), and it is the result of identifying the faces of
a solid cube in pairs as in Figure \ref{fig32b}. The filling Dehn
torus $\Sigma_{0}$ is the image in $M$ of the boundary of the cube
under this identification.

In Figure 33 is depicted how the Haken moves (except finger move $0$) for
immersions are reflected in the Johansson diagrams. These compose the
\textit{diagram moves}, and we denote them with the same names of the
corresponding moves of immersions. If we perform a diagram move in a filling
diagram, the move is \textit{filling-preserving} if the resulting diagram is
again a filling diagram.

\begin{figure}[b]
\centering \subfigure []{ \label{fig33a}
\includegraphics[ height=0.2\textheight
]{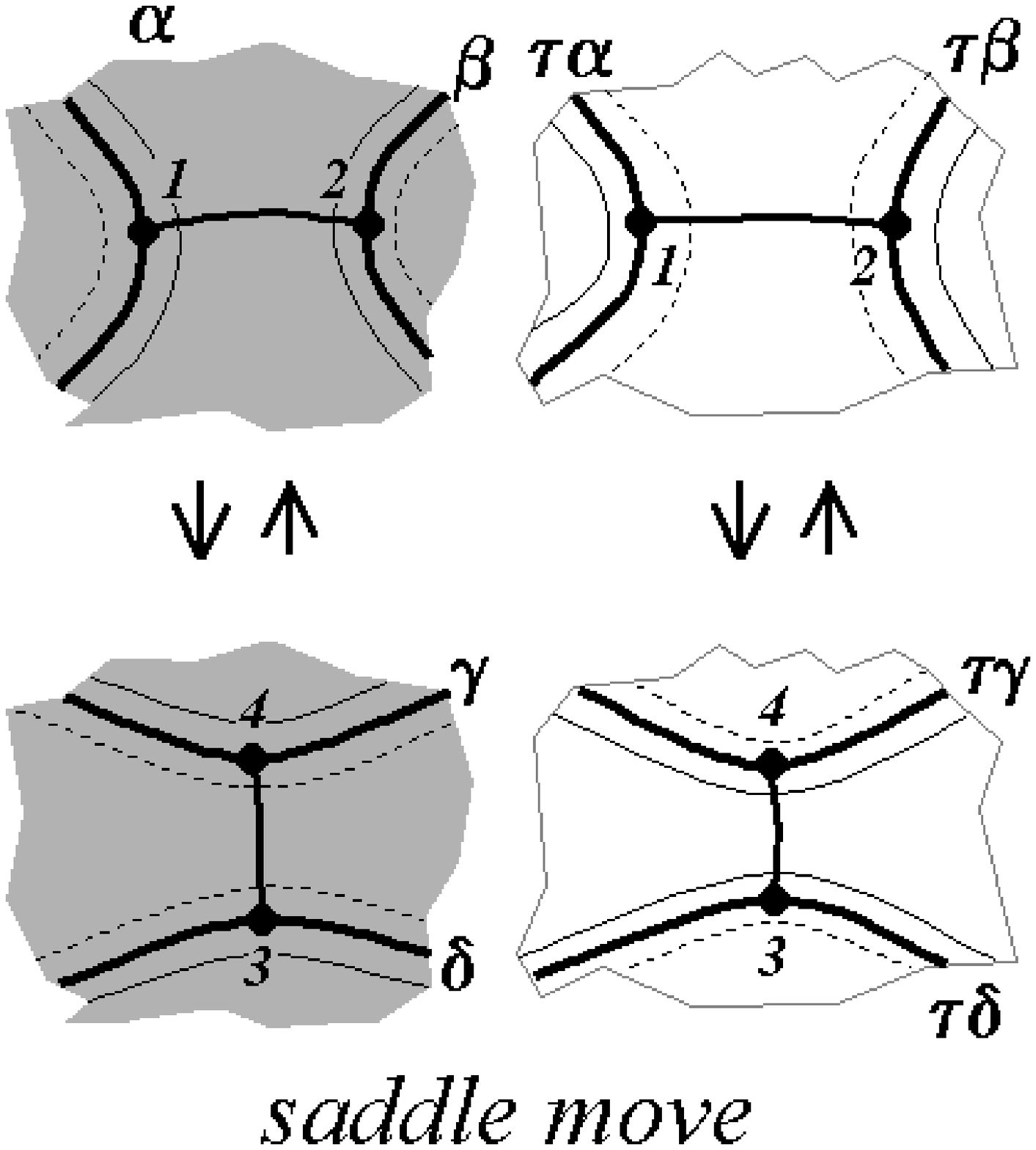}}\hspace{3cm} \subfigure[]{ \label{fig33b}
\includegraphics[ height=0.2\textheight
]{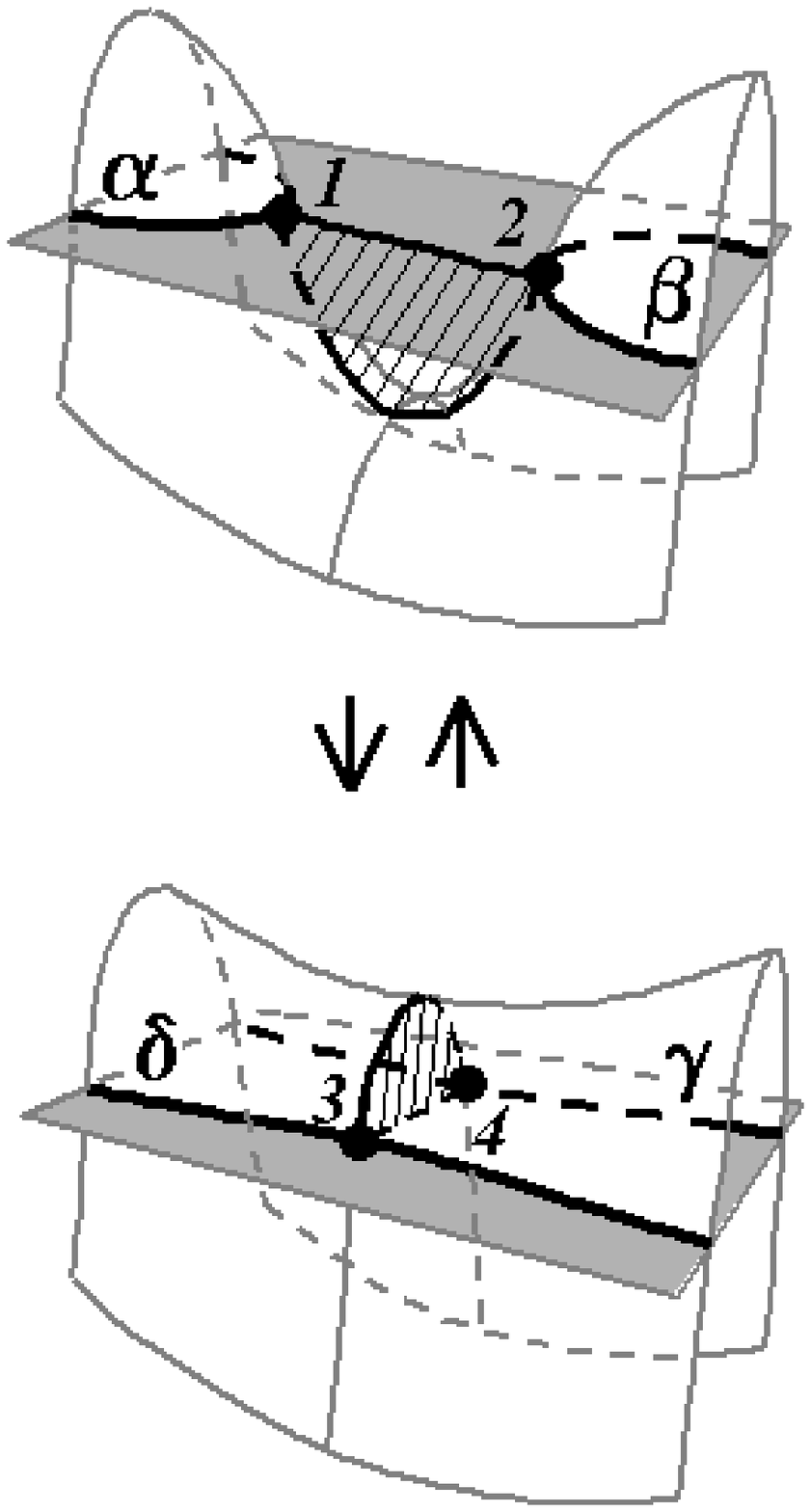}} \caption { }
\end{figure}
\addtocounter {figure}{-1}
\begin{figure}[h]
\centering \addtocounter {subfigure}{2} \subfigure []{
\label{fig33c}
\includegraphics[ height=0.2\textheight ]{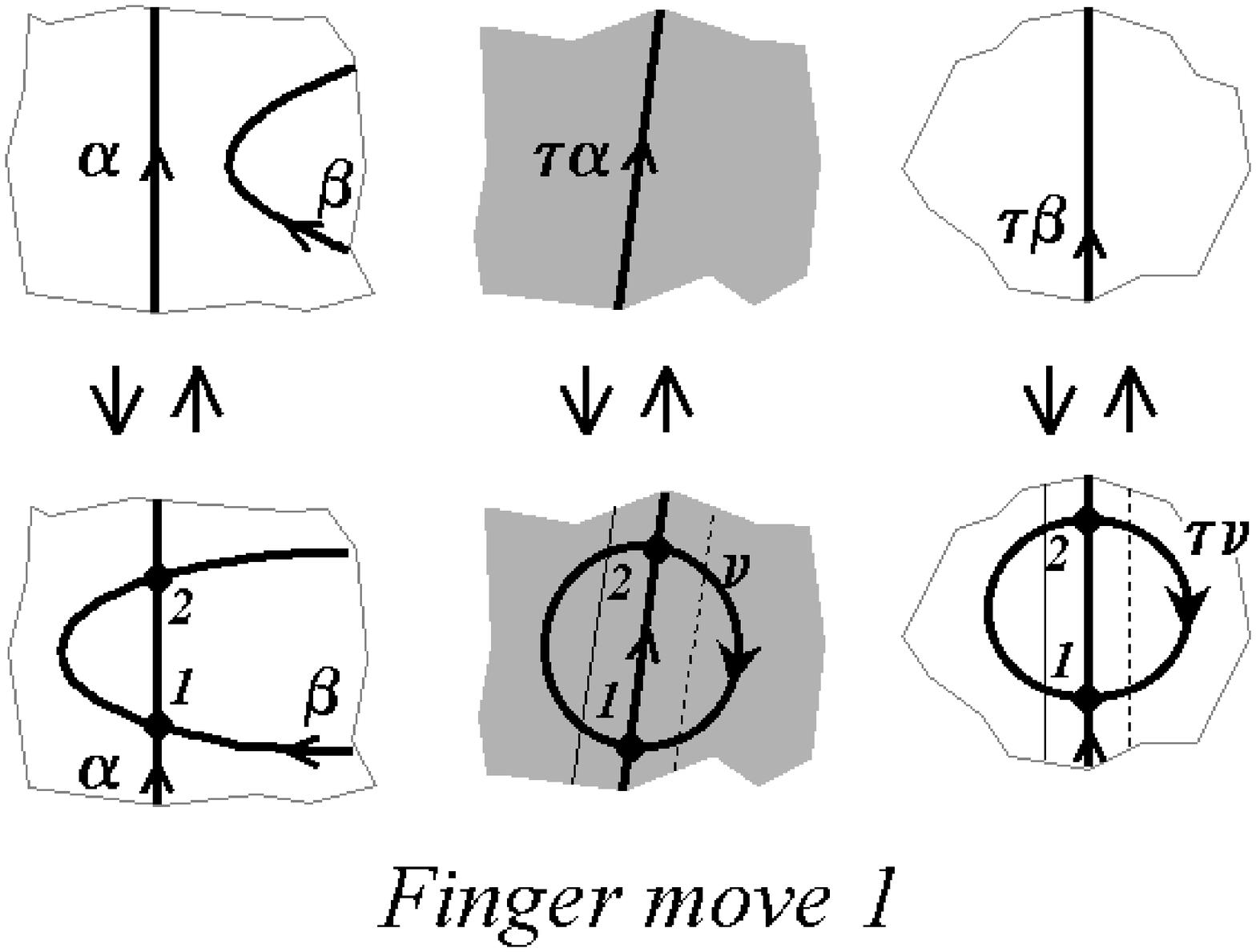}}\hspace {1.5cm} \subfigure[]{ \label{fig33d}
\includegraphics[ height=0.2\textheight ]{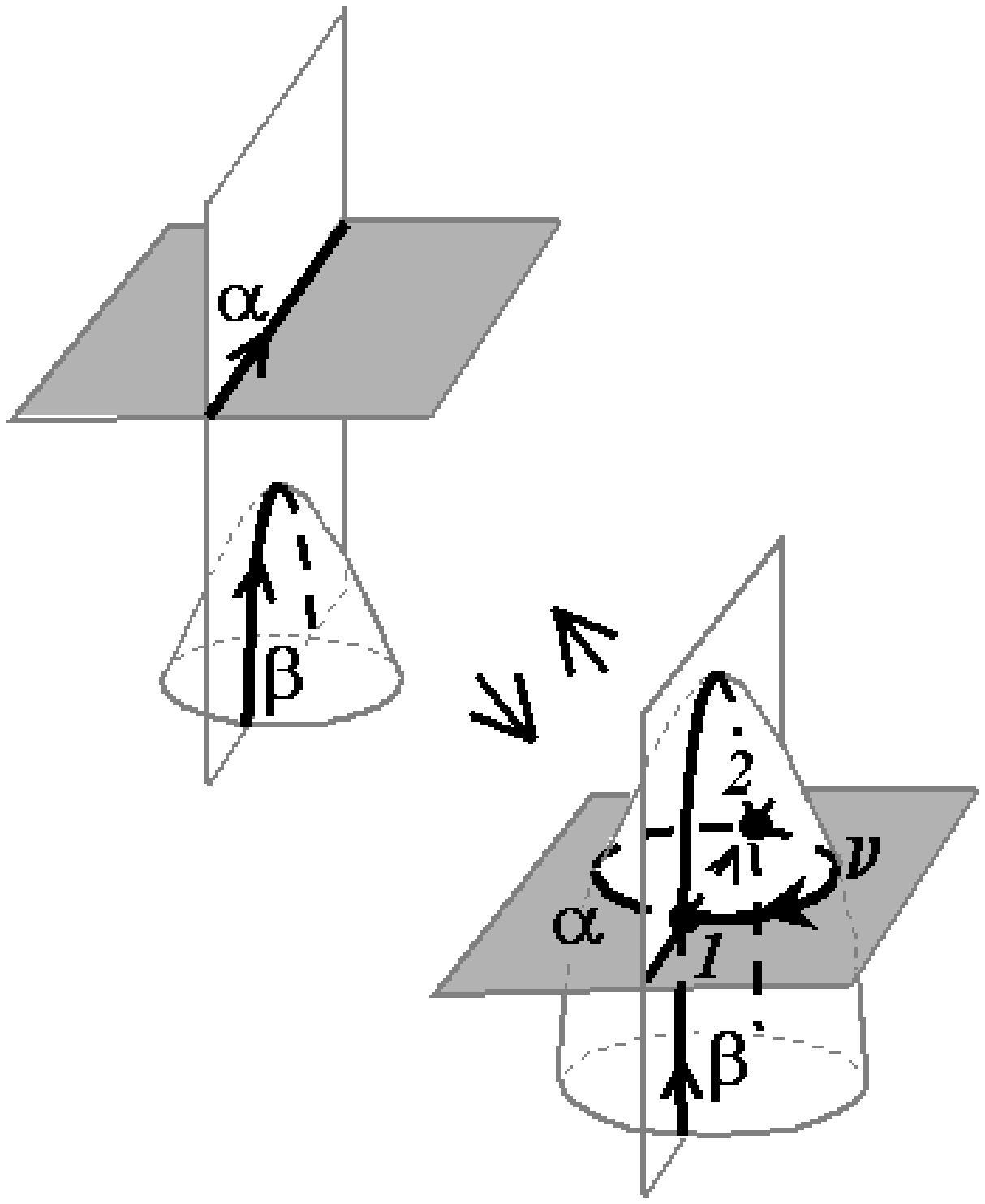}}\\
\subfigure []{ \label{fig33e}
\includegraphics[ height=0.2\textheight ]{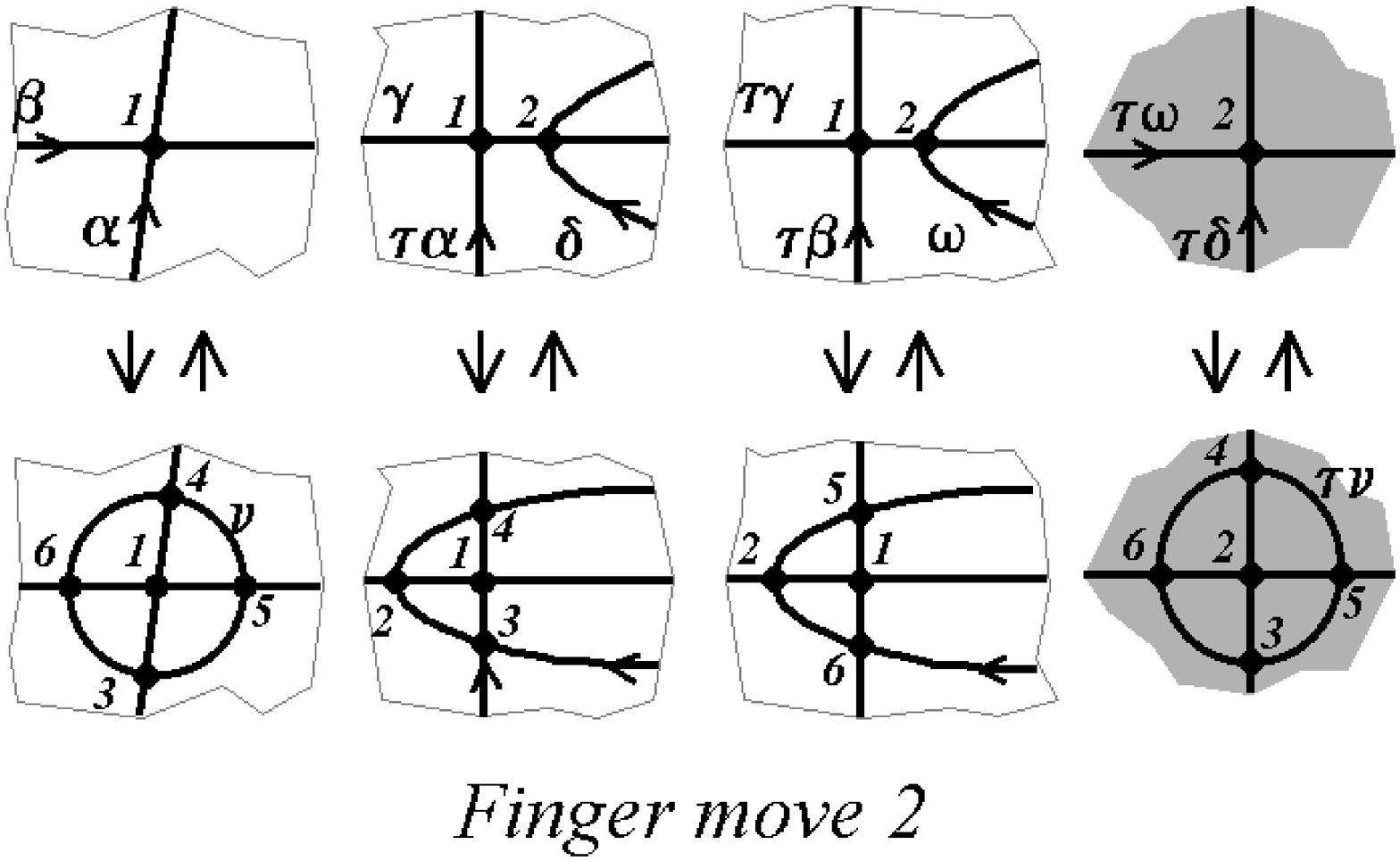}}\hfill \subfigure[]{ \label{fig33f}
\includegraphics[ height=0.2\textheight ]{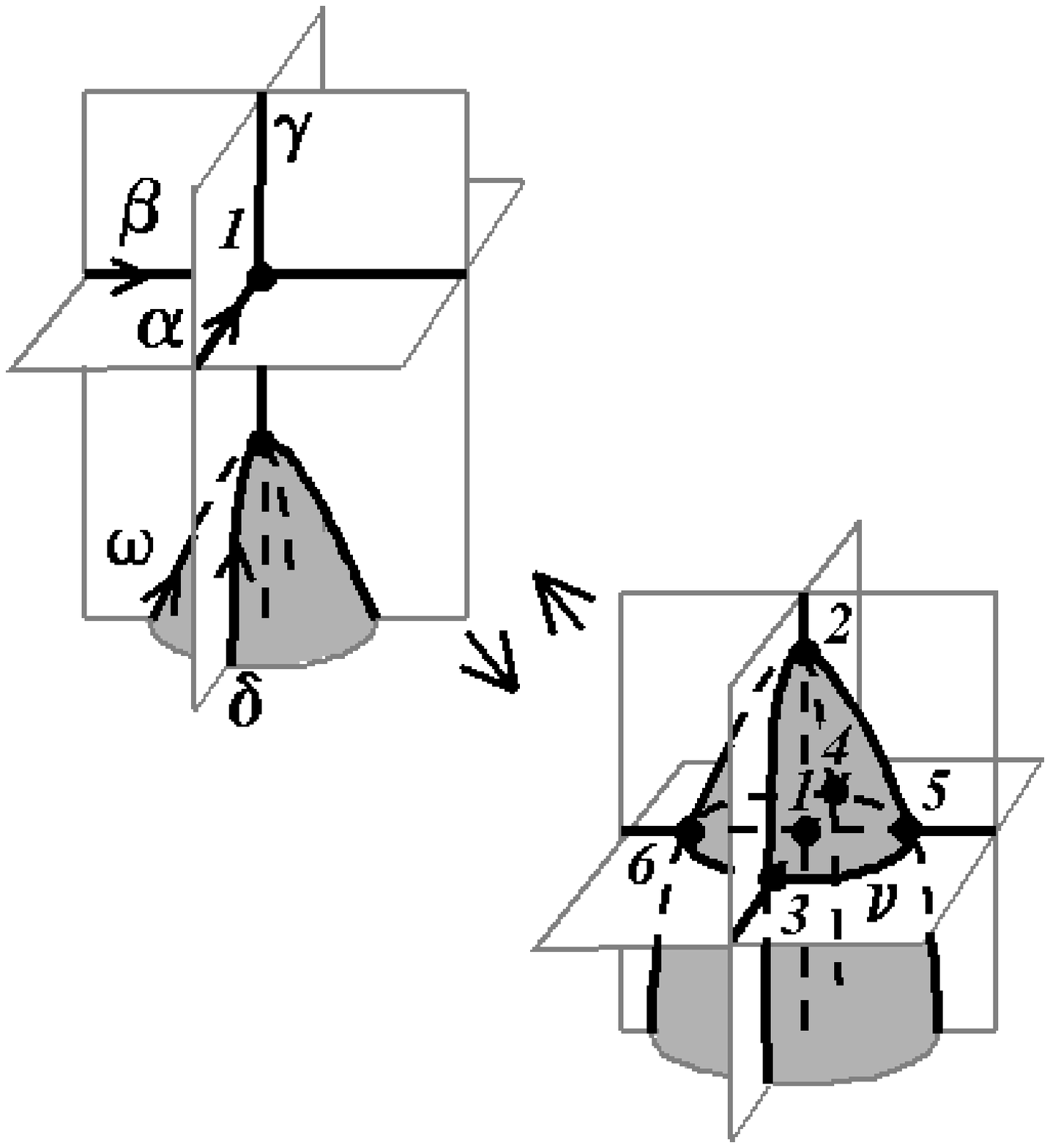}} \caption { }
\end{figure}

Let $\mathcal{D}$ be a filling diagram on the surface $S$, and let $i:\hat
{M}(\mathcal{D})\hookrightarrow M(\mathcal{D})$ the inclusion map. Let denote
simply by $f$ the immersion $i\circ\pi\circ j:S\rightarrow M(\mathcal{D})$
that realizes $\mathcal{D}$, and put $\Sigma=f(S)$.

If we perform a filling-preserving diagram move in $\mathcal{D}$, this move
will come from a filling-preserving move of $f$ and thus the new diagram
$\mathcal{D}^{\prime}$ we obtain verifies that $M(\mathcal{D}^{\prime
})=M(\mathcal{D})$. As it happened with Haken moves, a finger move $\pm$1 or
$\pm$2 in the filling diagram $\mathcal{D}$, will always be
filling-preserving, and a saddle move on $\mathcal{D}$ may be or not a
filling-preserving one.

The neighbouring curves of the diagram help us to perform the saddle move and
the finger move 1.

A saddle move can be performed in the diagram $\mathcal{D}$
everytime we have two arcs connecting related points of the
diagram as in Figure \ref{fig33a}. The two neighbouring curves of
$\mathcal{D}$ that intersect any of the two arcs must belong to
the same G-class. If we have such a pair of arcs, because
$\mathcal{D}$ is a filling diagram (the image by $f$ of) these
arcs must bound a 2-gon $w$ in $M(\mathcal{D})-\Sigma$ as in
Figure \ref{fig33b}, and we can perform a saddle move on $f$ by
pushing along $w$ any of the two sheets of $\Sigma$ that bound
$w$. This saddle move on $f$ is reflected in the saddle move of
the diagram $\mathcal{D}$ as depicted in Figure \ref{fig33a}.

If we perform a finger move +1 on the diagram $\mathcal{D}$ there
appear a new pair $\nu,\tau\nu$ of sister curves of the diagram.
The diagram $\mathcal{D}$ tells us how we must identify the new
double points (labeled $1$ and $2$ in Figure \ref{fig33c}) that
appear in the new diagram $\mathcal{D}^{\prime} $, but there is
some ambiguity (that do not occur for finger moves 2) because
there are two ways for identifying $\nu$ with $\tau\nu$ (for a
given orientation of $\nu$ there are two possible orientations of
$\tau\nu$). This ambiguity disappears when we draw the
neighbouring curves of the diagram $\mathcal{D}^{\prime}$, using
that related neighbouring points of the diagram must lie on
neighbouring curves of the same G-class.

We say that a filling diagram $\mathcal{D}$ on a surface $S$ is
\textit{nulhomotopic} if the immersion $f=i\circ\pi\circ j$ as
above is nulhomotopic. The diagram of Figure \ref{fig29a} is
nulhomotopic (every diagram representing $S^{3}$ must be
nulhomotopic), while the diagram of Figure \ref{fig29b} is not
nulhomotopic.

The following result is a Corollary of Theorem \ref{MAINtheorem}.

\begin{corollary}
\label{CORtoMAINtheorem}Two nulhomotopic filling diagrams on $S^{2}$
represents the same 3-manifold if and only if they are related by a finite
sequence of filling preserving moves.
\end{corollary}

\section{Duplication.\label{SECTION Duplication}}

It is possible to obtain algorithmicaly a nulhomotopic Johansson
representation of $M$ from \textit{any}, nulhomotopic or not, Johansson
representation $\mathcal{D}$ of $M$. We will call this process
\textit{duplication} of diagrams and it can be made using Johansson's
construction of $\hat{M}(\mathcal{D})$ as follows.

Let $f:S^{2}\rightarrow M$ be a filling immersion, and put
$\Sigma=f(S)$. Take a thickening $F:S^{2}\times\left[  -1,1\right]
\rightarrow M$ of $f$ as in the previous section, such that near a
double curve of $\Sigma$ the image of $F$ intersects itself as in
Figure \ref{fig31}, and consider the upper sheet
$\Sigma^{+}=F(S^{2}\times\left\{  1\right\}  )$ and the lower
sheet $\Sigma^{-}=F(S^{2}\times\left\{  -1\right\}  )$, which are
two filling Dehn spheres of $M$ parallel to $\Sigma$ at both sides
of $\Sigma$. The Johansson diagram $\mathcal{D}$ of $f$ has two
G-classes of neighbouring curves. We can take the neighbouring
curves of $\mathcal{D}$ such that their images by $f$ compose the
intersection of $\Sigma$ with $\Sigma^{-}\cup\Sigma^{+}$. We call
the \textit{upper} (\textit{lower}) G-class to the G-class of
neighbouring curves of $\mathcal{D}$ whose image by $f$ is
contained in the upper (lower) sheet $\Sigma^{+}$ ($\Sigma^{-}$).

\begin{figure}[htb]
\centering \subfigure []{ \label{fig34a}
\includegraphics[ height=0.2\textheight ] {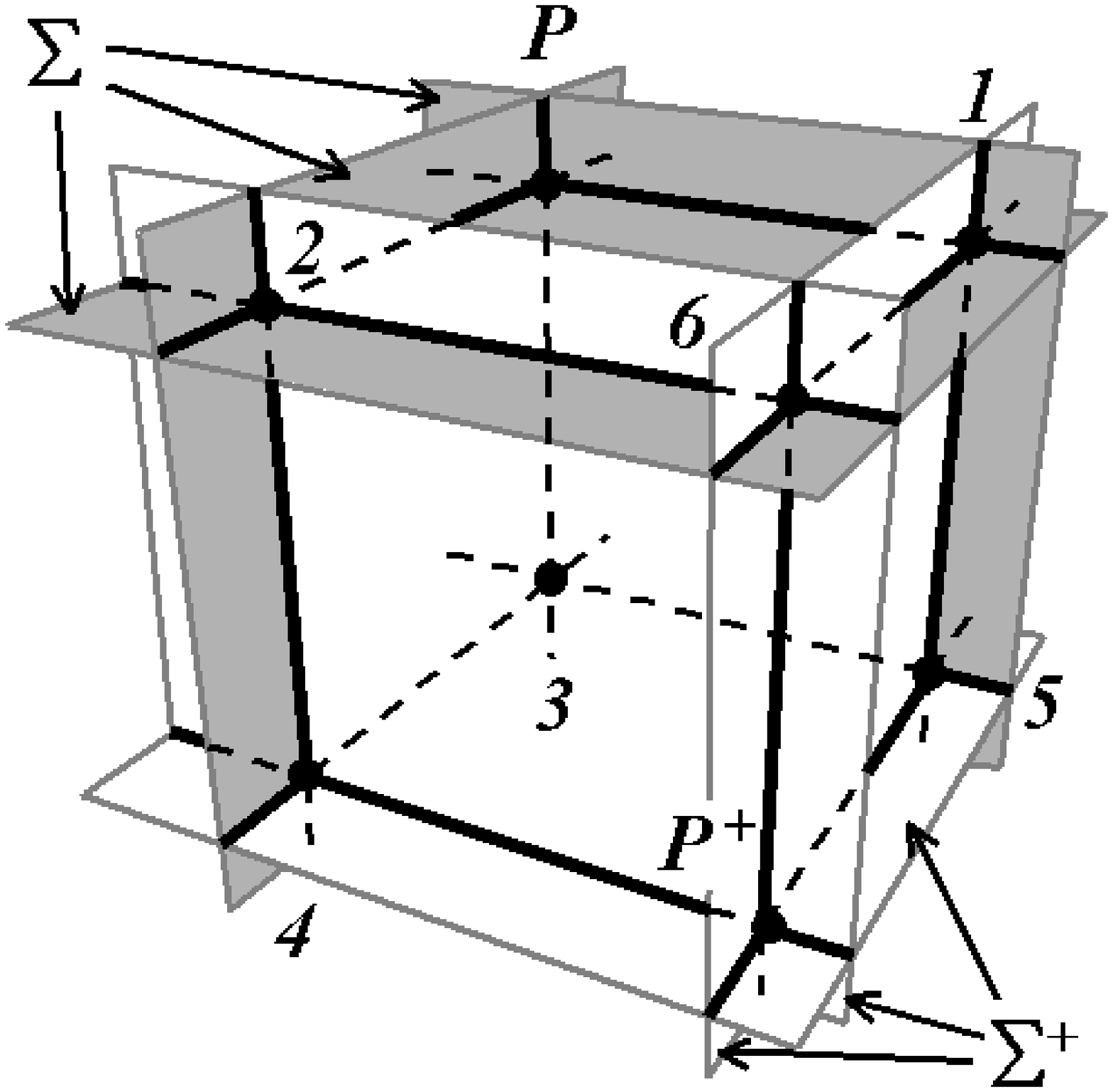}}\hspace{0.4cm} \subfigure {
\includegraphics[ height=0.2\textheight ] {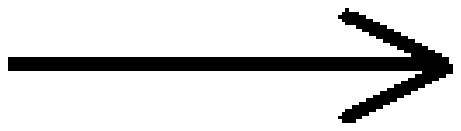} }\hspace{0.4cm} \addtocounter {subfigure}{-1}
\subfigure[]{ \label{fig34b}
\includegraphics[ height=0.2\textheight ] {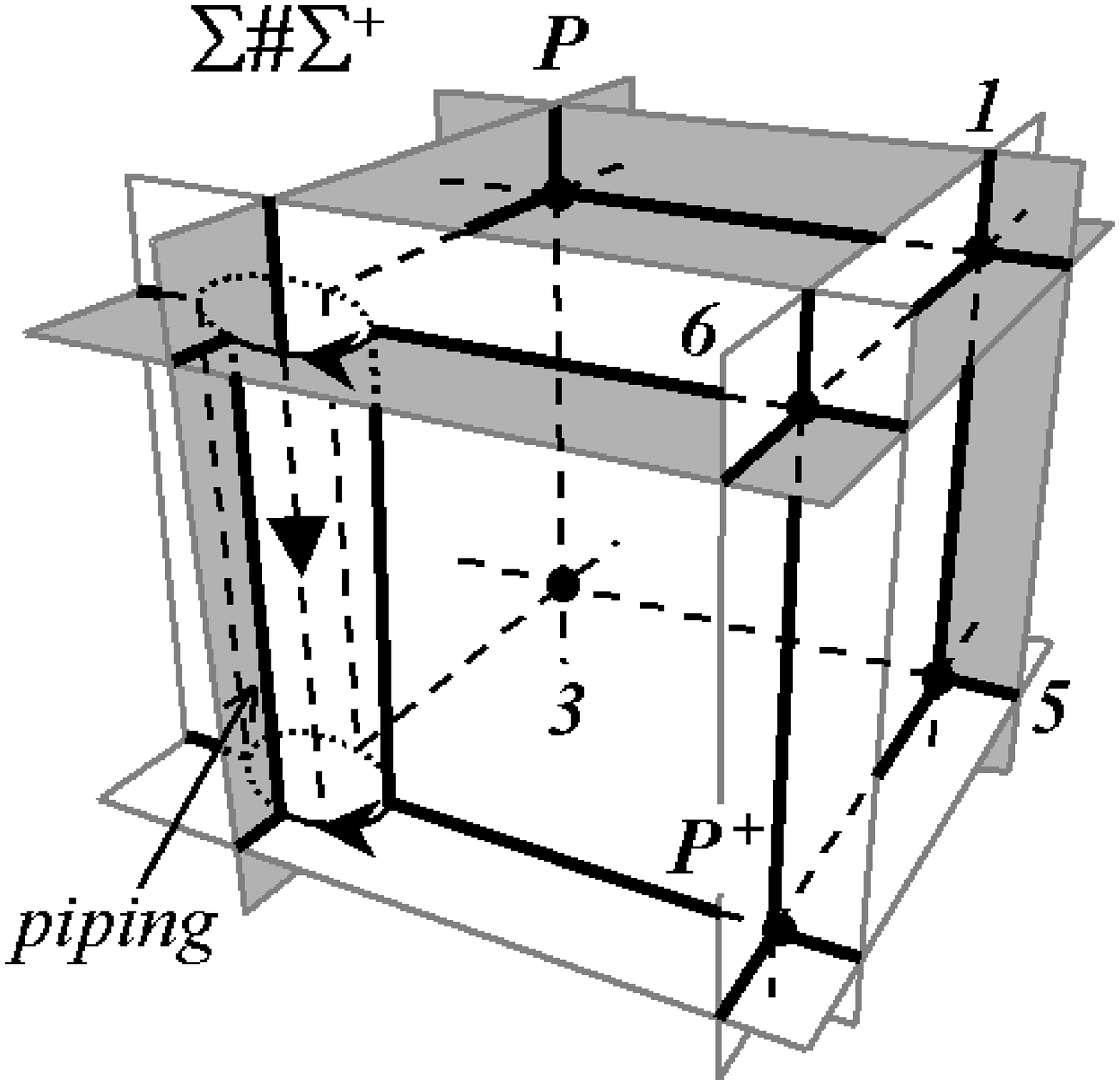}} \caption { }
\end{figure}

The Dehn spheres $\Sigma,\Sigma^{+}$ form a filling pair of
spheres in $M$. Near a triple point $P$ of $\Sigma$ there will be
eight triple points of the union $\Sigma\cup\Sigma^{+}$ as in
Figure \ref{fig34a}. In some situations the Dehn sphere
$\Sigma\#\Sigma^{+}$ that we obtain piping $\Sigma$ with
$\Sigma^{+}$ near $P$ as in Figure \ref{fig34b} is a filling Dehn
sphere of $M$. Assume that this is the case. This filling Dehn
sphere $\Sigma\#\Sigma^{+}$ is the image by $F$ of the boundary of
$S^{2}\times\left[  0,1\right]  $ with a small cylinder connecting
$S^{2}\times\left\{  0\right\}  $ with $S^{2}\times\left\{
1\right\}  $ removed, and thus it is nulhomotopic.

\begin{figure}[htb]
\centering \subfigure []{ \label{fig35a}
\includegraphics[ height=0.15\textheight ] {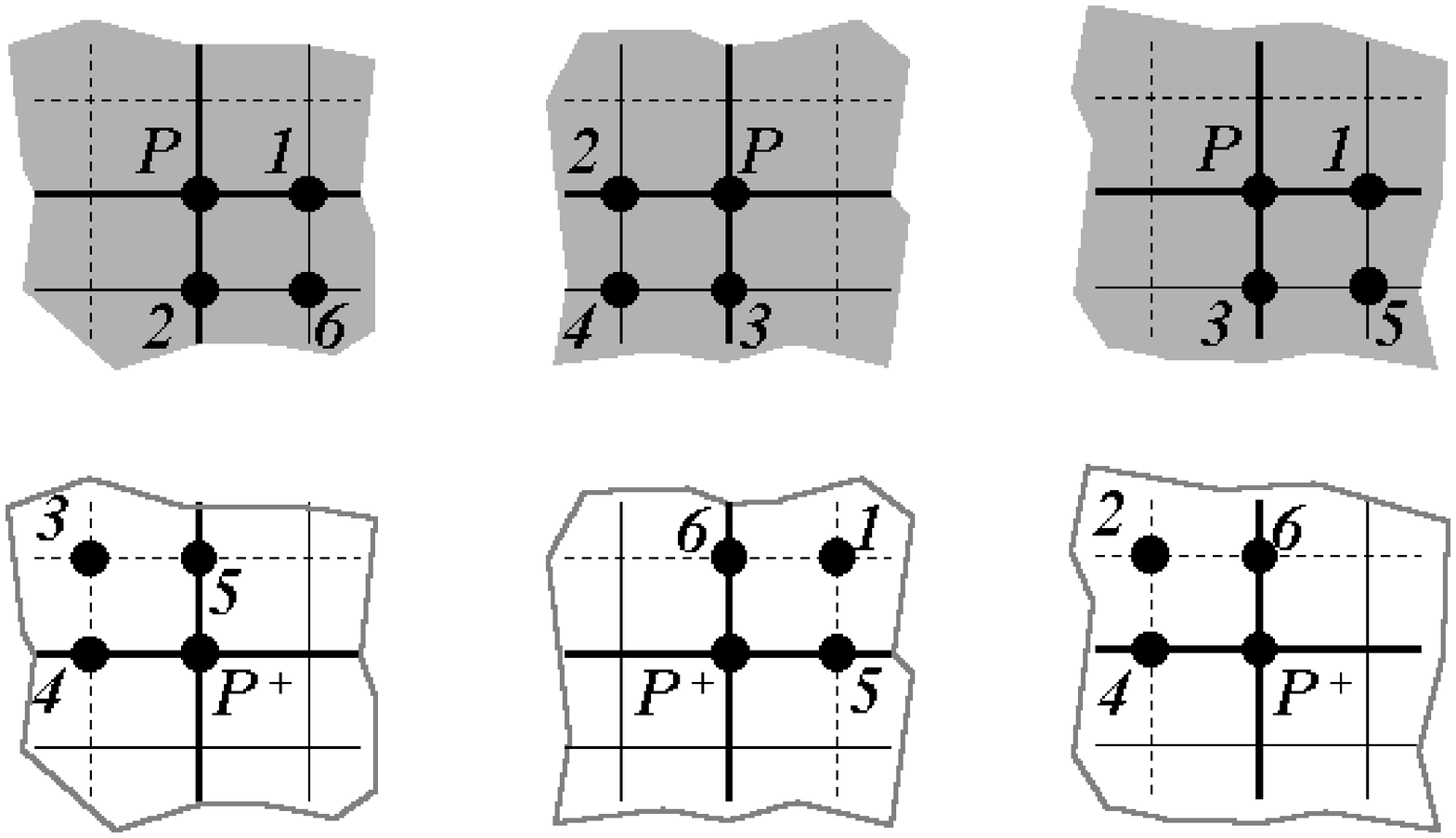}}\hfill \subfigure {
\includegraphics[ height=0.15\textheight ] {fig35ab} }\hfill
\addtocounter {subfigure}{-1} \subfigure[]{ \label{fig35b}
\includegraphics[ height=0.15\textheight ] {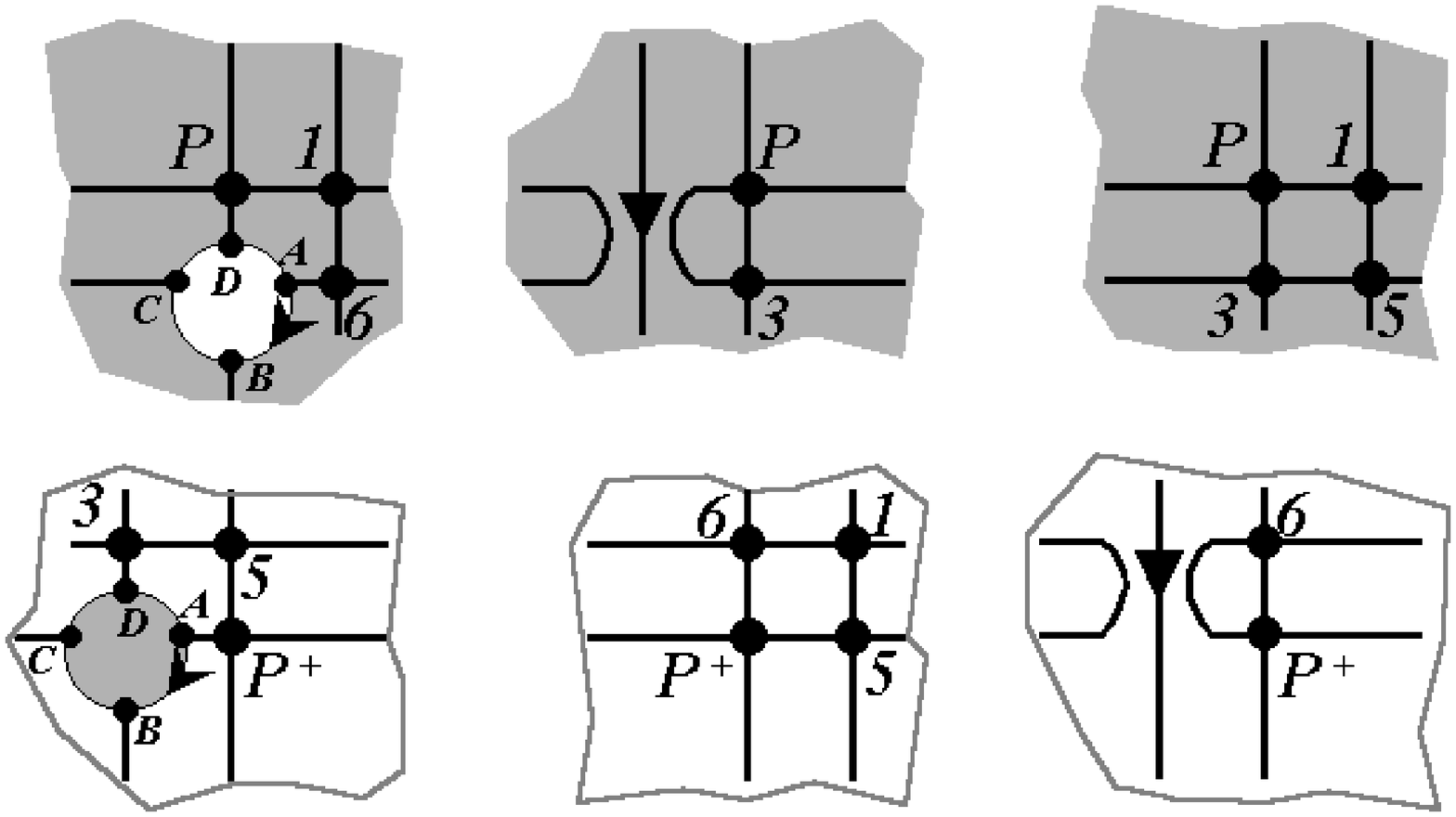}} \caption { }
\end{figure}

The Johansson diagram of $\Sigma\#\Sigma^{+}$ can be obtained
algorithmically from the Johansson diagram $\mathcal{D}$ of
$\Sigma$. Let $f^{+} :S^{2}\rightarrow M$ be the parametrization
of $\Sigma^{+}$ given by $f^{+}(X)=F(X,1)$.

First, remember that the upper G-class of neighbouring curves of
$\mathcal{D} $ composes the inverse image by $f$ of
$\Sigma\cap\Sigma^{+}$. A second observation is that the Johansson
diagram of $f^{+}$ is a copy $\mathcal{D} ^{+}$ of $\mathcal{D}$.
The third observation is the following

\begin{remark}
\label{REMARKparallelDehnSpheres}The relative position of $\Sigma$ with
respect to $\Sigma^{+}$ is exactly the same as the relative position of
$\Sigma^{-}$ with respect to $\Sigma$.
\end{remark}

This Remark implies that if we call the G-classes of
$\mathcal{D}^{+}$ with the same name as their respective copies in
$\mathcal{D}$, \textit{the lower G-class of
}$\mathcal{D}^{+}$\textit{\ composes the inverse image by }
$f$\textit{\ of }$\Sigma\cap\Sigma^{+}$.

The Johansson diagram $\mathcal{D}\#\mathcal{D}^{+}$ of
$\Sigma\#\Sigma^{+} $ is what we call a \textit{duplicate} of
$\mathcal{D}$, and it can be obtained from $\mathcal{D}$ and
$\mathcal{D}^{+}$ as indicated in Figure 35. In Figure
\ref{fig35a} we depict how the cube of Figure \ref{fig34a} is seen
in the Johansson diagrams of $\Sigma$ (grey) and $\Sigma^{+}$
(white). In both diagrams we draw the lower G-class and the upper
G-class similarly. As an example we apply this construction to the
diagram of Figure \ref{fig29b} in Figure 36.

It is not always possible to obtain a triple point of $\Sigma$
such that piping $\Sigma$ with $\Sigma^{+}$ near $P$ as in Figure
\ref{fig34b} the resulting Dehn sphere $\Sigma\#\Sigma^{+}$ fills
$M$, but \textit{we can always deform }$\Sigma$\textit{\ by
filling-preserving moves} to obtain another filling Dehn sphere
$\Sigma^{\prime}$ with a triple point where the
(filling-preserving) duplication is possible. This deformation can
be made very easily: for example, if $P$ is one of the two triple
points that appear after a finger move +1, then duplication is
possible at $P$.

\begin{figure}[b]
\centering \subfigure []{ \label{fig36a}
\includegraphics[ width=0.7\textwidth ] {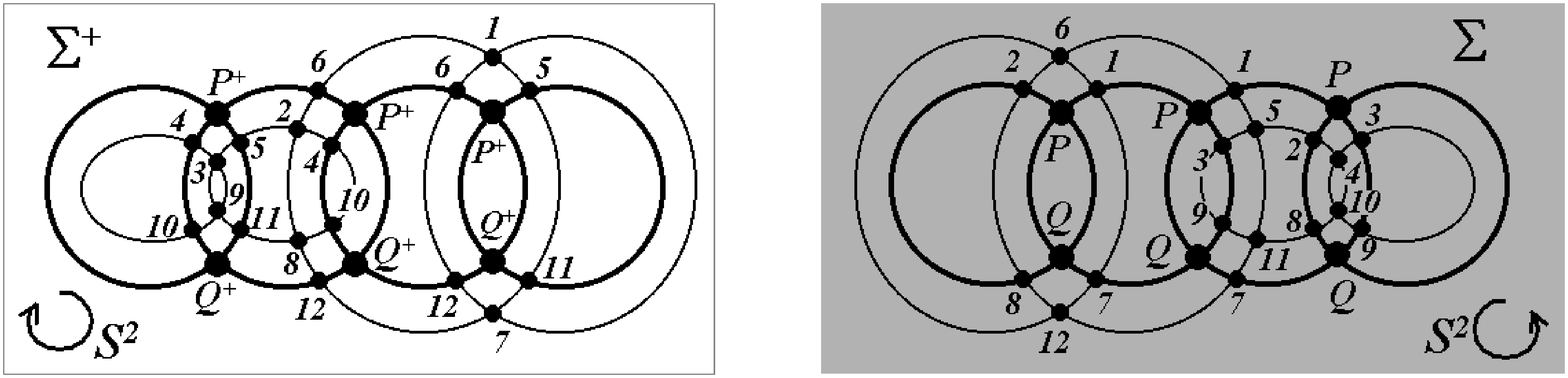}} \subfigure []{
\label{fig36b}
\includegraphics[ width=0.7\textwidth ] {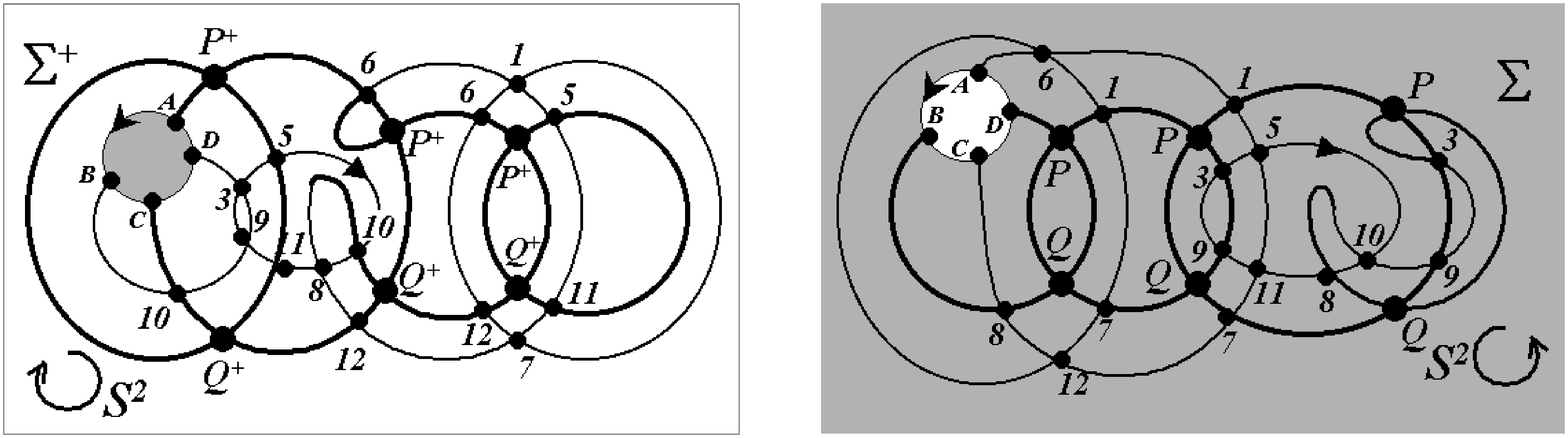}} \subfigure []{
\label{fig36c}
\includegraphics[ width=0.55\textwidth ]{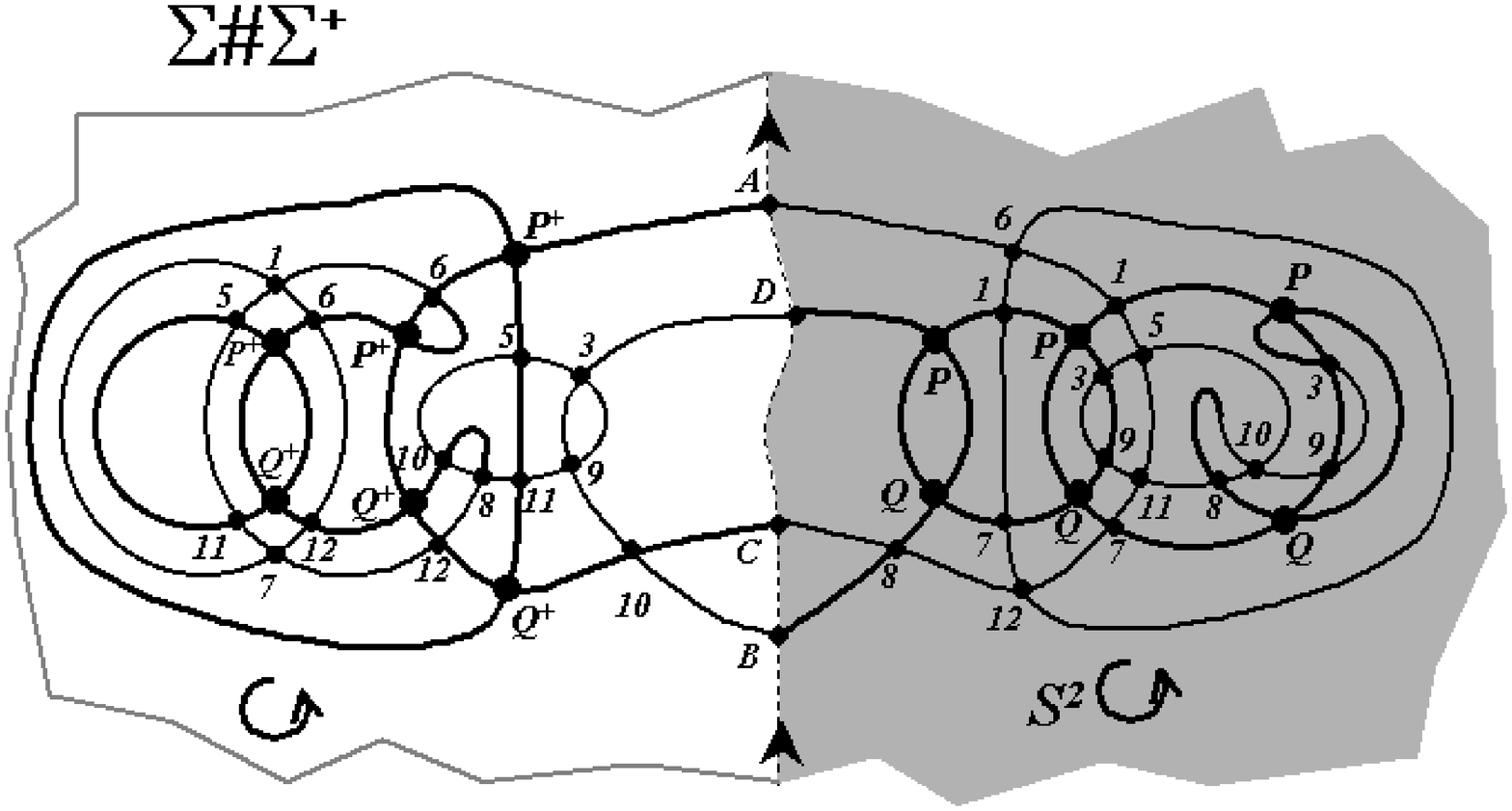}} \caption { }
\end{figure}

Duplication allows us to give a more general version of Corollary
\ref{CORtoMAINtheorem}:

\begin{theorem}
\label{THEOREMDuplicates}Two filling diagrams on $S^{2}$ represents the same
3-manifold if and only if their duplicates are related by a finite sequence of
filling preserving moves.
\end{theorem}

\section{Miscelany.\label{SECTION Miscelany}}

\subsection{Invariants of 3-manifolds.\label{SUBSECTION Invariants}}

The first immediate application of Corollary \ref{CORtoMAINtheorem} is the
search of invariants of 3-manifolds. If we could assign to each nulhomotopic
diagram on $S^{2}$ an object which remains invariant under filling-preserving
diagram moves, then this object defines a 3-manifold invariant. If $\varphi$
denotes such an invariant, for computing $\varphi$ for a given manifold $M$ we
would need a nulhomotopic Johansson representation of $M$. If we have an
arbitrary Johansson representation $\mathcal{D}$ of $M$ and we don't know if
it is a nulhomotopic one, duplicating $\mathcal{D}$ we will be able to compute
$\varphi$ from $\mathcal{D}$. Nevertheless, duplication produce very
complicated diagrams (the number of triple points of a duplication of
$\mathcal{D}$ is eight times the number of triple points of $\mathcal{D}$
minus $2$), and for this reason it should be interesting to know how to decide
if a given filling diagram in $S^{2}$ is nulhomotopic or not. This is an open
problem. In \cite{Montesinos} it is indicated an algorithm to obtain a
nulhomotopic Johansson representation of $M$ from any Heegaard diagram of $M$.
A simpler algorithm is studied in detail in \cite{Anewproof}.

\subsection{The diagram group.\label{SUBSECTION The diagram Group}}

Let $\mathcal{D}$ be a realizable diagram on $S^{2}$, and let
$f:S^{2} \rightarrow M$ be a transverse immersion parametrizing
$\mathcal{D}$. There is an easy way for obtaining a presentation
of the fundamental group of $\Sigma:=f(S^{2})$ in terms of the
diagram $\mathcal{D}$. If $\mathcal{D} =\left\{
\alpha_{1},...,\alpha_{n}\right\}  $, then we define the
\textit{diagram group}
\[
\pi\left(  \mathcal{D}\right)  =\left|
\alpha_{1},...,\alpha_{n}:\alpha
_{1}\cdot\tau\alpha_{1}=...=\alpha_{n}\cdot\tau\alpha_{n}=r_{1}=...=r_{k}
=1\right|  \text{ ,}
\]
where the relators $r_{1},...,r_{k}$ are given by the triplets
points of $\mathcal{D}$. If $P_{1},...,P_{k}$ are the triple
points of $\Sigma$, and $P_{j}$ is the triple point $P$ of Figure
\ref{fig28a}, which is reflected in the triplet of $\mathcal{D}$
of Figure \ref{fig28b}, then the associated relator is
\[
r_{j}=\alpha\beta\gamma\text{ .}
\]

It can be proved that $\pi_{1}\left(  \Sigma\right)  $ is isomorphic to the
diagram group \cite{Tesis} and thus, if $\Sigma$ fills $M$ the diagram group
gives a presentation of the fundamental group of $M$. This presentation is due
to W. Haken (see Problem 3.98 of \cite{KirbyProblems}), and we have not seen
it in printed form.

With this presentation, it can be proved the following (unpublished) theorem
of W. Haken (see also \cite{Haken1}). If $\mathcal{D}$ is a connected
realizable diagram in $S^{2}$ with only two double curves $\alpha,\tau\alpha$,
then both are simple or both are non-simple.

\begin{theorem}
\label{THeoremDtwocurvesSimplyConnectedorZ3}If $\alpha$ and $\tau\alpha$ are
simple, then $\pi\left(  \mathcal{D}\right)  \simeq\mathbb{Z}_{3}$, and if
both are non-simple, then $\pi\left(  \mathcal{D}\right)  \simeq1$.
\end{theorem}

\subsection{Checking fillingness.\label{SUBSECTION Checking fillingness}}

As we have said, to test if a realizable diagram $\mathcal{D}$ is a filling
diagram is to check if $\partial\hat{M}(\mathcal{D})$ is a collection of
2-spheres. Though the complete construction of the manifold with boundary
$\hat{M}(\mathcal{D})$ from the diagram $\mathcal{D}$ can be made in an
algorithmic way using Johansson's construction, it is interesting to have
faster methods for checking fillingness. The following result will give us a
method for saving time in this process. It can be proved easily using Euler's
characteristic techniques.

\begin{lemma}
\label{LEMMADfills if P+2 boundary com}A realizable diagram $\mathcal{D}$ on
the genus $g$ surface $S$ is a filling diagram if and only if it fills $S$ and
$\partial\hat{M}(\mathcal{D})$ has $p+\chi(S)$ connected components, where $p$
is the number of triplets (of pairwise related double points) of $\mathcal{D}$.
\end{lemma}

The diagram group can help us also for checking fillingness. In Lemma 4.9 of
\cite{Hempel}, it is proved that if a 3-manifold with boundary $\hat{M}$ has a
boundary component which is not a 2-sphere, then $\hat{M}$ has a double cover. Thus,

\begin{lemma}
\label{LEMMANosubgroupOfindexTWOfills}If $\mathcal{D}$ is a
realizable connected diagram in $S^{2}$ and the diagram group
$\pi\left(  \mathcal{D} \right)  $ has no subgroup of index $2$,
then $\mathcal{D}$ is a filling diagram.
\end{lemma}

This Lemma, together with Theorem \ref{THeoremDtwocurvesSimplyConnectedorZ3}
gives the following.

\begin{corollary}
\label{CORDIagramsWithtwoCurvesFill}If $\mathcal{D}$ is a realizable connected
diagram in $S^{2}$ with only two curves, then it is a filling diagram.\bigskip
\end{corollary}

\subsection{\label{SUBSECTION Eversion}Filling eversion.}

Theorem applies not only to filling Dehn spheres but to their
parametrizations. Let $f:S^{2}\rightarrow M$ be a filling
immersion, and consider also the immersion $g:S^{2}\rightarrow M$
given by $g=f\circ a$, where $a$ denotes now the antipodal map of
$S^{2}$. In this situation, the fact that $f(S^{2})$ is filling
homotopic to $g(S^{2})$ is trivial because they are the same Dehn
sphere of $M$. Theorem \ref{MAINtheorem} asserts that $f$ is
filling homotopic to $g$, which is now a non-trivial fact. This is
a filling version of the problem of the eversion of the 2-sphere
(see \cite{Morin-Petit}). In \cite{Max-Banchoff} it is proved that
every eversion of the 2-sphere in $S^{2}$ has at least one
quadruple point (compare also \cite{NowikTahl}). Using this, it
can be seen that any parametrization $f$ of the Johansson's sphere
and its antipodal parametrization $g$ cannot be taken into one
another by using only finger moves 1 and filling-preserving saddle
moves. This means that finger move 2 cannot be dispensed off from
the statement of Theorem \ref{MAINtheorem}.

\subsection{The non-orientable case.\label{SUBSECTION Non-orientable}}

In our discussion about diagrams in section \ref{SECTION Diagrams} we have
assumed that surfaces and 3-manifolds are orientable. If $f:S\rightarrow M$ is
a transverse immersion and $S$ or $M$ (or both) are non-orientable, then the
inverse image by $f$ of a double curve $\bar{\alpha}$ of $f$ may be a unique
closed curve $\alpha$ in $S$, such that $\alpha$ is a 2-fold covering of
$\bar{\alpha}$ (see \cite{Johansson1}). For this reason, in this general case
we need a more general definition of abstract diagram. The one given in
section \ref{SECTION Diagrams} can be adapted to this general case by allowing
a \textit{non-free} involution $\tau$ and requiring that the curves $\alpha$
of the diagram with $\tau\alpha=\alpha$ commute with the antipodal map of
$S^{1}$. Another fact is that in the general case, the immersion $f$ might be
\textit{1-sided}, and thus the proof of Theorem \ref{THM Johansson G-Clases}
breaks down.

Johansson proves in \cite{Johansson1} a Theorem characterizing
diagrams in $S^{2}$ which are realizable in 3-manifolds
(orientable or not). It is an interesting problem to generalize
this last theorem of Johansson to cover diagrams in any surface
and immersions in general 3-manifolds (orientable or not). This is
certainly not very difficult but it is an intermediate step to
generalize our theory of filling immersions to cover general Dehn
surfaces (orientable or not) immersed in general 3-manifolds
(orientable or not). This is an open program.

\subsection{A question of R. Fenn.\label{SUBSECTION Fenn's question}}

The following question was asked to us by R. Fenn:

\textit{Do filling Dehn surfaces in }$M$\textit{\ lift to embeddings in
}$M\times\left[  0,1\right]  $\textit{?}

We do not know the complete answer to this question. In
\cite{Giller} (see also \cite{Carter-Saito}) it is given an
algorithm for deciding if a Dehn $\Sigma$ surface in
$\mathbb{R}^{3}$ lift to an embedding in $\mathbb{R}^{4}$ in terms
of the Johansson diagram of $\Sigma$. In the same paper it is
given an example of a Dehn sphere $\Sigma_{1}$ in $\mathbb{R}^{3}$
that does not lift to an embedding in $\mathbb{R}^{4}$. The
Johansson diagram of $\Sigma _{1}$ has only two non-simple curves,
and by Corollary \ref{CORDIagramsWithtwoCurvesFill} $\Sigma_{1}$
will be a filling Dehn sphere of $S^{3}$. On the other hand,
Johansson's example of Figure \ref{fig29a} represents a liftable
(to an embedding in $\mathbb{R}^{4}$) filling Dehn sphere of
$S^{3}$ (see also \cite{Giller}). Thus in $S^{3}$ there are
liftable and non-liftable filling Dehn spheres. A generalization
for general 3-manifolds of the mentioned result of \cite{Giller}
could be applied for giving a complete answer to Fenn's question.

\subsection{The triple point spectrum.\label{SUBSECTION Triplepoint Spectrum}}

The minimal number of triple points of filling Dehn surfaces of a
3-manifold $M$ satisfying some particular property can be in some
cases a topological invariant of $M$. We define the \textit{triple
point number} $t(M)$ of a closed orientable 3-manifold $M$ as the
minimal number of triple points of all its filling Dehn surfaces
and the \textit{genus }$g$\textit{\ triple point number}
$t_{g}(M)$ of $M$ as the minimal number of triple points of all
its genus $g$ filling Dehn surfaces. The ordered collection
$\left( t_{0}(M),t_{1}(M),t_{2}(M),...\right)  $ of all the genus
$g$ triple point numbers of the 3-manifold $M$ for all $g\geq0$ is
what we call the \textit{triple point spectrum} $\frak{T}(M)$ of
$M$. We can make similar definitions imposing topological
restrictions on the filling Dehn surfaces considered. For example,
we can define the \textit{nulhomotopic triple point number} of $M$
as the minimal number of triple points of all its nulhomotopic
filling Dehn surfaces, and in a similar way the
\textit{nulhomotopic genus }$g$\textit{\ triple point number} or
the \textit{nulhomotopic triple point spectrum }can be defined.
All of them are topological invariants of the 3-manifold and give
a measure of the\textit{\ complexity} of the manifold in the same
way as the Heegaard genus, for example. If we have a filling Dehn
surface $\Sigma$ in a 3-manifold, using pipings as that of Figure
\ref{fig23b}, perhaps we can reduce the number of triple points of
$\Sigma$, but increasing the genus of the filling Dehn surface. So
there is some relation between the different genus $g$ triple
point numbers that would be interesting to clarify.

Any Dehn sphere in a closed orientable 3-manifold has an even
number of triple points (\cite{Haken1}, p. 105). This is not the
case for genus $g>0$ Dehn surfaces, as it can be seen in the
example given by Figure \ref{fig32a}. This means that if we want a
set of moves for relating \textit{all} Dehn surfaces (of any
genus) of any 3-manifold, the Homma-Nagase moves introduced here,
together with pipings, do not suffice because all of them are
operations that preserve the parity of the number of triple
points.

We define that a genus $g$ filling Dehn surface $\Sigma$ of a
3-manifold $M$ is \textit{minimal} if there is no other genus $g$
filling Dehn surface of $M $ with less triple points than
$\Sigma$. Minimal filling Dehn surfaces, in particular minimal
filling Dehn spheres, should have interesting properties, and
their classification is another interesting problem. The
classification of minimal Dehn spheres has been solved for $S^{3}$
in \cite{A.Shima2}. In that work, A. Shima gives in a different
context six examples of Dehn spheres in $S^{3}$ with only 2 triple
points. Three of these six examples fill $S^{3}$ (one of them is
Johansson's sphere of Figure \ref{fig29a}) and they are minimal
because, as we have said, any filling Dehn sphere must have at
least 2 triple points. It can be deduced by the main theorem of
\cite{A.Shima2} that these three examples are the unique possible
minimal filling Dehn spheres in $S^{3}$.

Finally, we want to introduce a later definition. We say that a filling Dehn
surface $\Sigma$ in a 3-manifold $M$ is \textit{irreducible} if the only
allowable filling preserving moves on $\Sigma$ are finger move +1 or +2. That
is, $\Sigma$ is irreducible if any Dehn surface $\Sigma^{\prime}$ which can be
obtained performing a filling-preserving move on $\Sigma$ has more triple
points than $\Sigma$. Johansson's sphere is not irreducible, while Example 1.3
of \cite{A.Shima2} is irreducible. This means that minimality does not imply
irreducibility. We are interested also in the converse question: are there
examples of non-minimal irreducible filling Dehn surfaces?

\end{document}